\documentclass[12pt]{amsart}%

\usepackage{amssymb}
\usepackage{amsmath}
\usepackage{amsthm}
\usepackage{amscd}
\usepackage{upgreek}
\usepackage[margin=1in]{geometry}
\usepackage{hyperref,xcolor,tikz-cd, bbold}

\usepackage{stackengine}

\theoremstyle{plain}
\newtheorem{theorem}{Theorem}[section]
\newtheorem{proposition}{Proposition}[section]
\newtheorem{corollary}{Corollary}[section]
\newtheorem{lemma}{Lemma}[section]

\theoremstyle{remark}

\newtheorem{remark}{Remark}[section]
\newtheorem{examples}{Example}[section]
\newtheorem{assumption}{Assumption}[section]

\DeclareMathOperator{\dist}{dist}
\DeclareMathOperator{\supp}{supp}
\DeclareMathOperator{\diam}{diam}

\DeclareMathOperator{\lin}{span}
\DeclareMathOperator{\cpct}{cap}
\DeclareMathOperator{\Dfcpct}{Cap}
\DeclareMathOperator{\diag}{diag}
\DeclareMathOperator{\lip}{Lip}

\DeclareMathOperator{\Alt}{Alt}
\DeclareMathOperator{\Sym}{Sym}
\DeclareMathOperator{\sgn}{sgn}

\DeclareMathOperator{\im}{\mathrm{im}}
\DeclareMathOperator{\inj}{\mathrm{inj}}

\begin{document}

\title{A tensor product approach to non-local differential complexes}
\author{Michael Hinz$^1$}
\address{$^1$ Fakult\"{a}t f\"{u}r Mathematik, Universit\"{a}t Bielefeld, Postfach 100131, 33501 Bielefeld,
Germany}
\email{mhinz@math.uni-bielefeld.de}
\author{J\"orn Kommer$^2$}
\address{$^2$ Fakult\"{a}t f\"{u}r Mathematik, Universit\"{a}t Bielefeld, Postfach 100131, 33501 Bielefeld,
Germany}
\thanks{$^1$, $^2$ Research supported by the DFG IRTG 2235: 'Searching for the regular in the irregular: Analysis of singular and random systems'. }
\email{jkommer@math.uni-bielefeld.de}

\date{\today}

\begin{abstract}
We study differential complexes of Kolmogorov-Alexander-Spanier type on metric measure spaces associated with unbounded non-local operators, such as operators of fractional Laplacian type. We define Hilbert complexes, observe invariance properties and obtain self-adjoint non-local analogues of Hodge Laplacians. For $d$-regular measures and operators of fractional Laplacian type we provide results on removable sets in terms of Hausdorff measures. We prove a Mayer-Vietoris principle and a Poincar\'e lemma and verify that in the compact Riemannian manifold case the deRham cohomology can be recovered. 
\tableofcontents
\end{abstract}

\maketitle

\section{Introduction}

In this article we study differential complexes of functions associated with unbounded non-local operators on metric measure spaces; particular examples are operators of fractional Laplacian type. 

The classical deRham cohomology theory describes how topological features of a smooth manifold are detected by cochain complexes of differential forms, \cite{BT82, dRh60, Warner}. There are numerous related cohomology theories, for instance those 
used for cell complexes and groups, \cite{Eckmann44, Lueck}, and those based on (commutative or non-commutative) algebras, \cite{CM90, GVF, Hochschild45}. Applications of related theories to data analysis, \cite{Carlsson, ChazaldeSilvaOudot, ELZ}, sparked new interest in scaled differential complexes and cohomologies on metric spaces, \cite{BSSS12, Genton, SS12}; ideas had been sketched earlier in \cite{Pansu96, Pansu04}. The complexes studied in these articles are of \emph{Kolmogorov-Alexander-Spanier type}, \cite{Alexander35, Kolmogorov36, Kolmogorov36a, Massey78}. They retrieve metric information and involve a parameter that determines the scale at which features are recognized. This is conceptually close to Vietoris-Rips complexes, \cite{Gromov87, Hausmann95, Latschev01, Vietoris27}, and linked to uniform structures, \cite[Chapter II]{Bourbaki}. Approaches to homology and homotopy involving a metric scale parameter can be found in \cite{BCW, Plaut18, PlautWilkins12}. 

If the given metric space is endowed with a suitable measure, one can define Kolmogorov-Alexander-Spanier complexes based on 
$L^2$-spaces of (classes of) functions; this was done in \cite{BSSS12, Genton, SS12}. Using the terminology of \cite{BL92} in a slightly wider sense, they may be seen as examples of \emph{Hilbert complexes}. This point of view is in line with the classical variational approach to Hodge theory, \cite{dRh60,dRh-K50, Gaffney55, Hodge41, Kodaira49}, and it emphasizes the link to elliptic partial differential equations. 

At order zero (that is, acting on scalar functions) the exterior derivative on a Riemannian manifold $M$, seen as a closed \emph{unbounded} operator between $L^2$-spaces, defines an unbounded local Dirichlet form, \cite{BH91, FOT94, Grigoryan2009}. Its infinitesimal generator is the self-adjoint Dirichlet Laplacian on $M$, clearly \emph{unbounded}. The $L^2$-complexes on metric measure spaces $X$ considered in \cite{BSSS12, Genton, SS12} involve \emph{bounded} coboundary operators, and at order zero one obtains a bounded (purely) non-local Dirichlet form. Its generator is a \emph{bounded} non-local operator. 

From the perspective of partial differential equations it is more natural to consider \emph{unbounded} non-local operators and their Dirichlet forms, and there is rich literature on such operators on Euclidean spaces, including fractional Laplacians. See for instance \cite{BucurValdinoci16, diNezza12, Garofalo19} for basic concepts and and applications, \cite{CaffarelliSilvestre07, DiCastroKuusiPalatucci14, DipierroRos-OtonValdinoci17, Grubb15, KassmannMimica17, Nazarov22, PalatucciSavinValdinoci13, Ros-OtonSerra12, Silvestre05} for a number of well-known results and \cite{AbatangeloValdinoci14, CaffarelliRoquejoffreSavin10, ChangGonzales11, DaLioRiviere11, Garofalo19, PalatucciKuusi18, SpenerWeberZacher} for connections to geometric analysis and recent developments. Moreover, there is a well-established theory of \emph{unbounded} purely non-local Dirichlet forms on metric measure spaces $(X,\varrho,\mu)$; prototype examples are quadratic forms 
\begin{equation}\label{E:fracLaplaceform}
\mathcal{E}(f)=\int_X\int_X (f(x)-f(y))^2\varrho(x,y)^{-d-\alpha}\mu(dx)\mu(dy)
\end{equation}
associated with non-local operators of fractional Laplace type, where $\mu$ is a $d$-regular measure on $X$ and $\alpha\in (0,2)$. See for instance \cite{BarlowGrigKumagai09, ChKu03, ChKu08, GrigHuHu18, GrigHuLau14, Stos00} and the references therein.

We propose an approach to higher order 'differential forms' based on \emph{unbounded} non-local operators on metric measure spaces and robust enough to include the complex induced by (\ref{E:fracLaplaceform}). To our knowledge, this has not been discussed anywhere else. For 'differential forms' of order one some aspects of a more general theory based on non-local Dirichlet forms have been investigated, see for instance \cite{CS03, H15,  HRT13, Sauvageot89, Sauvageot90}, but higher orders have not been studied.\\

We set ourselves two objectives:
\begin{enumerate}
\item[(a)] To define 'non-local' Hilbert complexes of Kolmogorov-Alexander-Spanier type, involving a metric scale and based on purely non-local Dirichlet forms with unbounded integral kernels.
\item[(b)] To prove that in the case of compact Riemannian manifolds the cohomologies of suitable non-local complexes are isomorphic to the deRham cohomologies.
\end{enumerate}

The passage from bounded kernels as in \cite{BSSS12} to unbounded kernels such as $j(x,y)=\varrho(x,y)^{-d-\alpha}$ in (\ref{E:fracLaplaceform}) must be carried out carefully. To meet objective (a), one has to introduce spaces of sufficiently regular multivariate functions that are compatible with applications of the coboundary operator. For (b) one needs partitions of unity, and this requires the boundedness of certain multipliers. 

Objective (a) gives differential complexes connected to unbounded non-local operators in a similar way as the deRham complex is connected to the Laplace-Beltrami operator. One motivation for (a) is a close connection to partial differential
equations. For instance, Poisson regularity was listed as an open problem in \cite{BSSS12}, and one may expect this problem to have a more familiar flavor in the context of unbounded non-local operators, cf. \cite{Silvestre05}. A second motivation for (a) is that non-local Dirichlet forms (\ref{E:fracLaplaceform}) and related operators can readily be defined on a wide class of metric measure spaces, including fractal spaces for which local Dirichlet forms or Laplacians are not known to exist (or definitely do not exist). In other words, even if no theory of local complexes is available, such non-local complexes and their cohomologies can be studied. A third motivation for (a) is the possibility to tune the sensitivity of the complex towards the removal of small closed sets. If the starting point is (\ref{E:fracLaplaceform}), then there is the 'new' (in comparison to \cite{BSSS12, Genton, SS12}) parameter $\alpha\in (0,2)$ that can be varied. We obtain Hilbert complexes with a metric scale parameter $\varepsilon$ and the parameter $\alpha$; finer details of $X$ are noted by the complex if $\varepsilon$ is decreased or $\alpha$ is increased. Removable sets have also been studied for local complexes, \cite[Theorem 4.4.]{BL92}, \cite[Theorem 5.1]{GT10}, but as for Sobolev spaces, \cite{AH96, Mazya85}, the possibility to parametrize the effect (while keeping the order of integrability fixed) is a feature of the fractional case. Apart from these motivations, objective (a) should also be seen as a key step towards a more general theory of 'differential forms' based on Dirichlet forms (local, non-local or mixed). In view of known results for first order forms, \cite{CS03, H15, HRT13, HTams, HinzTeplyaev18, Sauvageot89, Sauvageot90, W00} and related results for higher order forms under somewhat different hypotheses, \cite{Gigli17}, such a theory seems desirable. 

Objective (b) is set to indicate that our approach is not detached from classical theories.  The recovery of the deRham cohomologies on compact manifolds by non-local complexes may be seen as an analog of similar results proved in \cite{BSSS12} and \cite{Genton} for bounded kernels. A new obstacle is that a priori we have to distinguish between the cohomologies based on 'smooth functions or forms' and those based on closed operators (in the Hilbert complex sense). In the deRham case additional smoothing arguments are used to show that these cohomologies are isomorphic, \cite{BL92, Cheeger}. Such arguments seem rather out of bounds for the non-local complexes we consider. However, in the non-local case the Poincar\'e lemma involves a homotopy operator that combines well with operator closures, and this can be used to obtain a similar isomorphy. 

Our perspectives and results are as follows: We consider algebras of suitable bounded functions (eventually continuous, of bounded differentials and energy finite) and consider linear combinations of antisymmetric tensor products, which we call \emph{elementary $p$-functions}. The basic idea is standard, see for instance \cite{BlooreRoberts, CM90} or \cite[Section II.4]{RS80}. By a slight change of notation one can express the tensor products in a more 'simplicial' (or 'affine') manner, Proposition \ref{P:Vp}; for cochains of higher order the single difference in (\ref{E:fracLaplaceform}) is then replaced by a determinant of differences, Example \ref{Ex:derivations}. This perspective shows that if the individual scalar valued factors are sufficiently regular, the integrability of elementary $p$-functions with respect to unbounded kernels is guaranteed, cf. Lemma \ref{L:elemL2} and Example \ref{Ex:upperenergies}. The dominating geometric idea is that of simple vectors, but with 'infinitesimal displacements' replaced by differences (jumps) according to the kernel. For each given order $p$ we view the elementary $p$-functions as a core for the coboundary operators and related energy forms, a perspective in line with the classical Beurling-Deny theory, \cite{Allain75,BD58, BD59, FOT94}. To meet objective (a) we then assume that the kernels are absolutely continuous (as in (\ref{E:fracLaplaceform})) and pass to a Hilbert complex by taking closures of the coboundary operators, Theorem \ref{T:closable}. We address questions of invariance, Theorem \ref{T:invariance}, introduce coderivations, Corollary \ref{C:adjoint}, observe weak Hodge decompositions, Corollary \ref{C:weakHodge}, and obtain self-adjoint non-local analogs of Hodge Laplacians, Proposition \ref{P:HodgeLaplace}. We discuss removable sets, Theorem \ref{T:removable}, with results in terms of Hausdorff measures for the special case related to (\ref{E:fracLaplaceform}), Corollary \ref{C:catalogue}. Objective (b) we discuss for compact spaces, following the classical path, \cite{BT82, Weil52}. We prove a Mayer-Vietoris principle for finite open covers, Proposition \ref{P:MVlemma}. For suitable covers by open balls and under an assumption which is a variant of \cite[Hypothesis  ($\ast$), p. 34]{BSSS12} we establish a Poincar\'e lemma; see Assumption \ref{A:hypostar}, Proposition \ref{P:Poincare} and Corollary \ref{C:Poincare}. The arguments are sufficiently robust to work for both the cores and the domains of the closures, and they show that the respective cohomologies are isomorphic to the \v{C}ech cohomology of the cover, Corollary \ref{C:Poincare2}. For compact Riemannian manifolds specific finite covers by sufficiently small balls satisfy all assumptions. In this case both the cohomologies based on elementary functions and those of the Hilbert complex are isomorphic to the respective deRham cohomologies. 

We point out that under fairly general assumptions non-local complexes as studied here can be used to approximate local complexes in a natural way; this will be the subject of a follow up article, \cite{HinzKommer22+}.
 
 Section \ref{S:Prelim} contains preliminaries and basic notation. Sections \ref{S:elem} and \ref{S:kernels} discuss elementary $p$-functions and their integrability properties. Operator closures and Hilbert complexes are studied in Section \ref{S:Hilbert}. Open covers and the link to deRham cohomology in the compact manifold case are the subject of Section \ref{S:covers}, along with simple examples in Subsection \ref{SS:Examples}.

\section{Preliminaries}\label{S:Prelim}

Let $X$ be a nonempty set and $p\geq 1$ and integer. Let $\mathcal{S}_p$ denote the symmetric group of order $p$.  A function $F\colon X^p\to \mathbb{R}$ is called \emph{symmetric} if $F(x_{\sigma(1)}, \dotsc, x_{\sigma(p)}) = F(x_1, \dotsc, x_p)$ for all $\sigma\in \mathcal{S}_p$, $x_1, \dotsc, x_p\in X$, and \emph{antisymmetric} if $F(x_{\sigma(1)}, \dotsc, x_{\sigma(p)}) = (\sgn \sigma) F(x_1, \dotsc, x_p)$ for all $\sigma\in \mathcal{S}_p$, $x_1, \dotsc,x_p\in X$;  here $\sgn \sigma$ denotes the sign of the permutation $\sigma$.  The \emph{symmetrizer} $\Sym_p$ of order $p$ is defined by
\begin{equation}\label{E:def-symm}
	\Sym_p(F)(x_1, \dotsc, x_p): = \frac{1}{p!}\sum_{\sigma\in \mathcal{S}_p} F(x_{\sigma(1)}, \dotsc, x_{\sigma(p)});
\end{equation}
it takes a given function $F: X^p\to\mathbb{R}$ into a symmetric function $\Sym_p(F): X^p\to\mathbb{R}$. Clearly  $\Sym_1=\mathrm{id}$, $\Sym_p(F) = F$ if $F$ is symmetric, and $\Sym_p^2 = \Sym_p$. The \emph{pointwise product} $GF$ of a symmetric function $G$ on $X^p$ and a function $F$ on $X^p$ satisfies
\begin{equation}\label{E:prod-sym-sym}
	\Sym_p(GF) = G\Sym_p(F);
\end{equation}
if in addition $F$ is symmetric, then $FG$ is symmetric and equals (\ref{E:prod-sym-sym}). The \emph{antisymmetrizer} $\Alt_p$ of order $p$ is defined by
\begin{equation}\label{E:def-antisymm}
	\Alt_p(F)(x_1, \dotsc, x_p): = \frac{1}{p!}\sum_{\sigma\in \mathcal{S}_p} (\sgn \sigma) F(x_{\sigma(1)}, \dotsc, x_{\sigma(p)});
\end{equation}
it takes a given function $F: X^p\to\mathbb{R}$ into an antisymmetric function $\Alt_p(F): X^p\to\mathbb{R}$. Similarly as before, $\Alt_1=\mathrm{id}$, $\Alt_p(F) = F$ if $F$ is antisymmetric, and $\Alt_p^2 = \Alt_p$. The pointwise product $GF$ of a symmetric function $G$ on $X^p$ and a function $F$ on $X^p$ satisfies
\begin{equation}\label{E:prod-sym-alt}
	\Alt_p(GF) = G\Alt_p(F);
\end{equation}
if in addition $F$ is antisymmetric, then $FG$ is again antisymmetric and equals (\ref{E:prod-sym-alt}).

Given a function $F: X^p\to \mathbb{R}$,  we consider the function $\delta_{p-1} F: X^{p+1}\to\mathbb{R}$, defined by 
\begin{equation}\label{E:def-coboundary}
	\delta_{p-1} F(x_0, \dotsc, x_p):= \sum_{i=0}^p (-1)^i F(x_0, \dotsc, \hat{x}_i, \dotsc, x_p),
\end{equation}
where, as usual, $\hat{x}_i$ means that $x_i$ is omitted. Obviously $\delta_0\mathbf{1}=0$, and it is easy to check that
\begin{equation}\label{E:coboundary-square-zero}
\delta_p\circ\delta_{p-1}=0, \quad p\geq 1.
\end{equation}
The operators $\delta_p$ are the \emph{(Kolmogorov-Alexander-Spanier) coboundary operators}, \cite{Alexander35, Massey78, Spanier48, Spanier66}.
For functions of the form $F = f_1\otimes \dotsb\otimes f_p$ with $f_i\colon X\to\mathbb{R}$ (cf. \cite{BlooreRoberts, CM90}), we obtain
\begin{equation}\label{E:coboundary-tensor-prod}
	\delta_{p-1}(f_1\otimes\dotsb\otimes f_p) = \mathbf{1}\otimes f_1\otimes\dotsb\otimes f_p + \sum_{i=1}^p (-1)^i f_1\otimes\dotsb\otimes f_i\otimes \mathbf{1}\otimes f_{i+1}\otimes \dots\otimes f_p.
\end{equation} 
It is easily seen that for any integer $p\geq 1$ we have 
\begin{equation}\label{E:alt-coboundary-commute}
	\delta_{p-1}\circ \Alt_p = \Alt_{p+1}\circ\: \delta_{p-1}.
\end{equation}  
Identities (\ref{E:alt-coboundary-commute}) and (\ref{E:coboundary-tensor-prod}) together imply that 
\begin{equation}\label{E:coboundary-of-alt-is-alt}
	\delta_{p-1}\Alt_p(f_1\otimes\dotsb\otimes f_p) = (p+1)\Alt_{p+1}(\mathbf{1}\otimes f_1\otimes\dotsb\otimes f_p).
\end{equation}
This can be rewritten as
\begin{align}\label{E:determinant}
\delta_{p-1}\Alt_p(f_1\otimes...\otimes f_p)(x_0,x_1,...,x_p)&=\frac{1}{p!}\det[(f_i(x_j)-f_i(x_0))_{i,j=1}^p]\notag\\
&=\Alt_p(\delta_0 f_1(x_0,\cdot)\otimes...\otimes \delta_0f_p(x_0,\cdot))(x_1,...,x_p).
\end{align}

\section{Complexes of elementary functions}\label{S:elem}

Let $(X,\varrho)$ be a metric space and $\mathcal{C}$ an algebra of real valued functions on $X$. We consider complexes of antisymmetric tensor products of elements of $\mathcal{C}$. 

Let $\mathcal{C}^{a,0}:=\mathcal{C}$, and for $p\geq 1$ define
\begin{equation}\label{E:Cp}
\mathcal{C}^{a,p}:=\lin \{ \Alt_{p+1}(f_0\otimes f_1\otimes \dots\otimes f_p): f_0\in \mathcal{C}\oplus \mathbb{R}\ \text{and}\ f_1,\dots, f_p \in \mathcal{C}\}.
\end{equation}
To the elements of $\mathcal{C}^{a,p}$ we refer as \emph{elementary $p$-functions}. To save notation we will write $\mathcal{C}^p:=\mathcal{C}^{a,p}$, except where we wish to point out the antisymmetry explicitely.

We call a family $N_\ast=(N_p)_{p\geq 0}$ a \emph{system of diagonal neighborhoods} for $X$ if 
\begin{enumerate}
\item[(i)] the $N_p$, $p\geq 0$, are neighborhoods of the diagonal $\diag_p:=\{(x_0,...,x_0): x_0\in X\}$ in $X^{p+1}$, respectively, and either all open or all closed,
\item[(ii)] the $N_p$ are symmetric in the sense that for any $\pi\in \mathcal{S}_{p+1}$ and any $(x_{0},...,x_p)\in N_p$
we have $(x_{\pi(0)},...,x_{\pi(p)})\in N_p$, 
\item[(iii)] for any $p\geq 1$, any $(x_0,...,x_p)\in N_p$ and any $0\leq i\leq p$ we have $(x_0,...,\hat{x}_i,...,x_p)\in N_{p-1}$.
\end{enumerate}
Note that $N_0=X$. The choice of a system of diagonal neighborhoods determines a metric scale. Given systems 
$N_\ast=(N_p)_{p\geq 0}$ and $N'_\ast=(N'_p)_{p\geq 0}$ we write $N_\ast\prec N'_\ast$ if $N_p\subset N'_p$ for all $p\geq 0$.
(We use the symbol $\subset$ in the non-strict sense; equality is permitted.)

\begin{examples}\label{Ex:N}\mbox{}
\begin{enumerate}
\item[(i)] Setting $N_p:=X^{p+1}$, $p\geq 0$, we obtain a system of diagonal neighborhoods; it is the largest possible in the sense of $\prec$.
\item[(ii)] Given $\varepsilon>0$, the sets 
\begin{equation}\label{E:boundedrange}
N_p(\varepsilon):=\{(x_0,...,x_p)\in X^{p+1}: \max_{0\leq i<j\leq p}\varrho(x_i,x_j)<\varepsilon\}
\end{equation}
form a system $N_\ast(\varepsilon)$ of diagonal neighborhoods for $X$. The strict inequality in (\ref{E:boundedrange}) could be replaced by $\leq$. 
\item[(iii)] For $p\geq 0$ we write $\varrho_p((x_0,..,x_p),(y_0,...,y_p)):=\max_{0\leq i\leq p}\varrho(x_i,y_i)$ 
and denote the resulting distance between points and sets by $\dist_p$. Given $\varepsilon>0$ also the sets 
\begin{equation}\label{E:boundedrangehat}
\hat{N}_p(\varepsilon):=\{(x_0,...,x_p)\in X^{p+1}: \dist_p((x_0,...,x_p),\diag_p)\leq \varepsilon\}
\end{equation}
form a system $\hat{N}_\ast(\varepsilon)$ of diagonal neighborhoods for $X$. These sets were used in \cite[Section 7]{BSSS12}.
It is easily seen that $N_\ast(\varepsilon)\prec \hat{N}_\ast(\varepsilon)\prec N_p(3\varepsilon)$, cf. \cite[p. 17]{Genton}. 
\item[(iii)] Given an open cover $\mathcal{V}=\{V_\alpha\}_{\alpha\in I}$ of $X$, the sets 
\begin{equation}\label{E:comesfromcover}
N_p(\mathcal{V}):=\bigcup_{\alpha\in I} V_\alpha^{p+1}
\end{equation}
form a system $N_\ast(\mathcal{V})$ of diagonal neighborhoods, cf. \cite[p. 346]{CM90}.
\end{enumerate}
\end{examples}

Given a system of diagonal neighborhoods $N_\ast=(N_p)_{p\geq 0}$ and an integer $p\geq 0$, we use the notation
\begin{equation}\label{E:defasrestriction}
\mathcal{C}^p(N_p):=\mathcal{C}^p|_{N_p}\quad \text{(respectively $\mathcal{C}^{a,p}(N_p):=\mathcal{C}^{a,p}|_{N_p}$).}
\end{equation} 
Clearly $\mathcal{C}^p(X^{p+1})=\mathcal{C}^p$.

By (\ref{E:coboundary-square-zero}) and (\ref{E:coboundary-of-alt-is-alt}) the spaces $\mathcal{C}^p(N_p)$ are seen to form differential complexes; they are of \emph{Kolmogorov-Alexander-Spanier type}.

\begin{proposition}\label{P:elemcomplex} 
Let $N_\ast$ be a system of diagonal neighborhoods. Then the sequence
\begin{equation}\label{E:elemcomplex}
0\longrightarrow \mathcal{C}^0(N_0)\stackrel{\delta_0}{\longrightarrow}\mathcal{C}^1(N_1)\stackrel{\delta_1}{\longrightarrow} ... \stackrel{\delta_{p-1}}{\longrightarrow} \mathcal{C}^p(N_p)  \stackrel{\delta_p}{\longrightarrow} ...
\end{equation}
is a cochain complex.
\end{proposition}
We write $(\mathcal{C}^\ast(N_\ast), \delta_\ast)$ to denote this complex, and we refer to it as the \emph{elementary complex}. Here and in the following we agree to set $\delta_{-1}:=0$. For any integer $p\geq 0$ the $p$-th cohomology of (\ref{E:elemcomplex}) is 
\begin{equation}\label{E:elemcoho}
H^p\mathcal{C}^\ast(N_\ast):=\ker \delta_p|_{\mathcal{C}^p(N_p)} / \im \delta_{p-1}|_{\mathcal{C}^{p-1}(N_{p-1})}.
\end{equation}

\begin{remark} It is well-known that whether the cohomologies for $p\geq 1$ can be nontrivial or not depends on the choice of $N_\ast$, see Subsection \ref{SS:Examples} or \cite{BSSS12, Genton}. 
\end{remark}

We rewrite the spaces $\mathcal{C}^p(N_p)$ in a convenient way.

\begin{proposition}\label{P:Vp}
For any integer $p\geq 1$ we have 
\begin{equation}\label{E:Vp}
\mathcal{C}^{p}=\lin \{ \overline{g}\:\delta_{p-1}\Alt_p(f_1\otimes \dots\otimes f_p): g\in \mathcal{C}\oplus \mathbb{R} \ \text{and}\ f_1,\dots, f_p \in \mathcal{C}\},
\end{equation}
where $\overline{g}$ denotes the average 
\begin{equation}\label{E:symmetrization}
\overline{g}(x_0,x_1,...,x_p):=\frac{1}{p+1}\sum_{i=0}^p g(x_i),\quad g\in \mathcal{C}\oplus \mathbb{R}.
\end{equation}
Moreover, for any $g\in \mathcal{C}\oplus \mathbb{R}$ and $f_1, ..., f_p \in\mathcal{C}$ the identity
\begin{equation}\label{E:differentiate}
\delta_p(\bar{g}\delta_{p-1}\Alt_p(f_1\otimes \dots\otimes f_p)) = \delta_p\Alt_{p+1}(g\otimes f_1\otimes \dots\otimes f_p)
\end{equation}
holds.
\end{proposition}

\begin{proof}
Identity (\ref{E:Vp}) is immediate from the fact that for any $g\in \mathcal{C}\oplus \mathbb{R}$ and $f_1,\dots, f_p, \in \mathcal{C}$ we have, by (\ref{E:coboundary-of-alt-is-alt}) and (\ref{E:prod-sym-alt}), 
\begin{align}
\overline{g} \delta_{p-1}\Alt_p(f_1\otimes \cdots \otimes f_p)&=(p+1)\overline{g}\Alt_{p+1}(\mathbf{1}\otimes f_1\otimes\cdots\otimes f_p)\notag \\ 
&=\sum_{k=0}^p\Alt_{p+1}(\mathbf{1}\otimes f_1\otimes\cdots\otimes (gf_k)\otimes\dotsb\otimes f_p)\notag \\
&= \Alt_{p+1}(g\otimes f_1\otimes\cdots \otimes f_p)\notag\\
&\hspace{50pt} + \sum_{k=1}^p\Alt_{p+1}(\mathbf{1}\otimes f_1\otimes\dotsb\otimes (gf_k)\otimes\dotsb\otimes f_p)\notag \\
&= \Alt_{p+1}(g\otimes f_1\otimes\cdots \otimes f_p)\notag\\
&\hspace{50pt}  + \frac{1}{p+1}\sum_{k=1}^p\delta_{p-1}\Alt_p(f_1\otimes\dotsb\otimes (gf_k)\otimes\dotsb\otimes f_p).\notag
	\end{align}
By (\ref{E:coboundary-square-zero}) an application of $\delta_p$ yields (\ref{E:differentiate}).
\end{proof}

\begin{examples}\label{Ex:derivations}\mbox{}
\begin{enumerate}
\item[(i)] For $p=0$ we have 
\[\delta_0 f(x_0,x_1)=f(x_1)-f(x_0).\] 
\item[(ii)] For $p=1$ we have 
\[\delta_1 (\overline{f}_1\delta_0 f_2)(x_0,x_1,x_2)=\delta_1\Alt_2(f_1\otimes f_2)(x_0,x_1,x_2)=\frac12\det[((f_i(x_j)-f_i(x_0))_{i,j=1}^2].\] 
\end{enumerate}
\end{examples}

\begin{remark} The averaging (\ref{E:symmetrization}) of the coefficients amounts to an antisymmetrization of the product: 
Given $g\in \mathcal{C}\oplus \mathbb{R}$ and $F:X^{p+1}\to\mathbb{R}$, we define a function $g\cup F:X^{p+1}\to\mathbb{R}$ by $ 
g\cup F(x_0,...,x_p):=g(x_0)F(x_0,...,x_p)$. If $F$ is antisymmetric, then 
\begin{equation}\label{E:averageandantisymm}
\Alt_{p+1}(g\cup F)=\overline{g}F.
\end{equation}
\end{remark}

From (\ref{E:def-symm}) and (\ref{E:def-antisymm}) it is easily seen that for any integer $p\geq 1$, any $f_0,...,f_p,g_0,...,g_p\in \mathcal{C}\oplus \mathbb{R}$ we have 
\begin{multline}\label{E:symmult}
\Sym_{p+1}(g_0\otimes \cdots \otimes g_p)\Sym_{p+1}(f_0\otimes\cdots \otimes f_p)\\ = \frac{1}{(p+1)!} \sum_{\pi\in\mathcal{S}_{p+1}} \Sym_{p+1}(g_{\pi(0)}f_0\otimes\cdots \otimes g_{\pi(p)}f_p) 
\end{multline}
and
\begin{multline}\label{E:symmaltmult}
 \Sym_{p+1}(g_0 \otimes\cdots \otimes g_p)\Alt_{p+1}(f_0\otimes\cdots \otimes f_p)\\ = \frac{1}{(p+1)!}\sum_{\pi\in\mathcal{S}_{p+1}}\Alt_{p+1}( g_{\pi(0)}f_0\otimes\cdots \otimes g_{\pi(p)}f_p). 
 \end{multline}
Let $\mathcal{C}^{s,p}$ be the vector space spanned by $\{ \Sym_{p+1}(f_0\otimes f_1\otimes \dots\otimes f_p): f_0, f_1,\dots, f_p \in \mathcal{C}\oplus \mathbb{R}\}$ and $\mathcal{C}^{s,p}(N_p):=\mathcal{C}^{s,p}|_{N_p}$. From (\ref{E:symmult}) it follows that $\mathcal{C}^{s,p}(N_p)$ with pointwise multiplication is a subring of the ring of real valued functions on $N_p$. Using (\ref{E:symmaltmult}) it is seen that for any $F\in \mathcal{C}^{a,p}(N_p)$ and $\chi\in \mathcal{C}^{s,p}(N_p)$ we have $\chi F\in \mathcal{C}^{a,p}(N_p)$, and that this multiplicative action of $\mathcal{C}^{s,p}(N_p)$ makes $\mathcal{C}^{a,p}(N_p)$ a module.

For later use we record an observation about differentials of functions of a special form. Given $\chi\in \mathcal{C}\oplus \mathbb{R}$ and integer $p\geq 0$ we write
\begin{equation}\label{E:tensorpower}
\chi^{\otimes (p+1)}:=\chi\otimes \dots\otimes \chi,
\end{equation}
where $\chi$ appears $p+1$ times on the right hand side. 
\begin{lemma}\label{L:gradientsplitoff}
Let $p\geq 0$ be an integer and $f_0,...,f_p,\chi\in \mathcal{C}\oplus \mathbb{R}$. Then 
\begin{equation}\label{E:basicmultact}
\chi^{\otimes (p+1)}\Alt_{p+1}(f_0\otimes\dots\otimes f_p)=\Alt_{p+1}(\chi f_0\otimes\dots \otimes \chi f_p)
\end{equation}
and
\begin{multline}
\delta_{p}(\chi^{\otimes (p+1)}\Alt_{p+1}(f_0\otimes\dots\otimes f_p))(x_0, ...,x_{p+1})\notag\\
=\chi(x_1)\cdots \chi(x_{p+1})\:\delta_{p}\Alt_{p+1}(f_0\otimes\cdots\otimes f_p)(x_0,...,x_{p+1})\notag\\
+\sum_{k=1}^{p+1} (-1)^{k-1}\chi(x_1)\cdots\chi(x_{k-1})\delta_0\chi(x_0,x_k)\chi(x_{k+1})\cdots \chi(x_{p+1})\times\notag\\
\times\Alt_{p+1}(f_0\otimes\cdots\otimes f_p)(x_0,...,\hat{x}_k,...,x_{p+1}).
\end{multline}
\end{lemma}

\begin{proof} 
The first statement is obvious. To see the second, note that the evaluation of
\begin{align}
\delta_p\left(\chi^{\otimes (p+1)}\Alt_{p+1}(f_0\otimes\dotsm\otimes f_p)\right)&=(\mathbf{1}\otimes \chi\otimes\dotsm\otimes\chi)\delta_p\Alt_{p+1}(f_0\otimes\dotsm\otimes f_p)\notag\\
&\hspace{30pt} +\delta_p\left(\chi^{\otimes (p+1)}\Alt_{p+1}(f_0\otimes\dotsm\otimes f_p)\right)\notag\\
&\hspace{30pt}-(\mathbf{1}\otimes \chi\otimes\dotsm\otimes\chi)\delta_p\Alt_{p+1}(f_0\otimes\dotsm\otimes f_p)\notag
\end{align}
at $(x_0,...,x_{p+1})\in X^{p+2}$ gives
\begin{align}
\chi(x_1)&\dotsm\chi(x_{p+1})\delta_p\Alt_{p+1}(f_0\otimes\dotsm\otimes f_p)(x_0, \dots, x_{p+1})\notag\\
&-\sum_{k=1}^{p+1} (-1)^{k-1} \chi(x_0)\cdots\widehat{\chi(x_k)}\cdots\chi(x_{p+1})\Alt_{p+1}(f_0\otimes\dotsm\otimes f_p)(x_0, \dotsc, \hat{x}_k, \dotsc, x_{p+1}) \notag\\
&+\sum_{k=1}^{p+1} (-1)^{k-1} \chi(x_1)\dotsm\chi(x_{p+1})\Alt_{p+1}(f_0\otimes\dotsm\otimes f_p)(x_0, \dotsc, \hat{x}_k, \dotsc, x_{p+1}). \notag
\end{align}
\end{proof}

Suppose that $\tilde{X}$ is another metric space and $\varphi:\tilde{X}\to X$ is a map. Given $p\geq 0$ integer and a function $F:X^{p+1}\to \mathbb{R}$ we define a function $\varphi^\ast F:\tilde{X}^{p+1}\to \mathbb{R}$ by 
\[\varphi^\ast F(\tilde{x}_0,...,\tilde{x}_p):=F(\varphi(\tilde{x}_0),...,\varphi(\tilde{x}_p)),\quad (\tilde{x}_0,...,\tilde{x}_p)\in \tilde{X}^{p+1}.\]
We write $\varphi^\ast\mathcal{C}:=\{\varphi^\ast f: f\in \mathcal{C}\}$, and given a system $N_\ast=(N_p)_{p\geq 0}$ of diagonal neighborhoods on $X$, we set $\varphi^\ast N_p:=\{(\tilde{x}_0,...,\tilde{x}_p)\in \tilde{X}^{p+1}: (\varphi(\tilde{x}_0),...,\varphi(\tilde{x}_p))\in N_p\}$. The following is easily seen.

\begin{proposition}\label{P:pullback}
Suppose that $\tilde{\mathcal{C}}$ is an algebra of real valued functions on $\tilde{X}$ and $\varphi^\ast \mathcal{C}\subset \tilde{\mathcal{C}}$.
\begin{enumerate}
\item[(i)] The map $\varphi^\ast:\mathcal{C}\to \tilde{\mathcal{C}}$ is an algebra homomorphism and $\varphi^\ast \mathcal{C}$ is a subalgebra of $\tilde{\mathcal{C}}$. For each integer $p\geq 0$ we have $(\varphi^\ast\mathcal{C})^{a,p}=\varphi^\ast\mathcal{C}^{a,p}\subset \tilde{\mathcal{C}}^{a,p}$ and $(\varphi^\ast\mathcal{C})^{s,p}=\varphi^\ast\mathcal{C}^{s,p}\subset \tilde{\mathcal{C}}^{s,p}$, and $\varphi^\ast$ is a module homomorphism. If $\varphi$ is bijective and $\varphi^\ast \mathcal{C}= \tilde{\mathcal{C}}$, then the maps $\varphi^\ast$ are (algebra, module) isomorphisms.
\item[(ii)] If $\varphi$ is a homeomorphism from $\tilde{X}\to X$, then $\varphi^\ast N_\ast:=(\varphi^\ast N_p)_{p\geq 0}$ is a system of diagonal neighborhoods on $\tilde{X}$ and $\varphi^\ast$ is a cochain map from the complex $(\mathcal{C}^\ast(N_\ast),\delta_\ast)$ to the complex $(\tilde{\mathcal{C}}^\ast(\varphi^\ast N_\ast),\delta_\ast)$. In the case that $\varphi^\ast \mathcal{C}= \tilde{\mathcal{C}}$ the map $\varphi^\ast$ is an isomorphism of cochain complexes.
\end{enumerate}
\end{proposition}

\section{Kernels and measures}\label{S:kernels}

We introduce kernels and measures and related conditions on the algebra $\mathcal{C}$. Let $(X,\varrho)$ be a locally compact metric space. Let $\mathcal{B}(X)$ denote the Borel $\sigma$-algebra on $X$.

\begin{assumption}\label{A:C}\mbox{}
\begin{enumerate}
\item[(i)] We assume that  $j:X\times \mathcal{B}(X)\to [0,+\infty]$ is a kernel in the sense that for any $x\in X$ the map $A\mapsto j(x,A)$ is a Borel measure on $X$, locally finite on $X\setminus \{x\}$, and for any $A\in \mathcal{B}(X)$ the function $x\mapsto j(x,A)$ is Borel measurable. 
\item[(ii)] We assume that $\mathcal{C}$ is an algebra of bounded real valued Borel functions on $X$ such that
\begin{equation}\label{E:Gammabdj}
\sup_{x\in X}\int_{X}(f(x)-f(y))^2 j(x,dy)<+\infty,\quad f\in \mathcal{C}.
\end{equation}
\end{enumerate}
\end{assumption}

If a volume measure $\mu$ is given and satisfies the following assumption, then pointwise statements can be complemented by integrated versions.

\begin{assumption}\label{A:CE}\mbox{}
We assume that $\mu$ is a nonnegative Radon measure on $X$ with full support (that is, $\supp\mu=X$), $\mathcal{C}\subset L^2(X,\mu)$, and 
\begin{equation}\label{E:energyfinitej}
\int_X\int_X (f(x)-f(y))^2j(x,dy)\mu(dx)<+\infty,\quad f\in \mathcal{C}.
\end{equation}
\end{assumption}

\begin{remark} If in the presence of a measure $\mu$ as in Assumption \ref{A:CE} the supremum in (\ref{E:Gammabdj}) is replaced by an essential supremum, pointwise statements below remain true in the $\mu$-a.e. sense.
\end{remark}

Given a system $N_\ast=(N_p)_{p\geq 0}$ of diagonal neighborhoods and $x_0\in X$, we write
\[N_{p,x_0}:=\{(x_1,...,x_p)\in X^p:\ (x_0,x_1,...,x_p)\in N_p\},\quad p\geq 1.\]
We set $D_0:=\emptyset$,
\[D_{p}:=\left\lbrace (x_0,...,x_p) \in X^{p+1}: \ x_i=x_j\ \text{for some distinct $i,j\in \{0,...,p\}$}\right\rbrace,\quad p\geq 1,\]
and write 
\[D_{p,x_0}:=\{(x_1,...,x_p)\in X^p:\ (x_0,x_1,...,x_p)\in D_p\},\quad p\geq 1.\]

Suppose that Assumption \ref{A:C} is satisfied. We define kernels
\begin{equation}\label{E:kernelsjx0}
j_p(x_0,d(x_1,...,x_p)):=j(x_0,dx_1)\cdots j(x_0,dx_p),\quad x_0\in X,\quad p\geq 1,
\end{equation}
from $X$ to the $p$-fold product of $\mathcal{B}(X)$. Clearly $j_1=j$. For each fixed $x_0$ and $p$ the measures $j_p(x_0,\cdot)$ are Radon on $(X\setminus \{x_0\})^p$ and symmetric in the variables $x_1,...,x_p$. Below we consider the Hilbert spaces 
\[L^2(N_{p,x_0}\setminus D_{p,x_0},j_p(x_0,\cdot))\quad \text{with their natural norms}\quad  \left\|\cdot\right\|_{L^2(N_{p,x_0}\setminus D_{p,x_0},j_p(x_0,\cdot))}.\]

If in addition Assumption \ref{A:CE} is satisfied, we write $J_0:=\mu$ and define measures
\begin{equation}\label{E:Jp}
J_p(d(x_0,...,x_p)):=\frac{1}{p+1}\sum_{k=0}^p j_p(x_k,d(x_0,...,\hat{x}_k,...,x_p))\mu(dx_k),\quad p\geq 1.
\end{equation}
The measures $J_p$ are Radon on $X^{p+1}\setminus D_p$, respectively, and symmetric in the variables $x_0,...,x_p$. Below we consider the Hilbert spaces 
\[L^2(N_p\setminus D_p,J_p)\quad \text{with their natural norms}\quad \left\|\cdot\right\|_{L^2(N_p\setminus D_p,J_p)}.\]

\begin{remark}
Rewriting (\ref{E:Jp}) in terms of (\ref{E:kernelsjx0}), we observe that 
\[J_p(d(x_0,...,x_p)):=\frac{1}{p+1}\sum_{k=0}^p \prod_{\ell\neq k} j(x_k,dx_\ell)\mu(dx_k),\quad p\geq 1.\] 
A somewhat similar point of view has been pursued in \cite[Section 2]{BSSS12}, but with (bounded) kernels implemented into the definitions of operators, not spaces.
\end{remark}

\begin{remark}\label{R:trivialinclusions}
If $N_\ast$ and $N_\ast'$ are systems of diagonal neighborhoods and $N_\ast'\prec N_\ast$, then $\|\cdot\|_{L^2(N_p'\setminus D_p,J_p)}\leq \|\cdot\|_{L^2(N_p\setminus D_p,J_p)}$ and in particular, $L^2(N_p\setminus D_p,J_p)\subset L^2(N_p'\setminus D_p,J_p)$. Similarly for the spaces involving $x_0$.
\end{remark}

The maps $\Alt_{p}$ and $\Alt_{p+1}$ are linear and bounded on the spaces  $L^2(N_{p,x_0}\setminus D_{p,x_0},j_p(x_0,\cdot))$ and $L^2(N_p\setminus D_p,J_p)$, respectively. Since by the mentioned symmetries they are self-adjoint, they act as orthogonal projections. We write 
\begin{equation}\label{E:L2ax0}
L^2_a(N_{p,x_0}\setminus D_{p,x_0},j_p(x_0,\cdot)):=\Alt_{p}\Big(L^2(N_{p,x_0}\setminus D_{p,x_0},j_p(x_0,\cdot))\Big)
\end{equation}
and
\begin{equation}\label{E:L2a}
L^2_a(N_p\setminus D_p,J_p):=\Alt_{p+1}\Big(L^2(N_p\setminus D_p,J_p)\Big)
\end{equation}
for the closed subspaces obtained as their images.

Given $x_0\in X$ we write $\mathcal{C}^p(N_p,x_0)$ for the space of functions $(x_1,...,x_p)\mapsto F(x_0,x_1, ...,x_p)$ on $N_{p,x_0}$ with $F\in \mathcal{C}^p(N_p)$. To keep notation short we will denote the elements of $\mathcal{C}^p(N_p,x_0)$ also by $F$ instead of the correct $F(x_0,\cdot)$.

\begin{lemma}\label{L:elemL2} Suppose that Assumption \ref{A:C} is satisfied and
let $N_\ast$ be a system of diagonal neighborhoods. 
\begin{enumerate}
\item[(i)] For any $p\geq 1$ and $x_0\in X$ we have $\mathcal{C}^p(N_p,x_0) \subset L^2(N_{p,x_0}\setminus D_{p,x_0},j_p(x_0,\cdot))$. For any $g\in\mathcal{C}\oplus\mathbb{R}$ and $f_1,...,f_p\in\mathcal{C}$ the inequality
\begin{multline}\label{E:fixedbasepoint}
\left\|\overline{g}\delta_{p-1}\Alt_p(f_1\otimes \dots \otimes f_p)\right\|_{L^2(N_{p,x_0}\setminus D_{p,x_0},j_p(x_0,\cdot))}\\
\leq \left\|g\right\|_{\sup}\prod_{i=1}^p\left(\int_{N_{1,x_0}\setminus D_{1,x_0}} (f_i(y)-f_i(x_0))^2j(x_0,dy)\right)^{1/2}
\end{multline} 
holds. The functional $\left\|\cdot\right\|_{L^2(N_{p,x_0}\setminus D_{p,x_0},j_p(x_0,\cdot))}$ is a Hilbert seminorm on $\mathcal{C}^p(N_p,x_0)$; it is a norm if $\mathcal{C}\subset C(X)$ and $j(x_0,\cdot)$ has full support on $N_{1,x_0}$.
\item[(ii)] Suppose that also Assumption \ref{A:CE} is satisfied. Then for any integer $p\geq 1$ we have 
$\mathcal{C}^p(N_p) \subset L_a^2(N_p\setminus D_p,J_p)$. For any $g\in\mathcal{C}\oplus\mathbb{R}$ and $f_1,...,f_p\in\mathcal{C}$ and any $i=1,...,p$ the inequality
\begin{multline}\label{E:integrated}
\left\|\overline{g}\delta_{p-1}\Alt_p(f_1\otimes \dots \otimes f_p)\right\|_{L^2(N_p\setminus D_p,J_p)}\\
\leq \left\|g\right\|_{\sup}\prod_{k\neq i}\left(\sup_{x\in X}\int_{N_{1,x}\setminus D_{1,x}} (f_k(y)-f_k(x))^2j(x,dy)\right)^{1/2}\times\\
\times\left(\int_X\int_{N_{1,x}\setminus D_{1,x}} (f_i(y)-f_i(x))^2j(x,dy)\mu(dx)\right)^{1/2}
\end{multline}
holds. The functional $\left\|\cdot\right\|_{L^2(N_p\setminus D_p,J_p)}$ is a Hilbert seminorm on $\mathcal{C}^p(N_p)$; it is a norm if $\mathcal{C}\subset C(X)$ and for $\mu$-a.e. $x_0\in X$ the measure $j(x_0, \cdot)$ has full support on $N_{1,x_0}$.
\end{enumerate}
\end{lemma}

We write $B(x,r)$ to denote the open ball of radius $r>0$ centered at $x\in X$.

\begin{proof}
Using (\ref{E:determinant}),
\begin{align}
&\int_{N_{p,x_0}\setminus D_{p,x_0}} (\overline{g}\:\delta_{p-1}(\Alt_p(f_1\otimes \dots\otimes f_p))(x_0,...,x_p))^2\:j_p(x_0,d(x_1,...,x_p))\notag\\
&\leq \frac{\left\|g\right\|_{\sup}^2}{(p!)^2}\int_{N_{p,x_0}\setminus D_{p,x_0}} \left(\det((f_i(x_k)-f_i(x_0))_{i,k=1}^p)\right)^2 j_p(x_0,d(x_1,...,x_p))\notag\\
&\leq \frac{\left\|g\right\|_{\sup}^2}{(p!)^2}\sum_{\pi\in \mathcal{S}_p}\sum_{\sigma\in \mathcal{S}_p}\int_{N_{p,x_0}\setminus D_{p,x_0}} \prod_{i=1}^p |f_i(x_{\pi(i)})-f_i(x_0)|\times\notag\\
&\hspace{170pt}\times \prod_{k=1}^p |f_k(x_{\sigma(k)})-f_k(x_0)| j_p(x_0,d(x_1,...,x_p)).\notag
\end{align}
By Cauchy-Schwarz with respect to the measure $j_p(x_0,\cdot)$, its symmetry in $x_1,...,x_p$ and Fubini this is seen to be bounded by 
\begin{align}
\left\|g\right\|_{\sup}^2 & \int_{N_{p,x_0}\setminus D_{p,x_0}}  \prod_{i=1}^p  (f_i(x_i)-f_i(x_0))^2 j_p(x_0,d(x_1,...,x_p))\notag\\
& \leq \left\|g\right\|_{\sup}^2 \prod_{i=1}^p \left(\int_{N_{1,x_0}\setminus D_{1,x_0}} (f_i(x_i)-f_i(x_0))^2 j(x_0,dx_i)\right).\notag
\end{align}
This is (\ref{E:fixedbasepoint}), which by (\ref{E:Gammabdj}) is finite, so (i) follows. Using (\ref{E:energyfinitej}) and (\ref{E:Jp}) the inclusion in (ii) and estimate (\ref{E:integrated}) are seen similarly. 

The seminorm properties are clear. To see the last statement in (i) it suffices to note that if $j_1(x_0,\cdot)$ has full support on $N_{1,x_0}$, then $j_p(x_0,\cdot)$ has full support on $N_{p,x_0}\setminus D_{p,x_0}$: Given $U\subset N_{p,x_0}\setminus D_{p,x_0}$ nonempty and open, we can find $(x_0,...,x_p)\in N_p\setminus D_p$ such that $\varepsilon:=\frac12\min_{0\leq i<k\leq p}\varrho(x_i,x_k)>0$ and $B(x_1,\varepsilon)\times ...\times B(x_p,\varepsilon)$ is contained in $U$. By hypothesis this product of balls has positive $j_p(x_0,\cdot)$-measure, and therefore $j_p(x_0,U)>0$. The last statement in (ii) follows similarly.
\end{proof}

\begin{examples}\label{Ex:upperenergies}
Recall Examples \ref{Ex:derivations}.
\begin{enumerate}
\item[(i)] For $p=0$ we have 
\begin{align}
\big\|\overline{g}\delta_0 f\big\|_{L^2(N_1\setminus D_1,J_1)}^2
&=\int_{N_1\setminus D_1}\overline{g}(x_0,x_1)^2(f(x_1)-f(x_0))^2 J_1(d(x_0,x_1))\notag\\
&=\frac12 \int_X\int_X \overline{g}(x_0,x_1)^2(f(x_1)-f(x_0))^2\mathbf{1}_{N_1}(x_0,x_1)\times\notag\\
&\hspace{100pt}\times (j(x_0,dx_1)\mu(dx_0)+j(x_1,dx_0)\mu(dx_1)).\notag
\end{align}
\item[(ii)] For $p=1$ we have 
\begin{align}
\big\|\overline{g}\delta_1\Alt_2(f_1 &\otimes f_2)\big\|_{L^2(N_2\setminus D_2,J_2)}^2\notag\\
&=\frac13\int_{N_2\setminus D_2}\overline{g}(x_0,x_1,x_2)^2\det[(f_i(x_j)-f_i(x_0))_{i,j=1}^2]^2 J_2(d(x_0,x_1,x_2)).\notag
\end{align}
\end{enumerate}
\end{examples}

A class of typical examples arises from Dirichlet forms, \cite{FOT94, GrigHuLau14, MasamuneUemura2011}. 
\begin{examples}\label{Ex:nonlocalDf}
Suppose that $\mu$ is a nonnegative Radon measure on $X$ with full support,
\begin{equation}\label{E:nonlocalDf}
\mathcal{E}(f):=\int_X\int_X (f(x)-f(y))^2j(x,dy)\mu(dx),\quad f\in L^2(X,\mu),
\end{equation}
where $j$ is as in Assumption \ref{A:C} (i) and satisfies $j(x,dy)\mu(dx)=j(y,dx)\mu(dy)$. Let 
\begin{equation}\label{E:maximaldomain}
\mathcal{D}(\mathcal{E})=\{f\in L^2(X,\mu): \mathcal{E}(f)<+\infty\}.
\end{equation}
Then $(\mathcal{E},\mathcal{D}(\mathcal{E}))$ is a Dirichlet form on $X$. If in addition all balls have finite measure and there is some $\varepsilon>0$ such that 
\begin{equation}\label{E:kernelconditions}
\sup_{x\in X}\int_{B(x,\varepsilon)}\varrho(x,y)^2j(x,dy)<+\infty\quad \text{and}\quad \sup_{x\in X}\int_{B(x,\varepsilon)^c}j(x,dy)<+\infty,
\end{equation}
then the algebra $\lip_c(X)$ of compactly supported Lipschitz functions on $X$ satisfies Assumptions \ref{A:C} and \ref{A:CE} in place of $\mathcal{C}$; note that conditions (\ref{E:kernelconditions}) ensure (\ref{E:Gammabdj}). Given $f\in \lip_c(X)$, let $B$ be an open ball containing $K:=\supp f$ and large enough to have $\dist(K,B^c)>\varepsilon$. Then 
\[\mathcal{E}(f)=\int_B\int_X(f(y)-f(x))^2j(x,dy)\mu(dx)+\int_{B^c}\int_K f(y)^2j(x,dy)\mu(dx)<+\infty.\]
By (\ref{E:Gammabdj}) and since $\mu(B)<+\infty$ the first summand is finite. The second is bounded by $\|f\|_{\sup}^2\mu(K)\sup_{y\in K}\int_{B^c} j(y,dx)$, which is finite by (\ref{E:kernelconditions}). This, together with $\lip_c(X)\subset L^2(X,\mu)$, implies that $\lip_c(X)\subset \mathcal{D}(\mathcal{E})$.
\end{examples}

Recall Lemma \ref{L:gradientsplitoff}. We consider the linear extension of (\ref{E:basicmultact}) and record a related
norm estimate for derivatives.
\begin{lemma}\label{L:normsplitoff}
Suppose that Assumptions \ref{A:C} and \ref{A:CE} are satisfied and let $N_\ast$ be a system of 
diagonal neighborhoods. Let $p\geq 0$ be an integer, $F\in \mathcal{C}^{p}(N_{p})$ and $\chi\in \mathcal{C}\oplus\mathbb{R}$. Then 
\begin{multline}
\left\|\delta_{p}(\chi^{\otimes (p+1)}F)\right\|_{L^2(N_{p+1}\setminus D_{p+1},J_{p+1})}\leq \left\|\chi\right\|_{\sup}^{p+1}\left\|\delta_{p}F\right\|_{L^2(N_{p+1}\setminus D_{p+1},J_{p+1})}\notag\\
+(p+1)\left\|\chi\right\|_{\sup}^{p}\Big(\sup_{x\in X}\int_{N_{1,x}\setminus D_{1,x}}(\delta_0\chi(x,y))^2j(x,dy)\Big)^{1/2}\left\|F\right\|_{L^2(N_{p}\setminus D_{p},J_{p})}.
\end{multline}
\end{lemma}
Note that by Lemma \ref{L:elemL2} (ii) the norms of $F$ and $\delta_pF$ on the right hand side are finite.

\begin{proof}
Since $(\delta_{p}(\chi^{\otimes (p+1)}F)(x_0,...,x_{p+1}))^2$ is symmetric in $x_0,...,x_{p+1}$, we have 
\begin{equation}\label{E:simplifiedexp}
\left\|\delta_{p}(\chi^{\otimes (p+1)}F)\right\|_{L^2(N_{p+1}\setminus D_{p+1},J_{p+1})}=\Big(\int_X\left\|\delta_{p}(\chi^{\otimes (p+1)}F)\right\|_{L^2(N_{p+1,x_0}\setminus D_{p+1,x_0},j_{p+1}(x_0,\cdot))}^2\mu(dx_0)\Big)^{1/2};
\end{equation}
similarly for $\delta_p F$ or $F$ in place of $\delta_{p}(\chi^{\otimes (p+1)}F)$. Since an element $F$ of $\mathcal{C}^{p}(N_{p})$ is of the form $F=\sum_{i=1}^m\Alt_{p+1}(f_0^{(i)}\otimes\cdots\otimes f_{p}^{(i)})$, we can employ Lemma \ref{L:gradientsplitoff} to bound (\ref{E:simplifiedexp}) by 
\begin{multline}
\left\|\chi\right\|_{\sup}^{p+1}\left\|\delta_{p}F\right\|_{L^2(N_{p+1}\setminus D_{p+1},J_{p+1})}\notag\\
+\sum_{k=1}^{p+1}\left\|\chi\right\|_{\sup}^{p}\Big(\int_X\int_{N_{1,x_0}\setminus D_{1,x_0}}\cdots\int_{N_{1,x_0}\setminus D_{1,x_0}}(\delta_0\chi(x_0,x_k))^2(F(x_0,\dots,\hat{x}_k,\dots,x_{p+1}))^2\times\notag\\
\times j(x_0,dx_1)\cdots j(x_0,dx_{p+1})\:\mu(dx_0)\Big)^{1/2},\notag
\end{multline}
and each of the multiple integrals in brackets in the second summand can be estimated by 
\[\sup_{x_0\in X}\int_{N_{1,x_0}\setminus D_{1,x_0}}(\delta_0\chi(x_0,x_k))^2j(x_0,dx_k)\:\left\|F\right\|_{L^2(N_{p}\setminus D_{p},J_{p})}^2.\]
\end{proof}

\section{Non-local Hilbert complexes}\label{S:Hilbert}

\subsection{Regularity and density}

Let $(X,\varrho)$ be a locally compact metric space. We use the notation introduced in Sections \ref{S:elem} and \ref{S:kernels} and consider the following density assumption.

\begin{assumption}\label{A:regularity}
The algebra $\mathcal{C}\cap C_c(X)$ is a uniformly dense subalgebra of $C_c(X)$.
\end{assumption}

A generalization of \cite[Lemma 3.1]{H15} shows that the elementary forms are dense in the Hilbert spaces $L^2_a(N_p\setminus D_p,J_p)$.

\begin{proposition}\label{P:dense}
Let Assumptions \ref{A:C}, \ref{A:CE} and \ref{A:regularity} be satisfied, let $N_\ast$ be a system of diagonal neighborhoods and $p\geq 0$ an integer. 
\begin{enumerate}
\item[(i)] The space $\mathcal{C}^p(N_p)$ is dense in $L^2_a(N_p\setminus D_p,J_p)$. 
\item[(ii)] The operators $(\delta_p,\mathcal{C}^p(N_p))$ are densely defined as operators 
from $L^2_a(N_p\setminus D_p,J_p)$ into $L^2_a(N_{p+1}\setminus D_{p+1},J_{p+1})$.
\end{enumerate}
\end{proposition}

\begin{proof}
It suffices to show that if $F\in L^2_a(N_p\setminus D_p, J_p)$ is such that 
\begin{equation}\label{E:annihilate}
\left\langle \Alt_{p+1}(f_0\otimes f_1\otimes ...\otimes f_p),F\right\rangle_{L^2(N_p\setminus D_p,J_p)}=0
\end{equation}
for all $f_1,...,f_p\in \mathcal{C}$ and $f_0\in \mathcal{C}\oplus \mathbb{R}$, then we have 
\begin{equation}\label{E:desiredzero}
F\equiv 0\quad \text{in} \quad L^2(N_p\setminus D_p, J_p).
\end{equation}
Suppose that $f_0,f_1,...,f_p$ are as specified and that (\ref{E:annihilate}) holds. Since $\Alt_{p+1}$ acts as an orthogonal projection in $L^2(N_p\setminus D_p,J_p)$ and $\Alt_{p+1}F=F$, it follows that
\begin{equation}\label{E:intzero}
0=\int_{N_p\setminus D_p} f_0(x_0)f_1(x_1)\cdots f_p(x_p)\:F(x_0,...,x_p)J_p(d(x_0,...,x_p)).
\end{equation}
By Assumption \ref{A:regularity} and Stone-Weierstrass the algebra generated by products $f_0\otimes f_1\otimes ...\otimes f_p$ with $f_0,...,f_p\in \mathcal{C}\cap C_c(X)$ is dense in $C_c(X^{p+1}\setminus D_p)$ and therefore also in $L^2(N_p\setminus D_p,J_p)$. This, together with (\ref{E:intzero}), shows that 
\[0=\left\langle G, F\right\rangle_{L^2(N_p\setminus D_p,J_p)}\quad \text{for all $G\in L^2(N_p\setminus D_p,J_p)$,}\]
and this gives (\ref{E:desiredzero}).
\end{proof}

\subsection{Closed extensions} 

We proceed under the following absolute continuity condition.
\begin{assumption}\label{A:abscont}
The kernel $j$ as in Assumption \ref{A:C} is of the form
\[j(x,dy)=j(x,y)\mu(dy)\]
with a Borel function $j:X^2\to (0,+\infty]$, locally bounded on $X^2\setminus D_1$.
\end{assumption}

Given $p\geq 1$, let $\mu^p$ denote the $p$-fold product of $\mu$ on $X^p$. By Assumption \ref{A:abscont} we have
\[j_p(x_0,d(x_1,...,x_p))=\left(\prod_{\ell=1}^p j(x_0,x_\ell)\right)\mu^p(d(x_1,...,x_p)),\quad x_0\in X,\quad p\geq 1,\]
and 
\begin{equation}\label{E:Jpac}
J_p(d(x_0,...,x_p))=\left(\frac{1}{p+1}\sum_{k=0}^p \prod_{\ell\neq k} j(x_k,x_\ell)\right) \mu^{p+1}(d(x_0,...,x_{p})),\quad p\geq 1.
\end{equation}

Under Assumption \ref{A:abscont} we can extend the derivations $\delta_p$ to closed operators.

\begin{theorem}\label{T:closable} Let Assumptions \ref{A:C}, \ref{A:CE} and \ref{A:abscont} be satisfied.
Let $N_\ast$ be a system of diagonal neighborhoods and let $p\geq 0$ be an integer. 
\begin{enumerate}
\item[(i)] The operators $(\delta_p,\mathcal{C}^p(N_p))$ extend to 
closed linear operators $(\delta_p,\mathcal{D}(\delta_p,N_p))$ from $L^2_a(N_p\setminus D_p,J_p)$ into 
$L^2_a(N_{p+1}\setminus D_{p+1},J_{p+1})$. Given $F\in \mathcal{D}(\delta_p,N_p)$, we have
\begin{equation}\label{E:differenceopagain}
\delta_{p} F(x_0, \dotsc, x_{p+1})= \sum_{i=0}^{p+1} (-1)^i F(x_0, \dotsc, \hat{x}_i, \dotsc, x_{p+1})\quad \text{}
\end{equation}
in the $J_{p+1}$-a.e. sense. If in addition Assumption \ref{A:regularity} holds, then the operators $(\delta_p,\mathcal{D}(\delta_p,N_p))$ are densely defined.
\item[(ii)] We have 
\begin{equation}\label{E:chainprop}
\im \delta_{p}|_{\mathcal{D}(\delta_{p},N_{p})} \subset \mathcal{D}(\delta_{p+1},N_{p+1})\quad \text{and}\quad \delta_{p+1}\circ \delta_p=0.
\end{equation}
The sequence 
\begin{equation}\label{E:Hilbertcomplex}
0\longrightarrow \mathcal{D}(\delta_0,N_0)\stackrel{\delta_0}{\longrightarrow} \mathcal{D}(\delta_1,N_1)\stackrel{\delta_1}{\longrightarrow} \mathcal{D}(\delta_2,N_2) \stackrel{\delta_2}{\longrightarrow} ...
\end{equation}
is a cochain complex.
\end{enumerate}
\end{theorem}

Similarly as in \cite{BL92} one could call the complex $(\mathcal{D}(\delta_\ast,N_\ast),\delta_\ast)$ in (\ref{E:Hilbertcomplex}) a \emph{Hilbert complex}. It is 'non-local' in the sense that the operators $\delta_p$ are non-local. 

For any integer $p\geq 0$ the $p$-th cohomology of (\ref{E:Hilbertcomplex}) is 
\begin{equation}\label{E:closedcoho}
H^p\mathcal{D}(\delta_\ast,N_\ast):=\ker \delta_p|_{\mathcal{D}(\delta_p,N_p)} / \im \delta_{p-1}|_{\mathcal{D}(\delta_{p-1},N_{p-1})}.
\end{equation}

\begin{remark}\mbox{}
\begin{enumerate}
\item[(i)] In \cite{BL92} the ambient Hilbert spaces of sufficiently high order $p$ were assumed to be trivial; we do not make this assumption here.
\item[(ii)] Obviously the complex (\ref{E:elemcomplex}) of elementary functions is a subcomplex of (\ref{E:Hilbertcomplex}).
\end{enumerate} 
\end{remark}

To prove Theorem \ref{T:closable} we consider quadratic forms associated with the operators $\delta_p$. Note first that for each $p\geq 0$ the application of $\delta_p$ to an element $F$ of $L^2(N_{p}\setminus D_{p},J_{p})$ gives a well-defined $\mu^{p+2}$-equivalence class $\delta_pF$ of Borel functions on $N_{p+1}\setminus D_{p+1}$; this is straightforward from (\ref{E:def-coboundary}) and Assumption \ref{A:abscont}. Therefore
\begin{equation}\label{E:quadformdelta}
Q_p(F):=\left\|\delta_p F\right\|_{L^2(N_{p+1}\setminus D_{p+1},J_{p+1})}^2,\quad F\in L^2(N_{p}\setminus D_{p},J_{p}),
\end{equation}
defines a quadratic form
\[Q_p:L^2(N_p\setminus D_{p},J_{p})\to [0,+\infty].\]
By Lemma \ref{L:elemL2} we have $Q_{p}(F)<+\infty$ for any $F\in \mathcal{C}^p(N_p)$.

\begin{remark}\label{R:lsc}
Recall that if $I:B\to [0,+\infty]$ is a lower semicontinuous functional on a Banach space $B$ and $V$ is a subspace of $B$ on which $I$ is finite, then $(I,V)$ is closable. 
\end{remark}

\begin{proposition}\label{P:closable}
 Let Assumptions \ref{A:C}, \ref{A:CE} and \ref{A:abscont} be satisfied. Let $N_\ast$ be a system of diagonal neighborhoods and let $p\geq 0$ be an integer. 
 
Then the quadratic form $Q_p$ defined in (\ref{E:quadformdelta}) is lower semicontinuous on $L^2_a(N_p\setminus D_p, J_p)$. Its restriction $(Q_p,\mathcal{C}^p(N_p))$ to $\mathcal{C}^p(N_p)$ is closable.
\end{proposition}

\begin{proof}
Suppose that $(F_j)_j\subset L^2_a(N_p\setminus D_p, J_p)$ converges to $F$ in $L^2_a(N_p\setminus D_p, J_p)$ and let $(j_k)_k$ be such that 
\begin{equation}\label{E:select}
\lim_{k\to \infty} Q_p(F_{j_k})=\liminf_{j\to \infty}Q_p(F_j);
\end{equation}
we may assume the right hand side is finite. Passing to a further subsequence of $(F_{j_k})_k$ and relabelling, we may assume that 
\[\lim_{k\to \infty} F_{j_k}(x_0,...,x_p)=F(x_0,...,x_p)\]
for all $(x_0,...,x_p)\in (N_p\setminus D_p)\setminus Z_{p}$, where $Z_{p}\subset N_p$ is a set of measure zero for $J_p$. Since by (\ref{E:Jpac}) the density of $J$ with respect to $\mu^{p+1}$ is strictly positive, $Z_{p}\subset N_p$ is also of zero measure for $\mu^{p+1}$. Now let 
\[Z_{p+1}':=\bigcup_{i=0}^{p+1}\left\lbrace (x_0,...,x_{p+1})\in N_{p+1}\setminus D_{p+1}: (x_0,...,\hat{x}_i,...,x_{p+1})\in Z_{p}\right\rbrace.\]
Then $\mu^{p+2}(Z_{p+1}')=0$, and for all $(x_0,...,x_{p+1})\in (N_{p+1}\setminus D_{p+1})\setminus Z_{p+1}'$ we have 
\[\lim_{k\to \infty} \delta_pF_{j_k}(x_0,...,x_{p+1})=\delta_pF(x_0,...,x_{p+1}).\]
Since $J_{p+1}\ll \mu^{p+2}$ by (\ref{E:Jpac}), Fatou's lemma now implies that 
\begin{align}
\lim_{k\to \infty} Q_p(F_{j_k})&=\liminf_{k\to \infty} \int_{N_{p+1}\setminus D_{p+1}}(\delta_p F_{j_k}(x_0,...,x_{p+1}))^2 J_{p+1}(d(x_0,...,x_{p+1}))\notag\\
&\geq \int_{N_{p+1}\setminus D_{p+1}}(\delta_p F(x_0,...,x_{p+1}))^2 J_{p+1}(d(x_0,...,x_{p+1}))\notag\\
&= Q_p(F),\notag
\end{align}
and combining with (\ref{E:select}), the lower semicontinuity of $Q_p$ is observed. By Remark \ref{R:lsc} the form $(Q_p,\mathcal{C}^p(N_p))$ is closable. 
\end{proof}

The preceding now yields a quick proof of Theorem \ref{T:closable}.

\begin{proof}
Since the closability of $(\delta_p,\mathcal{C}^p(N_p))$ is equivalent to the closability of $(Q_p,\mathcal{C}^p(N_p))$, the first statement of Theorem \ref{T:closable} is immediate from Proposition \ref{P:closable}.

To see (ii), let $F\in \mathcal{D}(\delta_p,N_p)$ and let $(F_n)_n\subset \mathcal{C}^p(N_p)$ be such that $\lim_{n\to \infty}F_n=F$ in $L^2(N_p\setminus D_p,J_p)$ and $\lim_{n\to \infty}\delta_pF_n=\delta_pF$ in $L^2(N_{p+1}\setminus D_{p+1},J_{p+1})$. The lower semicontinuity of $Q_{p+1}$ and (\ref{E:coboundary-square-zero}) imply that 
\begin{multline}
\left\|\delta_{p+1}(\delta_pF)\right\|_{L^2(N_{p+2}\setminus D_{p+2},J_{p+2})}^2=Q_{p+1}(\delta_pF)
\leq \liminf_{n\to \infty} Q_{p+1}(\delta_pF_n)\notag\\
=\liminf_{n\to \infty}\left\|\delta_{p+1}(\delta_pF_n)\right\|_{L^2(N_{p+2}\setminus D_{p+2},J_{p+2})}^2=0.
\end{multline}
Consequently $\delta_pF\in \ker \delta_{p+1}\subset \mathcal{D}(\delta_{p+1},N_p)$, and $\delta_{p+1}\circ \delta_p=0$.
\end{proof}

We write $(\mathcal{Q}_p, \mathcal{D}(\delta_{p},N_p))$ to denote the closure of $(Q_p,\mathcal{C}^p(N_p))$ in $L^2_a(N_p\setminus D_p, J_p)$. By closedness the domain $\mathcal{D}(\delta_{p},N_p)$, endowed with the scalar product
\begin{equation}\label{E:Qp1}
\left\langle F,G\right\rangle_{\mathcal{D}(\delta_{p},N_p)}:=\left\langle F,G\right\rangle_{L^2(N_p\setminus D_p,J_p)}+\mathcal{Q}_p(F,G),\quad F,G\in \mathcal{D}(\delta_{p},N_p),
\end{equation}
is a Hilbert space. We write $\left\|\cdot\right\|_{\mathcal{D}(\delta_{p},N_p)}$ for the associated Hilbert norm.

\begin{remark}\label{R:notnecessarilymax}
We point out that in general $\mathcal{D}(\delta_{p},N_p)$ may be smaller than the maximal domain $\{F\in L^2_a(N_p\setminus D_p, J_p):\ Q_p(F)<+\infty\}$ of $Q$.
\end{remark}

\begin{corollary}\label{C:closed}
Let Assumptions \ref{A:C}, \ref{A:CE}, \ref{A:regularity} and \ref{A:abscont} be satisfied. Let $N_\ast$ be a system of diagonal neighborhoods and let $p\geq 0$ be an integer. Then 
$(\mathcal{Q}_p, \mathcal{D}(\delta_{p},N_p))$ 
is a densely defined closed quadratic form on $L^2_a(N_p\setminus D_p, J_p)$, and 
\begin{equation}\label{E:extendedform}
\mathcal{Q}_p(F)=\left\|\delta_p F\right\|_{L^2(N_{p+1}\setminus D_{p+1},J_{p+1})}^2,\quad F \in \mathcal{D}(\delta_{p},N_p).
\end{equation}
\end{corollary}

\begin{remark}\label{R:QisDirichletform}
Let the hypotheses of Corollary \ref{C:closed} be satisfied, assume that $\mathcal{C}\subset C_c(X)$ and that $X$ is separable. Then polarization of the quadratic form $(\mathcal{Q}_0, \mathcal{D}(\delta_{0},N_0))$ yields a regular Dirichlet form on $X=N_0$ in the sense of \cite{FOT94}, and this form is purely non-local. Its generator, that is, the unique nonpositive definite self-adjoint operator $(\mathcal{L}_0,\mathcal{D}(\mathcal{L}_0))$ on $L^2(X,\mu)$ satisfying 
\[\left\langle \mathcal{L}_0f,g\right\rangle_{L^2(X,\mu)}=-\mathcal{Q}_0(f,g),\quad f\in \mathcal{D}(\mathcal{L}_0),\quad g\in \mathcal{D}(\delta_{0},N_0),\]
is of the form 
\begin{equation}\label{E:nonlocop}
\mathcal{L}_0f(x)=\int_{N_1\setminus \{x\}}(f(y)-f(x))(j(x,y)+j(y,x))\mu(dy),\quad f\in \mathcal{D}(\mathcal{L}_0),
\end{equation}
understood in the $\mu$-a.e. sense. 
\end{remark}

\begin{examples} Suppose that $j$ is symmetric and 
\begin{equation}\label{E:frackerneld}
c^{-1}\varrho(x,y)^{-d-\alpha}\leq  j(x,y) \leq c\:\varrho(x,y)^{-d-\alpha},\quad x,y\in X,
\end{equation}
with some fixed $d>0$, $\alpha\in (0,2)$ and $c>1$. If $\varepsilon>0$ and $N_\ast=N_\ast(\varepsilon)$, then $\mathcal{Q}_0$ is comparable to the quadratic form 
\[f\mapsto \int\int_{\varrho(x,y)<\varepsilon} (f(x)-f(y))^2\varrho(x,y)^{-d-\alpha}\mu(dx)\mu(dy),\]
and $\mathcal{L}_0$ in (\ref{E:nonlocop}) is comparable to a \emph{truncated fractional Laplacian of order $\alpha/2$} on $X$.
\end{examples}

Recall Lemmas \ref{L:gradientsplitoff} and \ref{L:normsplitoff}. We observe that the linear extension of (\ref{E:basicmultact}) becomes a multiplier on $\mathcal{D}(\delta_p,N_p)$.
\begin{corollary}\label{C:Sobomult}
Let Assumptions \ref{A:C}, \ref{A:CE} and \ref{A:abscont} be satisfied, let $N_\ast$ be a system of diagonal neighborhoods and $\chi\in \mathcal{C}\oplus \mathbb{R}$. Then for any $p\geq 0$ the multiplication $F\mapsto \chi^{\otimes (p+1)}F$, $F\in\mathcal{C}^p(N_p)$, extends to a bounded linear operator on $\mathcal{D}(\delta_p,N_p)$, and we have
\begin{equation}\label{E:Sobomult}
\left\|\chi^{\otimes(p+1)}F\right\|_{\mathcal{D}(\delta_p,N_p)}\leq c_{p,\chi}\:\left\|F\right\|_{\mathcal{D}(\delta_p,N_p)},\quad F\in \mathcal{D}(\delta_p,N_p),
\end{equation}
where
\[c_{\chi,p}=\left\|\chi\right\|_{\sup}^p\Big(1+\left\|\chi\right\|_{\sup}+(p+1)\sup_{x\in X}\Big(\int_{N_{1,x}\setminus D_{1,x}}(\delta_0\chi(x,y))^2j(x,y)\mu(dy)\Big)^{1/2}\Big).\]
\end{corollary}

\begin{proof}
By the trivial bound 
\[\left\|\chi^{\otimes (p+1)}F\right\|_{L^2(N_p\setminus D_p,J_p)}\leq \left\|\chi\right\|_{\sup}^{p+1}\left\|F\right\|_{L^2(N_p\setminus D_p,J_p)}\] 
and Lemma \ref{L:normsplitoff}, estimate (\ref{E:Sobomult}) is seen to hold for all $F\in \mathcal{C}^p(N_p)$. Now let $F\in \mathcal{D}(\delta_p,N_p)$ and let $(F_n)_n\subset \mathcal{C}^p(N_p)$ be such that $\lim_{n\to \infty} \left\|F_n-F\right\|_{\mathcal{D}(\delta_p,N_p)}=0$. Clearly 
\[\lim_{n\to \infty}\chi^{\otimes (p+1)}F_n=\chi^{\otimes (p+1)}F\quad\text{in $L^2(N_p\setminus D_p,J_p)$.}\] 
By Lemma \ref{L:normsplitoff} the sequence $(\delta_p(\chi^{\otimes (p+1)}F_n))_n$ is Cauchy in $L^2(N_{p+1}\setminus D_{p+1},J_{p+1})$, so it has a limit $G$. Since $\delta_p$ is closed, we have $\chi^{\otimes (p+1)}F\in \mathcal{D}(\delta_p,N_p)$ and $\delta_p(\chi^{\otimes (p+1)}F)=G$. Estimate (\ref{E:Sobomult}) now follows easily.
\end{proof}

\subsection{Remarks on invariance} 

Choices of comparable kernels $j$ in Assumption \ref{A:abscont} lead to the same closures and therefore to the same complexes.
\begin{proposition}
Suppose that both $j$ and $j'$ are kernels such that Assumptions \ref{A:C}, \ref{A:CE} and \ref{A:abscont} are satisfied and let $N_\ast$ be a system of diagonal neighborhoods. If there is a constant $c>1$ such that $c^{-1}j(x,y)\leq j'(x,y)\leq c\:j(x,y)$, $(x,y)\in N_1\setminus D_1$, then the resulting complexes $(\mathcal{D}(\delta_\ast,N_\ast),\delta_\ast)$ and $(\mathcal{D}'(\delta_\ast,N_\ast),\delta_\ast)$
coincide.
\end{proposition}

Recall Proposition \ref{P:pullback} and the notation used there. Suppose that $(\tilde{X},\tilde{\varrho})$ is another (locally compact) metric space with associated data $\tilde{\mathcal{C}}$, $\tilde{j}$, $\tilde{\mu}$ and that $\varphi:\tilde{X}\to X$ is a given Borel map. Under suitable conditions the resulting Hilbert complexes on $X$ and $\tilde{X}$ are isomorphic. We write $\varphi^\ast j(\tilde{x},\tilde{y}):=j (\varphi(\tilde{x}),\varphi(\tilde{y}))$ and $\varphi_\ast \tilde{\mu}:=\tilde{\mu}\circ \varphi^{-1}$, and we denote the Hilbert spaces obtained from $\tilde{\mathcal{C}}^p(\varphi^\ast N_p)$ by taking closures as in Theorem \ref{T:closable} by $\tilde{\mathcal{D}}(\delta_p,\varphi^\ast N_p)$.

\begin{theorem}\label{T:invariance}
Suppose that Assumptions \ref{A:C}, \ref{A:CE} and \ref{A:abscont} are satisfied for $X$, $\mathcal{C}$, $j$, $\mu$ and also for $\tilde{X}$, $\tilde{\mathcal{C}}$, $\tilde{j}$, $\tilde{\mu}$. Let $N_\ast$ be a system of diagonal neighborhoods for $X$. Assume that $\varphi$ is a homeomorphism from $\tilde{X}$ onto $X$, 
\begin{equation}\label{E:kernelsequiv}
c_j^{-1}\varphi^\ast j\leq \tilde{j}\leq c_j\:\varphi^\ast j,
\end{equation}
and
\begin{equation}\label{E:measuresequiv}
c_\mu^{-1}\mu\leq \varphi_\ast \tilde{\mu}\leq c_\mu \mu
\end{equation}
with constants $c_j>1$, $c_\mu>1$, and that $\varphi^\ast \mathcal{C}\subset \tilde{\mathcal{C}}$ is dense in $\tilde{\mathcal{C}}$ with respect to $\left\|\cdot\right\|_{\tilde{\mathcal{D}}(\delta_0,\varphi^\ast N_0)}$. Then for any integer $p\geq 0$ the map $\varphi^\ast:\mathcal{D}(\delta_p,N_p)\to \tilde{\mathcal{D}}(\delta_p,\varphi^\ast N_p)$ is an isomorphism of equivalently normed spaces, and 
\begin{equation}\label{E:comparenorms}
c_p^{-1}\left\|F\right\|_{\mathcal{D}(\delta_p,N_p)}\leq \left\|\varphi^\ast F\right\|_{\tilde{\mathcal{D}}(\delta_p,\varphi^\ast N_p)}\leq c_p\:\left\|F\right\|_{\mathcal{D}(\delta_p,N_p)},\quad F\in \mathcal{D}(\delta_p,N_p),
\end{equation}
with a constant $c_p>1$. Moreover, $\varphi^\ast:\mathcal{D}(\delta_\ast,N_\ast)\to \tilde{\mathcal{D}}(\delta_\ast,\varphi^\ast N_\ast)$ is an isomorphism of cochain complexes.  
\end{theorem}

\begin{proof}
It is quickly seen that $\tilde{D}_p$, defined similarly as $D_p$, equals $\varphi^\ast D_p:=\{(\tilde{x}_0,...,\tilde{x}_p)\in \tilde{X}^{p+1}: (\varphi(\tilde{x}_0),...,\varphi(\tilde{x}_p))\in D_p\}$. Using (\ref{E:kernelsequiv}), change of variables and (\ref{E:measuresequiv}) we obtain
\begin{align}
\int_{\varphi^\ast N_p\setminus \varphi^\ast D_p}&(\varphi^\ast F(\tilde{x}_0,...,\tilde{x}_p))^2\tilde{j}(\tilde{x}_0,\tilde{x}_1)\cdots \tilde{j}(\tilde{x}_0,\tilde{x}_p)\tilde{\mu}^{p+1}(d(\tilde{x}_0,...,\tilde{x}_p))\notag\\
&\leq c_j^p \int_{\varphi^\ast N_p\setminus \varphi^\ast D_p}(\varphi^\ast F(\tilde{x}_0,...,\tilde{x}_p))^2\varphi^\ast j(\tilde{x}_0,\tilde{x}_1)\cdots \varphi^\ast j(\tilde{x}_0,\tilde{x}_p)\tilde{\mu}^{p+1}(d(\tilde{x}_0,...,\tilde{x}_p))\notag\\
&=c_j^p \int_{N_p\setminus D_p} (F(x_0,...,x_p))^2 j(x_0,x_1)\cdots j(x_0,x_p)\varphi_\ast \tilde{\mu}^{p+1}(d(x_0,...,x_p))\notag\\
&\leq c_\mu^{p+1}c_j^p \int_{N_p\setminus D_p} (F(x_0,...,x_p))^2 j(x_0,x_1)\cdots j(x_0,x_p)\mu^{p+1}(d(x_0,...,x_p))\notag
\end{align}
and an analogous lower bound with the reciprocal constant. Since $\delta_p\varphi^\ast F=\varphi^\ast\delta_p F$, similar estimates are seen to hold with $\delta_p \varphi^\ast F$ and $\delta_p F$ in place of $\varphi^\ast F$ and $F$, respectively. As in (\ref{E:simplifiedexp}) these estimates carry over to the $L^2$-norms with respect to $J_p$ and $J_{p+1}$ by the symmetry of the integrands. Combining, we arrive at (\ref{E:comparenorms}), and this implies that $\varphi^\ast$ preserves the property to be a Cauchy sequence and therefore induces isomorphisms as stated.
\end{proof}

Let $d>0$. Recall that a nonnegative Radon measure $\mu$ on $X$ is said to be \emph{(Ahlfors) $d$-regular}, \cite[p. 62]{Heinonen}, if there is a constant $c>1$ such that 
\begin{equation}\label{E:dregular}
c^{-1}r^d\leq \mu(B(x,r))\leq c\:r^d,\quad x\in \supp \mu,\ 0<r<\diam(\supp \mu).
\end{equation}
If in this situation $\supp\mu=X$, then $d$ equals the Hausdorff dimension of $X$.

A special case of Theorem \ref{T:invariance} gives the following Lipschitz invariance property.
\begin{corollary}\label{C:invariance}
Suppose that $\varphi$ is a bi-Lipschitz map from $\tilde{X}$ onto $X$ and that both $\mu$ and $\tilde{\mu}$ are $d$-regular and of full support. Let $\alpha\in (0,2)$ and assume that $j$ is as in (\ref{E:frackerneld}) and $\tilde{j}$, too, but with $\tilde{\varrho}$ in place of $\varrho$. 
\begin{enumerate}
\item[(i)] If $\mathcal{C}$ satisfies Assumptions \ref{A:C} (ii) and \ref{A:CE}, then so does $\varphi^\ast\mathcal{C}$.
\item[(ii)] If $\tilde{\mathcal{C}}$ satisfies Assumptions \ref{A:C} (ii), \ref{A:CE} and \ref{A:abscont}, $N_\ast$ is a system of diagonal neighborhoods for $X$ and $\varphi^\ast \mathcal{C}$ is dense in $\tilde{\mathcal{C}}$ with respect to $\left\|\cdot\right\|_{\tilde{\mathcal{D}}(\delta_0,\varphi^\ast N_0)}$, then the norm estimates (\ref{E:comparenorms}) hold and the cochain complexes $\mathcal{D}(\delta_\ast,N_\ast)$ and $\tilde{\mathcal{D}}(\delta_\ast,\varphi^\ast N_\ast)$ are isomorphic.
\end{enumerate}
\end{corollary}
\begin{proof}
To verify (\ref{E:measuresequiv}) one can use the fact that $\mu$ is comparable to the $d$-dimensional Hausdorff measure, \cite[Exercise 8.11]{Heinonen}, the other hypotheses of Theorem \ref{T:invariance} are easily seen.
\end{proof}

\subsection{Non-local Hodge Laplacians}

The next result on associated coderivations $\delta_p^\ast$ follows from general theory, \cite[Theorem 5.3]{Weidmann}.

\begin{corollary}\label{C:adjoint} Let Assumptions \ref{A:C}, \ref{A:CE}, \ref{A:regularity} and \ref{A:abscont} be satisfied. Let $N_\ast$ be a system of diagonal neighborhoods and let $p\geq 0$ be an integer. The adjoint $(\delta_p^\ast,\mathcal{D}(\delta_p^\ast,N_{p+1}))$ of $(\delta_p,\mathcal{D}(\delta_p,N_p))$ is a densely defined closed linear operator
from $L^2_a(N_{p+1}\setminus D_{p+1},J_{p+1})$ into $L^2_a(N_{p}\setminus D_p,J_p)$.
It is characterized by the identity 
\[\left\langle \delta_p F,G\right\rangle_{L^2(N_{p+1}\setminus D_{p+1},J_{p+1})}=\left\langle F, \delta_p^\ast G\right\rangle_{L^2(N_{p}\setminus D_{p},J_{p})},\quad F\in \mathcal{D}(\delta_p, N_p),\ G\in \mathcal{D}(\delta_p^\ast,N_{p+1}),\]
and satisfies
\begin{equation}\label{E:cochainprop}
\im \delta_p^\ast|_{\mathcal{D}(\delta_p^\ast,N_{p+1})}\subset \mathcal{D}(\delta_{p-1}^\ast,N_p)\quad\text{ and }\quad \delta_{p-1}^\ast\circ \delta_p^\ast=0.
\end{equation}
Moreover, $(\mathcal{Q}_p^\ast, \mathcal{D}(\delta_{p}^\ast,N_{p+1}))$, defined by 
\begin{equation}\label{E:extendedcoform}
\mathcal{Q}_p^\ast(F)=\left\|\delta_p^\ast F\right\|_{L^2(N_{p}\setminus D_{p},J_{p})}^2,\quad F \in \mathcal{D}(\delta_{p}^\ast,N_{p+1}),
\end{equation}
is a densely defined closed quadratic form on $L^2(N_{p+1}\setminus D_{p+1},J_{p+1})$.
\end{corollary}

The complex 
\[0\longleftarrow \mathcal{D}(\delta_{-1},N_0)\stackrel{\delta_0^\ast}{\longleftarrow} \mathcal{D}(\delta_0^\ast,N_1)\stackrel{\delta_1^\ast}{\longleftarrow} \mathcal{D}(\delta_1^\ast,N_2) \stackrel{\delta_2^\ast}{\longleftarrow} ...\]
is the dual complex of (\ref{E:Hilbertcomplex}).

\begin{examples}
Recall (\ref{E:nonlocop}). Given $f\in \mathcal{D}(\delta_0,N_0)$ we have $\delta_0 f\in \mathcal{D}(\delta_0^\ast,N_1)$ if and only if $f\in \mathcal{D}(\mathcal{L}_0)$. In this case, we have $\mathcal{L}_0f=-\delta_0^\ast\delta_0 f$. 
\end{examples}

As usual we set $\delta_{-1}:=0$; then also $\delta_{-1}^\ast=0$. For any integer $p\geq 0$ let 
\[\mathbb{H}_p:=\ker \delta_p\cap \ker \delta_{p-1}^\ast.\]
The following \emph{weak Hodge decomposition}, \cite{Gaffney55, Kodaira49}, into orthogonal closed subspaces is straightforward, see \cite[Lemma 2.1]{BL92}.

\begin{corollary}\label{C:weakHodge}
Let Assumptions \ref{A:C}, \ref{A:CE}, \ref{A:regularity} and \ref{A:abscont} be satisfied. Let $N_\ast$ be a system of diagonal neighborhoods.
For any integer $p\geq 0$ the space $L^2_a(N_{p}\setminus D_p,J_p)$ admits the orthogonal decomposition
\begin{equation}\label{E:Hodgedecomp}
L^2_a(N_{p}\setminus D_p,J_p)=\mathbb{H}_p\oplus \overline{\mathrm{im}\:\delta_{p-1}}\oplus \overline{\mathrm{im}\:\delta_{p}^\ast};
\end{equation}
here the closures are taken with respect to the norm in $L^2(N_{p}\setminus D_p,J_p)$.
\end{corollary}

The preceding can be used to introduce non-local analogs of Hodge Laplacians. Recall
(\ref{E:extendedform}) and (\ref{E:extendedcoform}). For any integer $p\geq 0$ consider the subspace
\[\mathcal{D}(\mathbb{D}_p,N_p):=\mathcal{D}(\delta_p,N_p)\cap \mathcal{D}(\delta_{p-1}^\ast,N_p)\]
of $L^2_a(N_{p}\setminus D_p,J_p)$ and set 
\begin{equation}\label{E:Dirichletint}
\mathbb{D}_p(F):=\mathcal{Q}_p(F)+ \mathcal{Q}_{p-1}^\ast(F),\quad F \in \mathcal{D}(\mathbb{D}_p,N_p).
\end{equation}

The quadratic forms  in (\ref{E:Dirichletint}) may be viewed as non-local analogs of the \emph{Dirichlet integrals} for differential forms, \cite{Gaffney55}. The following is a consequence of Corollaries \ref{C:closed} and \ref{C:adjoint}.

\begin{proposition}\label{P:HodgeLaplace}
Let Assumptions \ref{A:C}, \ref{A:CE}, \ref{A:regularity} and \ref{A:abscont} be satisfied and let $N_\ast$ be a system of diagonal neighborhoods.
The quadratic form $(\mathbb{D}_p,\mathcal{D}(\mathbb{D}_p,N_p))$ is closed and densely defined.
\end{proposition}

By polarization $\mathbb{D}_p$ may be seen as a closed and densely defined nonnegative definite symmetric bilinear form on $L^2_a(N_p\setminus D_p,J_p)$. For any $p\geq 0$ let $(\mathbb{L}_p,\mathcal{D}(\mathbb{L}_p))$ be the unique nonnegative definite self-adjoint operator associated with $(\mathbb{D}_p,\mathcal{D}(\mathbb{D}_p))$ in the sense that 
\[\left\langle \mathbb{L}_p F, G\right\rangle_{L^2(X^{p+1}\setminus D_p,J_p)}=\mathbb{D}_p(F,G),\quad F \in \mathcal{D}(\mathbb{L}_p),\quad G\in \mathcal{D}(\mathbb{D}_p).\]
The operator $(\mathbb{L}_p,\mathcal{D}(\mathbb{L}_p))$ may be viewed as a \emph{non-local analog of the Hodge Laplacian}, \cite{Kodaira49, Gaffney55}. Formally, it satisfies the identity
\[\mathbb{L}_p=\delta_{p-1}\delta_{p-1}^\ast+\delta_p^\ast\delta_p.\]
Note that for $p=0$ we have $\mathbb{L}_0=-\mathcal{L}_0$ with $\mathcal{L}_0$ as defined in (\ref{E:nonlocop}).

\begin{proof}[Proof of Proposition \ref{P:HodgeLaplace}]
Note first that by Proposition \ref{P:closable} and Corollary \ref{C:adjoint} the forms $(\mathcal{Q}_p, \mathcal{D}(\delta_{p},N_p))$ 
and $(\mathcal{Q}_{p-1}^\ast, \mathcal{D}(\delta_{p-1}^\ast,N_p))$ are densely defined closed quadratic forms. Their closedness implies that also $(\mathbb{D}_p,\mathcal{D}(\mathbb{D}_p,N_p))$ must be closed. 

It remains to show that $\mathcal{D}(\mathbb{D}_p,N_p)$ is dense in $L^2_a(N_p\setminus D_p,J_p)$. To see this, note first that the nonnegative definite self-adjoint operators uniquely associated with $\mathcal{Q}_p^\ast$ and $\mathcal{Q}_{p-1}$ are 
\[(\delta_p\delta_p^\ast,\mathcal{D}(\delta_p\delta_p^\ast,N_{p+1}))\quad \text{and}\quad (\delta_{p-1}^\ast\delta_{p-1},\mathcal{D}(\delta_{p-1}^\ast\delta_{p-1},N_{p-1})),\]
respectively. In particular, $\mathcal{D}(\delta_{p}\delta_{p}^\ast,N_{p+1})$ is dense in the Hilbert space $\mathcal{D}(\delta_{p}^\ast,N_{p+1})$ and $\mathcal{D}(\delta_{p-1}^\ast\delta_{p-1},N_{p-1})$ is dense in the Hilbert space $\mathcal{D}(\delta_{p-1},N_{p-1})$.
This implies that 
\begin{equation}\label{E:denseimages}
\text{$\delta^\ast_p(\mathcal{D}(\delta_{p}\delta_{p}^\ast,N_{p+1}))$ is dense in $\overline{\mathrm{Im}\:\delta_p^\ast}$ and $\delta_{p-1}(\mathcal{D}(\delta_{p-1}^\ast\delta_{p-1},N_{p-1}))$ is dense in $\overline{\mathrm{Im}\:\delta_{p-1}}$   }
\end{equation}
with respect to the norm in $L^2(N_p\setminus D_p,J_p)$. 

Given $G\in \mathcal{D}(\delta_{p}^\ast,N_{p+1})$, the co-exact function $\delta_{p}^\ast G$ is in $\mathcal{D}(\delta_{p},N_p)$ if and only if $G$ is in $\mathcal{D}(\delta_{p}\delta_{p}^\ast,N_{p+1})$. Given $G\in \mathcal{D}(\delta_{p-1},N_{p-1})$, the exact function $\delta_{p-1} G$ is in $\mathcal{D}(\delta_{p-1}^\ast,N_p)$ if and only if $G$ is in $\mathcal{D}(\delta_{p-1}^\ast\delta_{p-1},N_{p-1})$.
Together with (\ref{E:chainprop}) and (\ref{E:cochainprop}) this implies that 
\begin{equation}\label{E:inclusions}
\delta_{p}^\ast(\mathcal{D}(\delta_p\delta_p^\ast,N_{p+1}))\subset \mathcal{D}(\mathbb{D}_p,N_p)\quad \text{and}\quad \delta_{p-1}(\mathcal{D}(\delta_{p-1}^\ast\delta_{p-1},N_{p-1}))\subset \mathcal{D}(\mathbb{D}_p,N_p).
\end{equation}

Now suppose that $F\in L^2_a(N_p\setminus D_p,J_p)$ is such that 
\[\left\langle F,G\right\rangle_{L^2(N_{p}\setminus D_p,J_p)}=0\quad \text{for all $G\in \mathcal{D}(\mathbb{D}_p,N_p)$}.\] 
Since obviously $\mathbb{H}_p\subset \mathcal{D}(\mathbb{D}_p,N_p)$, we have 
\begin{equation}\label{E:killharmonics}
\left\langle F,G\right\rangle_{L^2(N_p\setminus D_p,J_p)}=0\quad \text{for all $G\in \mathbb{H}_p$.}
\end{equation}
By (\ref{E:inclusions}) also 
\[\left\langle F, \delta_{p}^\ast G\right\rangle_{L^2(N_p\setminus D_p,J_p)}=0 \quad \text{for all $G\in \mathcal{D}(\delta_p\delta_p^\ast,N_{p+1})$}\]
and 
\[\left\langle F, \delta_{p-1} G\right\rangle_{L^2(N_p\setminus D_p,J_p)}=0  \quad \text{for all $G\in \mathcal{D}(\delta_{p-1}^\ast\delta_{p-1},N_{p-1})$.}\]
By (\ref{E:denseimages}) this means that 
\begin{equation}\label{E:killimages}
\left\langle F,G\right\rangle_{L^2(N_p\setminus D_p,J_p)}=0\quad \text{for all $G\in \overline{\mathrm{Im}\:\delta_p^\ast}$ and all $G\in \overline{\mathrm{Im}\:\delta_{p-1}}$.}
\end{equation}
The combination of (\ref{E:Hodgedecomp}), (\ref{E:killharmonics}) and (\ref{E:killimages}) shows that $F=0$ in $L^2_a(N_p\setminus D_p,J_p)$.
\end{proof}

\subsection{Approximation results}

The following estimate complements (\ref{E:integrated}).

\begin{lemma}\label{L:L2dom}
Let Assumptions \ref{A:C}, \ref{A:CE} and \ref{A:abscont} be satisfied and let $N_\ast$ be a system of diagonal neighborhoods. Then for any integer $p\geq 1$ and any $g\in \mathcal{C}\oplus \mathbb{R}$, $f_1,...,f_p\in\mathcal{C}$ we have 
\begin{multline}\label{E:integratedL2}
\left\|\overline{g}\delta_{p-1}\Alt_p(f_1\otimes...\otimes f_p)\right\|_{L^2(N_p\setminus D_p,J_p)}\\
\leq \frac{3p+1}{p+1}\left\|g\right\|_{\mathcal{D}(\delta_0,N_0)}\:\prod_{i=1}^p\Big(\left\|f_i\right\|_{\sup}+\big(\sup_{x\in X}\int_{N_{1,x}\setminus D_{1,x}}(f_i(y)-f_i(x))^2j(x,y)\mu(dy)\big)^{1/2}\Big).
\end{multline}
\end{lemma}

\begin{proof}
As in (\ref{E:simplifiedexp}) we can use the symmetry of $(\overline{g}\delta_{p-1}\Alt_p(f_1\otimes...\otimes f_p))^2$ to see that the left hand side of (\ref{E:integratedL2}) equals 
\begin{multline}
\left(\int_X\left\|\overline{g}\delta_{p-1}\Alt_p(f_1\otimes\cdots\otimes f_p)\right\|_{L^2(N_{p,x_0}\setminus D_{p,x_0},j_p(x_0,\cdot)}^2\mu(dx_0)\right)^{1/2}\notag\\
\leq \frac{1}{p+1}\sum_{k=0}^p\Big(\int_X\int_{N_{p,x_0}\setminus D_{p,x_0}}g(x_k)^2(\Alt_p(\delta_0 f_1(x_0,\cdot)\otimes\cdots\otimes \delta_0 f_p(x_0,\cdot))(x_1,...,x_p))^2\times\notag\\
\times j_p(x_0,d(x_1,...,x_p))\mu(dx_0)\Big)^{1/2};
\end{multline}
the estimate uses (\ref{E:determinant}) and the triangle inequality. Similarly as in Lemma \ref{L:elemL2} it follows that the summand for $k=0$ does not exceed 
\begin{equation}\label{E:x0est}
\left\|g\right\|_{L^2(X,\mu)}\prod_{\ell=1}^p \Big(\sup_{x\in X}\int_{N_{1,x}\setminus D_{1,x}}(f_\ell(y)-f_\ell(x))^2j(x,y)\mu(dy)\Big)^{1/2}.
\end{equation}
For fixed $k\neq 0$ the corresponding summand is bounded by
\begin{multline}
\Big(\int_X\int_{N_{p,x_0}\setminus D_{p,x_0}}g(x_0)^2(\Alt_p(\delta_0 f_1(x_0,\cdot)\otimes\cdots\otimes \delta_0 f_p(x_0,\cdot))(x_1,...,x_p))^2\times \notag\\
\times j_p(x_0,d(x_1,...,x_p))\mu(dx_0)\Big)^{1/2}\notag\\
+\Big(\int_X\int_{N_{p,x_0}\setminus D_{p,x_0}}(g(x_k)-g(x_0))^2(\Alt_p(\delta_0 f_1(x_0,\cdot)\otimes\cdots\otimes \delta_0 f_p(x_0,\cdot))(x_1,...,x_p))^2\times
\notag\\
\times j_p(x_0,d(x_1,...,x_p))\mu(dx_0)\Big)^{1/2}.\notag
\end{multline}
The first of these two terms admits the estimate (\ref{E:x0est}) as before. The square of the second is bounded by
\begin{multline}
\frac{1}{(p!)^2}\sum_{\pi\in \mathcal{S}_p}\sum_{\sigma\in \mathcal{S}_p}\int_X\int_{N_{p,x}\setminus D_{p,x}}(\delta_0g(x_0,x_k))^2\prod_{\ell=1}^p|\delta_0f_{\pi(\ell)}(x_0,x_\ell)|\prod_{m=1}^p|\delta_0f_{\sigma(m)}(x_0,x_m)|\times\notag\\
 \hspace{150pt}\times j_p(x_0,d(x_1,...,x_p))\mu(dx_0)\notag\\
\leq \frac{4}{(p!)^2}\sum_{\pi\in \mathcal{S}_p}\sum_{\sigma\in \mathcal{S}_p}\left\|f_{\pi(k)}\right\|_{\sup}\left\|f_{\sigma(k)}\right\|_{\sup}\int_X\Big(\int_{N_{1,x_0}\setminus D_{1,x_0}}\cdots \int_{N_{1,x_0}\setminus D_{1,x_0}}\prod_{\ell\neq k}|\delta_0f_{\pi(\ell)}(x_0,x_\ell)|\times\notag\\
\times \prod_{m\neq k}|\delta_0f_{\sigma(m)}(x_0,x_m)|\:j(x_0,x_1)\cdots \widehat{j(x_0,x_k)}\cdots j(x_0,x_p)\mu(dx_1)\cdots\widehat{\mu(dx_k)}\cdots \mu(dx_p)\Big)\times\notag\\
\times \int_{N_{1,x_0}\setminus D_{1,x_0}}(\delta_0g(x_0,x_k))^2j(x_0,x_k)\mu(dx_k)\mu(dx_0).
\end{multline}
Using Cauchy-Schwarz with respect to the variables $x_1,...,\hat{x}_k,...,x_p$ in the inner integral and Fubini, we can bound the preceding by
\begin{multline}
\frac{4}{(p!)^2}\sum_{\pi\in \mathcal{S}_p}\sum_{\sigma\in \mathcal{S}_p}\left\|f_{\pi(k)}\right\|_{\sup}\left\|f_{\sigma(k)}\right\|_{\sup}\int_X \prod_{\ell\neq k}\Big(\int_{N_{1,x_0}\setminus D_{1,x_0}}(\delta_0f_{\pi(\ell)}(x_0,x_\ell))^2j(x_0,x_\ell)\mu(dx_\ell)\Big)^{1/2}\times\\
\times \prod_{m\neq k}\Big(\int_{N_{1,x_0}\setminus D_{1,x_0}}(\delta_0f_{\sigma(m)}(x_0,x_m))^2j(x_0,x_m)\mu(dx_m)\Big)^{1/2}\times\\
\hspace{100pt}\times \int_{N_{1,x_0}\setminus D_{1,x_0}}(\delta_0g(x_0,x_k))^2j(x_0,x_k)\mu(dx_k)\mu(dx_0)\notag\\
\leq 4\:\prod_{\ell=1}^p \Big(\left\|f_\ell\right\|_{\sup}+\big(\sup_{x\in X}\int_{N_{1,x}\setminus D_{1,x}}(\delta_0f_\ell(x,y))^2j(x,y)\mu(dy)\big)^{1/2}\Big)^2\times\notag\\
\times \int_X\int_{N_{1,x}\setminus D_{1,x}}(\delta_0g(x,y))^2j(x,y)\mu(dy)\mu(dx).
\end{multline}
\end{proof}

For certain applications it can be convenient to restrict attention to smaller cores. We have the following approximation result for the norms induced by (\ref{E:Qp1}). 

\begin{corollary}\label{C:densesubalgebra}
Let Assumptions \ref{A:C}, \ref{A:CE} and \ref{A:abscont} be satisfied and let $N_\ast$ be a system of diagonal neighborhoods. Let $\mathcal{C}'$ be a subalgebra of $\mathcal{C}$, dense in $\mathcal{C}$ with respect to $\|\cdot\|_{\mathcal{D}(\delta_0,N_0)}$.
Then for any integer $p\geq 0$ the space $\mathcal{C}'^p(N_p)$ is dense in $\mathcal{C}^p(N_p)$ with respect to $\|\cdot\|_{\mathcal{D}(\delta_p,N_p)}$.
\end{corollary}

\begin{proof}
For $p=0$ the result is immediate. Suppose that $p\geq 1$ and let 
\[F=\sum_{i=1}^m\overline{g}^{(i)}\delta_{p-1}\Alt_p(f_1^{(i)}\otimes \dots \otimes f_p^{(i)})\] 
be an element of $\mathcal{C}^p(N_p)$. Let $(g_n^{(i)})_n$ and $(f_{j,n_j}^{(i)})_{n_j}$ be sequences from $\mathcal{C}'$ such that for all $i=1,...,m$ and $j=1,...,p$ we have 
$\lim_n \|g_n^{(i)}-g^{(i)}\|_{L^2(X,\mu)}=0$ and $\lim_{n_j}\|f_{j,n_j}^{(i)}-f_j^{(i)}\|_{\mathcal{D}(\delta_0,N_0)}=0$. Consider the elementary $p$-functions $F_{n,n_1,...,n_p}\in\mathcal{C}'^p(N_p)$ defined by
\[F_{n,n_1,...,n_p}:=\sum_{i=1}^m\overline{g}_n^{(i)}\delta_{p-1}\Alt_p(f_{1,n_1}^{(i)}\otimes \dots\otimes f_{p,n_p}^{(i)}).\]
We claim that 
\begin{equation}\label{E:oneaftertheother}
\lim_n\lim_{n_1}\dots\lim_{n_p}\big\|F-F_{n,n_1,...,n_p}\big\|_{\mathcal{D}(\delta_p,N_p)}=0.
\end{equation}
For fixed $i$ we have 
\begin{align}
\big\|\overline{g}^{(i)}\delta_{p-1}\Alt_p&(f_1^{(i)}\otimes \dots \otimes f_p^{(i)})-\overline{g}_n^{(i)}\delta_{p-1}\Alt_p(f_{1,n_1}^{(i)}\otimes \dots\otimes f_{p,n_p}^{(i)})\big\|_{L^2(N_p\setminus D_p,J_p)}\notag\\
&\leq \big\|(\overline{g}^{(i)}-\overline{g}_n^{(i)})\delta_{p-1}\Alt_p(f_1^{(i)}\otimes \dots \otimes f_p^{(i)})\big\|_{L^2(N_p\setminus D_p,J_p)}\notag\\
&\hspace{50pt} + \big\|\overline{g}_n^{(i)}\delta_{p-1}\Alt_p((f_{1}^{(i)}-f_{1,n_1}^{(i)})\otimes \dots \otimes f_p^{(i)})\big\|_{L^2(N_p\setminus D_p,J_p)}\notag\\
&\hspace{50pt} +\dots\notag\\
&\hspace{50pt} + \big\|\overline{g}_n^{(i)}\delta_{p-1}\Alt_p(f_{1,n_1}^{(i)}\otimes \dots \otimes (f_{p}^{(i)}-f_{p,n_p}^{(i)}))\big\|_{L^2(N_p\setminus D_p,J_p)}.\notag
\end{align}
The first summand goes to zero as $n\to \infty$ by Lemma \ref{L:L2dom}. For fixed $n$ the second summand goes to zero as $n_1\to \infty$ by Lemma \ref{L:elemL2} (ii), and so on. Similar estimates hold for
\[\big\|\delta_p\Alt_{p+1}(g^{(i)}\otimes f_1^{(i)}\otimes\dots\otimes f_p^{(i)})-\delta_p\Alt_{p+1}(g^{(i)}_n\otimes f_{1,n_1}^{(i)}\otimes\dots\otimes f_{p,n_p}^{(i)})\big\|_{L^2(N_{p+1}\setminus D_{p+1},J_{p+1})},\]
and using the triangle inequality we arrive at (\ref{E:oneaftertheother}).
\end{proof}

\subsection{Removable sets}

In this subsection we assume that $(X,\varrho)$ is locally compact and separable. 

The Hilbert complex $(\mathcal{D}(\delta_\ast,N_\ast),\delta_\ast)$ in (\ref{E:Hilbertcomplex}) was obtained from the elementary complex $(\mathcal{C}^\ast(N_\ast),\delta_\ast)$ in (\ref{E:elemcomplex}) by taking closures. Now suppose that $\Sigma\subset X$ is a closed set and let $\mathring{X}:=X\setminus \Sigma$. Let $\mathring{\mathcal{C}}$ denote the ideal in $\mathcal{C}$ consisting of all elements whose support is contained in the open set $\mathring{X}$ and let $\mathring{\mathcal{C}}^p(N_p)$ be defined as in (\ref{E:Vp}) but with $\mathring{\mathcal{C}}$ in place of $\mathcal{C}$.
Proposition \ref{P:elemcomplex} remains true and $(\mathring{\mathcal{C}}^\ast(N_\ast),\delta_\ast)$ 
is a subcomplex of $(\mathcal{C}^\ast(N_\ast),\delta_\ast)$. If $\Sigma$ is of $\mu$-measure zero, then $\mu^{p+1}$ in (\ref{E:Jpac}) can be replaced by $(\mu|_{\mathring{X}})^{p+1}$ without changing $J_p$, and taking the closures of the operators $(\delta_p,\mathring{\mathcal{C}}^p(N_p))$ in the spaces $L^2(N_p\setminus D_p,J_p)$ gives a Hilbert complex $(\mathring{\mathcal{D}}(\delta_\ast,N_\ast),\delta_\ast)$. Clearly it is a subcomplex of $(\mathcal{D}(\delta_\ast,N_\ast),\delta_\ast)$. One expects that for small enough $\Sigma$ the entire complexes coincide,
\begin{equation}\label{E:fullcoincidence}
(\mathring{\mathcal{D}}(\delta_\ast,N_\ast),\delta_\ast)=(\mathcal{D}(\delta_\ast,N_\ast),\delta_\ast),
\end{equation}
in the sense that for all integers $p\geq 0$ we have 
\begin{equation}\label{E:domaincoincidence}
\mathring{\mathcal{D}}(\delta_p,N_p)=\mathcal{D}(\delta_p,N_p).
\end{equation} 
For too large $\Sigma$ they should differ. This can be discussed in terms of capacities. 

Given $K\subset X$ compact, let
\begin{equation}\label{E:capscompact}
\cpct_{N_0}(K):=\inf\left\lbrace \|u\|_{\mathcal{D}(\delta_0,N_0)}^2: \ \text{$u\in\mathcal{C}$ and $u=1$ on a neighborhood of $K$}\right\rbrace
\end{equation}
with the convention that $\inf\emptyset:=+\infty$. For general $E\subset X$ we set
\begin{equation}\label{E:capsgen}
\cpct_{N_0}(E):=\sup\left\lbrace \cpct_{N_0}(K):\ \text{$K\subset E$ and $K$ compact}\right\rbrace.
\end{equation}

The following is related to \cite[Theorem 4.4]{BL92}. For the order $p=0$ the result is classical, \cite{Mazya85}, and by Corollary \ref{C:densesubalgebra} the order $p=0$ determines whether (\ref{E:fullcoincidence}) holds or not.

\begin{theorem}\label{T:removable}
Suppose that $\mathcal{C}\subset C_c(X)$, Assumptions \ref{A:C}, \ref{A:CE} and \ref{A:abscont} are satisfied and $\Sigma\subset X$ is closed.
Then $\cpct_{N_0}(\Sigma)=0$ if and only if $\mu(\Sigma)=0$ and (\ref{E:fullcoincidence}) holds.
\end{theorem}

The proof follows a standard pattern, we provide it for convenience.

\begin{proof} We write $\cpct:=\cpct_{N_0}$. Suppose that $\mu(\Sigma)=0$ and (\ref{E:fullcoincidence}) holds. By (\ref{E:capsgen}) we can find a sequence $(\Sigma_i)_i$ of compact sets $\Sigma_i\subset \Sigma$ such that 
\begin{equation}\label{E:capsup}
\cpct(\Sigma)=\sup_i \cpct(\Sigma_i).
\end{equation}
Let $i$ be fixed. By (\ref{E:domaincoincidence}) with $p=0$ any $u\in \mathcal{C}$ with $u=1$ on a neighborhood of $\Sigma_i$ can be approximated in $\|\cdot\|_{\mathcal{D}(\delta_0,N_0)}$ by a sequence $(u_n)_n\subset \mathring{\mathcal{C}}$, and this implies that $\cpct(\Sigma_i)\leq \lim_n \|u-u_n\|_{\mathcal{D}(\delta_0,N_0)}^2=0$. Now $\cpct(\Sigma)=0$ follows using (\ref{E:capsup}).

If $\cpct(\Sigma)=0$, then by the inner regularity of $\mu$, (\ref{E:capscompact}) and since $\|\cdot\|_{L^2(X,\mu)}\leq \|\cdot\|_{\mathcal{D}(\delta_0, N_0)}$ we must have $\mu(\Sigma)=0$. Clearly $\mathring{\mathcal{D}}(\delta_p,N_p)\subset \mathcal{D}(\delta_p,N_p)$. We claim that $\mathring{\mathcal{C}}$ is dense in $\mathcal{C}$ with respect to $\|\cdot\|_{\mathcal{D}(\delta_0, N_0)}$. If so, then Corollary \ref{C:densesubalgebra} implies that for any integer $p\geq 0$ the space $\mathring{\mathcal{C}}^p(N_p)$ is dense in $\mathcal{D}(\delta_p,N_p)$, and this proves (\ref{E:domaincoincidence}). To see the claim, let $f\in \mathcal{C}$. Then $\Sigma_f:=\Sigma\cap \supp f$ is compact and $\cpct(\Sigma_f)\leq \cpct(\Sigma)=0$. By (\ref{E:capscompact}) there is a sequence $(u_n)_n\subset \mathcal{C}$ such that $\lim_n\|u_n\|_{\mathcal{D}(\delta_0,N_0)}=0$ and $u_n=1$ on a neighborhood of $\Sigma_f$, respectively. The functions $f_n:=(1-u_n)f$ are in $\mathring{\mathcal{C}}$, and $\lim_n \|f-f_n\|_{\mathcal{D}(\delta_0,N_0)}=\lim_n\|f u_n\|_{\mathcal{D}(\delta_0,N_0)}=0$ since $\|u_n f\|_{L^2(X,\mu)}\leq \|f\|_{\sup}\|u_n\|_{L^2(X,\mu)}$ and 
\[\mathcal{Q}_0(u_n f)\leq \Big(\sup_{y\in X}\int_{N_{1,y}\setminus D_{1,y}}(f(y)-f(x))^2j(x,y)\mu(dx)\Big)\:\left\|u_n\right\|_{L^2(X,\mu)}^2+\left\|f\right\|_{\sup}^2\mathcal{Q}_0(u_n).\]
\end{proof}

In applications one can start from a non-local Dirichlet form as in Example \ref{Ex:nonlocalDf}. Recall that if  $(\mathcal{E},\mathcal{D}(\mathcal{E}))$ is a regular Dirichlet form on $L^2(X,\mu)$, \cite{FOT94}, then the associated capacity $\Dfcpct^\mathcal{E}$ is defined by 
\begin{equation}\label{E:Dfcap}
\Dfcpct^\mathcal{E}(A):=\inf\left\lbrace \left\|u\right\|_{\mathcal{D}(\mathcal{E})}^2:\ \text{$u\in \mathcal{D}(\mathcal{E})$ and $u\geq 1$ $\mu$-a.e. on A}\right\rbrace
\end{equation}
with $\inf\emptyset:=+\infty$ for $A\subset X$ open and by
\[\Dfcpct^\mathcal{E}(E):=\inf\left\lbrace \Dfcpct^\mathcal{E}(A):\ \text{$E\subset A$, $A$ open}\right\rbrace\]
for general $E\subset X$. Here $u\mapsto \left\|u\right\|_{\mathcal{D}(\mathcal{E})}:=(\mathcal{E}(u)+\left\|u\right\|_{L^2(X,\mu)}^2)^{1/2}$ denotes the norm in the Hilbert space $\mathcal{D}(\mathcal{E})$. See \cite[Section 2.1]{FOT94}. 

\begin{lemma}\label{L:quasilocal}
Let $(\mathcal{E},\mathcal{D}(\mathcal{E}))$ be a purely non-local regular Dirichlet form on $L^2(X,\mu)$ with representation
\begin{equation}\label{E:jumpformrep}
\mathcal{E}(f)=\int_X\int_X (f(x)-f(y))^2j(x,y)\mu(dx)\mu(dy)
\end{equation}
for all $f\in \mathcal{D}(\mathcal{E})$, where $j$ is a symmetric density. Assume that (\ref{E:kernelconditions}) holds and that $\mathcal{C}:=\lip_c(X)$ is dense in $\mathcal{D}(\mathcal{E})$.  Let $\varepsilon>0$ and let $N_\ast=N_\ast(\varepsilon)$ be as in (\ref{E:boundedrange}).
\begin{enumerate}
\item[(i)] We have $\mathcal{D}(\delta_0,N_0)=\mathcal{D}(\mathcal{E})$ and
\[\mathcal{Q}_0(f)\leq \mathcal{E}(f)\leq \mathcal{Q}_0(f)+4\left\|f\right\|_{L^2(X,\mu)}^2\sup_{x\in X}\int_{B(x,\varepsilon)^c}j(x,y)\mu(dy), \quad f\in \mathcal{D}(\mathcal{E}).\]
\item[(ii)] There is a constant $c(\varepsilon)>1$ such that for any Borel set $E\subset X$ we have
\begin{equation}\label{E:equivcap}
\cpct_{N_0(\varepsilon)}(E)\leq \Dfcpct^\mathcal{E}(E)\leq c(\varepsilon)\:\cpct_{N_0(\varepsilon)}(E).
\end{equation}
\end{enumerate}
\end{lemma}

Let 
\begin{equation}\label{E:T1}
T_1(s):=\min(\max(s,0),1),\quad s\in \mathbb{R}.
\end{equation}

\begin{remark}\label{R:goodbumps}
To $T_1$ one refers as the unit contraction, \cite{Allain75, BD58, BD59, FOT94, LeJan78}. Clearly $u\in \lip_c(X)$ implies $T_1(u)\in \lip_c(X)$.
\end{remark}


\begin{proof}
Item (i) follows from \cite[Proposition 4.1]{GrigHuLau14}. For (ii), note that since $\Dfcpct^\mathcal{E}$ is inner regular, \cite[Theorem 2.1.1]{FOT94}, it suffices to verify (\ref{E:equivcap}) for compact $K$. Let $K\subset X$ be compact and $\eta>0$. By (\ref{E:capscompact}) we can find $u\in \mathcal{C}$ such that $u=1$ on a neighborhood of $K$ and $\left\|u\right\|_{\mathcal{D}(\delta_0,N_0)}^2<\cpct_{N_0}(K)+\eta$.
By (i) there is a constant $c(\varepsilon)>0$ such that $\left\|u\right\|_{\mathcal{D}(\mathcal{E})}^2\leq c(\varepsilon)\:\left\|u\right\|_{\mathcal{D}(\delta_0,N_0)}^2$. Therefore $\Dfcpct^\mathcal{E}(K)\leq \left\|u\right\|_{\mathcal{D}(\mathcal{E})}^2\leq c(\varepsilon)\:(\cpct_{N_0}(K)+\eta)$,
and letting $\eta\to 0$ gives the inequality on the right hand side of (\ref{E:equivcap}) for $K$ in place of $E$. To see the inequality on the left hand side, let $K$ again be compact. Since $\mathcal{C}=\lip_c(X)$ is a special standard core, \cite[Lemma 2.2.7 (ii)]{FOT94} gives
\[\Dfcpct^\mathcal{E}(K)=\inf\{\left\|u\right\|_{\mathcal{D}(\mathcal{E})}^2:\ u\in \mathcal{C}\quad \text{and $u\geq 1$ on $K$}\},\]
and we may replace $\geq$ by $>$. For any $\eta>0$ we can therefore find some $u\in \mathcal{C}$ with $u>1$ on a neighborhood $U$ of $K$ and $\left\|u\right\|_{\mathcal{D}(\mathcal{E})}^2<\Dfcpct^\mathcal{E}(K)+\eta$. Clearly $T_1(u)=1$ on $U$, $T_1(u)\in\mathcal{C}$ and $\left\|T_1(u)\right\|_{\mathcal{D}(\mathcal{E})}\leq \left\|u\right\|_{\mathcal{D}(\mathcal{E})}$, and by (i) therefore 
\[\cpct_{N_0}(K)\leq \left\|T_1(u)\right\|_{\mathcal{D}(\delta_0,N_0)}^2\leq  \left\|T_1(u)\right\|_{\mathcal{D}(\mathcal{E})}^2<\Dfcpct^\mathcal{E}(K)+\eta,\]
Letting $\eta\to 0$ gives the left part of (\ref{E:equivcap}).
\end{proof}

\begin{remark}
The expression (\ref{E:jumpformrep}) makes sense for all $f\in L^2(X,\mu)$. The norm $\|\cdot\|_{\mathcal{D}(\mathcal{E})}$ extends to the maximal domain $\{f\in L^2(X,\mu): \mathcal{E}(f)<+\infty\}$ of $\mathcal{E}$ and makes it a Hilbert space.
\end{remark}

Recall that the $1$-resolvent of a Dirichlet form $(\mathcal{E},\mathcal{D}(\mathcal{E}))$ is the bounded linear operator $G_1^\mathcal{E}$ on $L^2(X,\mu)$ uniquely determined by 
\[\left\langle G_1^\mathcal{E}f,u\right\rangle_{\mathcal{D}(\mathcal{E})}=\left\langle f,u\right\rangle_{L^2(X,\mu)},\quad f\in L^2(X,\mu),\quad u\in \mathcal{D}(\mathcal{E}).\] 
We say that $G_1^\mathcal{E}$ has a pointwise defined density $g_1^\mathcal{E}$ (with respect to $\mu$) if $g_1^\mathcal{E}$ is a Borel function from $X^2$ to $(0,+\infty]$ and for all $f\in L^2(X,\mu)$ we have 
\[G_1^\mathcal{E}f(x)=\int_Xg_1^\mathcal{E}(x,y)f(y)\mu(dy)\quad \text{for $\mu$-a.e. $x\in X$}.\]
If metric estimates for $g_1^\mathcal{E}$ are available, the critical size of $\Sigma$ can be characterized in terms of Hausdorff measures. Given $s\geq 0$, we write $\mathcal{H}^s$ to denote the $s$-dimensional Hausdorff measure, and given a Hausdorff function $h$, we write $\mathcal{H}^h$ for the corresponding generalized Hausdorff measure, see \cite[Section 5.1]{AH96}, \cite[Section 2.10]{Federer} or \cite[Section 7.2.3]{Mazya85}.

\begin{corollary}\label{C:catalogue}
Assume that $X$ is complete. Let $\mu$ be a nonnegative Radon measure on $X$ with full support and $d$-regular as in (\ref{E:dregular}). Let $\alpha\in (0,2)$ and let $j$ be a symmetric kernel satisfying (\ref{E:frackerneld}). Let $\mathcal{E}$ be as in (\ref{E:jumpformrep}), and let $\mathcal{D}(\mathcal{E})$ be the closure of $\lip_c(X)$ in the maximal domain of $\mathcal{E}$. Then $(\mathcal{E},\mathcal{D}(\mathcal{E}))$ is a purely non-local regular Dirichlet form and its $1$-resolvent admits a pointwise defined density $g_1^\mathcal{E}(x,y)$. If in addition $\varepsilon>0$ and $N_\ast=N_\ast(\varepsilon)$ is as in (\ref{E:boundedrange}), then we observe the following for a closed set $\Sigma\subset X$ of measure zero: 
\begin{enumerate}
\item[(i)] For $d<\alpha$ the density $g_1^\mathcal{E}(x,y)$ is bounded on $X$. In this case (\ref{E:fullcoincidence}) holds if and only if $\Sigma=\emptyset$.
\item[(ii)] For $d>\alpha$ the density $g_1^\mathcal{E}(x,y)$ is bounded by $c\varrho(x,y)^{\alpha-d}$ on $X$ and lower bounded by $c^{-1}\varrho(x,y)^{\alpha-d}$ on $\{\varrho(x,y)<1\}$, where $c>1$ is a universal constant. In this case $\mathcal{H}^{d-\alpha}(\Sigma)<+\infty$ implies (\ref{E:fullcoincidence}); and if (\ref{E:fullcoincidence}) holds, we have $\mathcal{H}^{d-\alpha+\eta}(\Sigma)=0$ for any $\eta>0$. 
\item[(iii)] For $d=\alpha$ the density $g_1^\mathcal{E}(x,y)$ is comparable to $(1+(-\log \varrho(x,y))_+)$ on $\{\varrho(x,y)<1\}$. In this case $\mathcal{H}^{h}(\Sigma)<+\infty$ with $h=(1+(-\log \varrho(x,y))_+)^{-1}$ implies (\ref{E:fullcoincidence}); and if $X$ is bounded and (\ref{E:fullcoincidence}) holds, we have $\mathcal{H}^{k}(\Sigma)=0$ for any Hausdorff function $k$ satisfying $\int_0^1(-\log r)\:dk(r)<+\infty$.
\end{enumerate}
\end{corollary}

\begin{proof} By \cite[Proposition 2.2]{ChKu08} the Dirichlet form  $(\mathcal{E},\mathcal{D}(\mathcal{E}))$ is regular.
By \cite[Theorem 1.12]{GrigHuHu18} and the remarks following it $(\mathcal{E},\mathcal{D}(\mathcal{E}))$ admits a heat kernel that obeys typical two sided estimates. Standard calculations and well-known potential theoretic arguments then imply (i), (ii) and (iii), see \cite{AH96, FOT94, Mazya85}. 
\end{proof}

\begin{remark}\label{R:completeness}
Since $X$ is assumed to be complete, we may apply Frostman's lemma to go from zero capacity to zero Hausdorff measure in Corollary \ref{C:catalogue}. The other implication (finite Hausdorff measure giving zero capacity) does not need completeness.
\end{remark}

\begin{remark}\label{R:adapted}
Lemma \ref{L:quasilocal} and Corollary \ref{C:catalogue} can also be adapted to $\hat{N}_\ast=\hat{N}_\ast(\varepsilon)$ as in (\ref{E:boundedrangehat}).
\end{remark}

\begin{examples}
Corollary \ref{C:catalogue} applies in particular if $0<d\leq n$ and $X$ is a compact $d$-set in $\mathbb{R}^n$, \cite{ChKu03}. 
\end{examples}

Another class of examples can be obtained by subordination, \cite{Bertoin96, Garofalo19, Jacob01}. Suppose that $(P_t)_{t>0}$ is a strongly continuous Markov semigroup on $L^2(X,\mu)$, symmetric in the sense that 
\[\left\langle P_tf,g\right\rangle_{L^2(X,\mu)}=\left\langle f,P_tg\right\rangle_{L^2(X,\mu)}\] 
for all $t>0$ and $f,g\in L^2(X,\mu)$, \cite{BH91, FOT94}. We say that $(P_t)_{t>0}$ has a \emph{heat kernel} $p_t(x,y)$ (with respect to $\mu$) if $(t,x,y)\mapsto p_t(x,y)$ is a real valued function on $(0,+\infty)\times X^2$, Borel measurable in $(x,y)$ for each $t>0$ and such that for all $t>0$ and $f\in L^2(X,\mu)$ we have 
\[P_tf(x)=\int_X p_t(x,y) f(y)\mu(dy)\quad \text{for $\mu$-a.e. $x\in X$.}\] 

\begin{examples}\label{Ex:subordination}
Suppose that $\mu$ is $d$-regular and $(P_t)_{t>0}$ is a symmetric strongly continuous Markov semigroup on $L^2(X,\mu)$ with 
heat kernel  $p_t(x,y)$ admitting two-sided Gaussian estimates of the form 
\begin{equation}\label{E:Gaussian}
c_1^{-1}\:t^{-d/2}\:\exp\left\lbrace -c_2\frac{\varrho(x,y)^2}{t}\right\rbrace\leq p_t(x,y)\leq c_1\:t^{-d/2}\:\exp\left\lbrace -\frac{\varrho(x,y)^2}{c_2t}\right\rbrace,\quad t>0,\ x,y\in X;
\end{equation}
here $c_1>1$ and $c_2>1$ are fixed constants. Given $\alpha\in (0,2)$ let 
\[\nu_{\alpha/2}(dt)=\frac{\alpha\:dt}{2\mathbb{\Gamma}(1-\alpha/2)\:t^{\alpha/2+1}}\]
denote the L\'evy measure of the strictly $\alpha/2$-stable subordinator on $(0,\infty)$, see \cite[Section 3.1]{Bertoin96} or
\cite[Example 3.9.16]{Jacob01}; here $\mathbb{\Gamma}$ denotes the Euler Gamma function. The L\'evy jump density
\[j(x,y):=\int_{0+}^\infty p_t(x,dy)\nu_{\alpha/2}(dt),\quad x\in X,\]
satisfies (\ref{E:frackerneld}) and in particular, fits Assumption \ref{A:C} (i) and Assumption \ref{A:abscont}. Let $\mathcal{E}$ be the form having representation (\ref{E:jumpformrep}) with this density $j$ and endowed with the natural domain $\mathcal{D}(\mathcal{E})$ defined using $(P_t)_{t>0}$ and spectral theory. The heat kernel of $(\mathcal{E},\mathcal{D}(\mathcal{E}))$ then is 
\[p_t^\mathcal{E}(x,y)=\int_0^\infty p_s(x,y)q_t^{\alpha/2}(ds),\] 
where $q_t^{\alpha/2}$ denotes the law of the strictly $\alpha/2$-stable subordinator at time $t$, and the behaviour of $g_1^\mathcal{E}$ can be read off from 
\[g_1^\mathcal{E}(x,y)=\int_0^\infty e^{-t}p_t^\mathcal{E}(x,y)\:dt\] 
using the estimates (\ref{E:Gaussian}) for $p_t(x,y)$. If in addition $\lip_c(X)$ is dense in $\mathcal{D}(\mathcal{E})$, then Corollary \ref{C:catalogue} applies. This is, for instance, the case if $X=M$ is a complete smooth Riemannian manifold of dimension $d$ or $X$ is an $\mathrm{RCD}^\ast(0,d)$ space, \cite{AGS14, EKS15, Gigli17}.
\end{examples}

\section{Covers and cohomology}\label{S:covers}
Throughout this section $(X,\varrho)$ is a compact metric space and $\mathcal{C}$ is a subalgebra of $C(X)$. 

\subsection{Partitions of unity}
Let $I$ be a finite ordered set of cardinality $|I|$. Suppose that $\mathcal{U}=\{U_\alpha\}_{\alpha\in I}$ is an open cover of $X$ and $\mathcal{V}=\{V_\alpha\}_{\alpha\in I}$ is a shrinking of $\mathcal{U}$ in the sense that $\mathcal{V}$ is an open cover of $X$ and for any $\alpha\in I$ we have $\overline{V}_\alpha\subset U_\alpha$. Consider the system $N_\ast(\mathcal{V})=(N_p(\mathcal{V}))_{\alpha\in I}$ of diagonal neighborhoods 
\begin{equation}\label{E:comesfromcover2}
N_p(\mathcal{V}):=\bigcup_{\alpha \in I} V_\alpha^{p+1}
\end{equation}
as in (\ref{E:comesfromcover}) and suppose that $N_\ast=(N_p)_{p\geq 0}$ is a system of diagonal neighborhoods such that $N_\ast\prec N_\ast(\mathcal{V})$. For each integer $p\geq 0$ let 
\begin{equation}\label{E:cutcover}
U_\alpha^{(p+1)}:=U_\alpha^{p+1}\cap N_p; 
\end{equation}
note that since $N_0=X$, we have $U_\alpha^{(1)}=U_\alpha$, $\alpha \in I$. The family 
\[\mathcal{U}^{(p+1)}:=\{U_\alpha^{(p+1)}\}_{\alpha \in I}\]
is a finite open cover of $N_p$. 

We make an assumption on the existence of bump functions in $\mathcal{C}$. 
\begin{assumption}\label{A:Plateaus}
For any open set $U\subset X$ and any compact $K\subset U$ we can find some $f\in \mathcal{C}$ with $0\leq f\leq 1$, $\supp f\subset U$ and $f=1$ on $K$.
\end{assumption}

\begin{remark}\label{R:StoneWeierstrass}
By Stone-Weierstrass Assumption \ref{A:Plateaus} implies Assumption \ref{A:regularity}. If $\mathcal{C}$ is stable under composition with the unit contraction $T_1$ as defined in (\ref{E:T1}), then Assumption \ref{A:regularity} also implies Assumption \ref{A:Plateaus}.
\end{remark}

Under Assumption \ref{A:Plateaus} one can find related symmetric partitions of unity.
\begin{lemma}\label{L:partofunity}
Suppose that Assumption \ref{A:Plateaus} holds. Let $\mathcal{V}$, $\mathcal{U}$ and $N_\ast$ be as above. Then for each integer $p\geq 0$ there is a partition of unity $\{\chi_{p,\alpha}\}_{\alpha\in I}$ on $N_p$ subordinate to $\mathcal{U}^{(p+1)}$ and such that each $\chi_{p,\alpha}$ is a finite linear combination of functions of form $\chi^{\otimes(p+1)}$ with $\chi\in \mathcal{C}$ and, in particular, $\chi_{p,\alpha}\in \mathcal{C}^{s,p}$.
\end{lemma}

\begin{proof} We use a variant of standard arguments, \cite[Theorem 2.2]{Grigoryan2009}: Assumption \ref{A:Plateaus} ensures that for each $\alpha \in I$ we can find $\varphi_\alpha\in \mathcal{C}$ such that $0\leq \varphi_\alpha\leq 1$, $\supp\varphi_\alpha\subset U_\alpha$ and $\varphi_\alpha\equiv 1$ on $V_\alpha$. Given an integer $p\geq 0$ and $\alpha\in I$ we write
\[\varphi_{\alpha}^{\otimes (p+1)}:=\varphi_\alpha\otimes\cdots\otimes\varphi_\alpha\]
for the $(p+1)$-fold tensor product of $\varphi_\alpha$ with itself. Then
\begin{equation}\label{E:prodprops}
0\leq \varphi_{\alpha}^{\otimes (p+1)}\leq 1,\quad \supp\varphi_{\alpha}^{\otimes (p+1)}\subset U_\alpha^{p+1}\quad \text{and}\quad  \varphi_\alpha^{\otimes (p+1)}\equiv 1\quad \text{on}\quad V_\alpha^{p+1}.
\end{equation}  
Writing $I=\{\alpha_0,...,\alpha_{|I|-1}\}$, we consider the functions
\begin{multline}\chi_{p,\alpha_0}:=\varphi_{\alpha_0}^{\otimes (p+1)},\ \chi_{p,\alpha_1}:=\varphi_{\alpha_1}^{\otimes (p+1)}(\mathbf{1}^{\otimes (p+1)}-\varphi_{\alpha_0}^{\otimes (p+1)}),\ \dots, \notag\\
\chi_{p,\alpha_{|I|-1}}:=\varphi_{\alpha_{|I|-1}}^{\otimes (p+1)}(\mathbf{1}^{\otimes (p+1)}-\varphi_{\alpha_{|I|-2}}^{\otimes (p+1)})\cdots (\mathbf{1}^{\otimes (p+1)}-\varphi_{\alpha_0}^{\otimes (p+1)})
\end{multline}
on $X^{p+1}$; note that $\mathbf{1}^{\otimes (p+1)}\equiv 1$. All functions $\chi_{p,\alpha}$ are finite linear combinations of $(p+1)$-fold tensor powers of functions from $\mathcal{C}$, they satisfy $0\leq \chi_{p,\alpha}\leq 1$ and $\supp\chi_{p,\alpha}\subset U_{\alpha}^{p+1}$. We have 
\[1-\chi_{p,\alpha_0}-\dots -  \chi_{p,\alpha_{|I|-1}}=(\mathbf{1}^{\otimes (p+1)}-\varphi_{\alpha_0}^{\otimes (p+1)})\cdots (\mathbf{1}^{\otimes (p+1)}-\varphi_{\alpha_{|I|-1}}^{\otimes (p+1)}),\]
and by the last item in (\ref{E:prodprops}) the right hand side of this equality is zero. It follows that $\sum_{\alpha\in I}\chi_{p,\alpha}\equiv 1$. Seen as a function on $N_p$, each $\chi_{p,\alpha}$ has support in $U_{\alpha}^{(p+1)}$.
\end{proof}

\subsection{Mayer-Vietoris sequences}
As in the classical case of the deRham complex on a smooth manifold, \cite[Chapter II, Section 8]{BT82}, one can obtain generalized Mayer-Vietoris sequences. Let $\mathcal{V}$, $\mathcal{U}$ and $N_\ast$ be as above.

It is not difficult to see that for each $\alpha\in I$, the family $(U_\alpha^{(p+1)})_{p\geq 0}$ is a system of diagonal neighborhoods for $U_\alpha$; this is a slight variation of Examples \ref{Ex:N} (i). More generally, given an integer $q\geq 0$ and distinct $\alpha_0,...,\alpha_q\in I$, we write 
\begin{equation}\label{E:intersection}
U_{\alpha_0\cdots\alpha_q}^{(p+1)}:=\bigcap_{k=0}^q U_{\alpha_k}^{(p+1)};
\end{equation}
again $U_{\alpha_0\cdots\alpha_q}^{(1)}=U_{\alpha_0\cdots\alpha_q}$. Note that $U_{\alpha_0\cdots \alpha_q}^{(p+1)}=U_{\alpha_0\cdots \alpha_q}^{p+1}\cap N_{p}$. The family $(U_{\alpha_0\cdots\alpha_q}^{(p+1)})_{p\geq 0}$ is a system of diagonal neighborhoods for $U_{\alpha_0\cdots\alpha_q}$. We continue to use the shortcut notation (\ref{E:defasrestriction}), that is, $\mathcal{C}^p(U_{\alpha_0\cdots\alpha_q}^{(p+1)}):=\mathcal{C}^p|_{U_{\alpha_0\cdots\alpha_q}^{(p+1)}}$. By Proposition \ref{P:elemcomplex} the sequence 
\begin{equation}\label{E:elemcomplexU}
0\longrightarrow \mathcal{C}^0(U_{\alpha_0\cdots\alpha_q})\stackrel{\delta_0}{\longrightarrow}\mathcal{C}^1(U_{\alpha_0\cdots\alpha_q}^{(2)})\stackrel{\delta_1}{\longrightarrow} ... \stackrel{\delta_{p-1}}{\longrightarrow} \mathcal{C}^p(U_{\alpha_0\cdots\alpha_q}^{(p+1)})  \stackrel{\delta_p}{\longrightarrow} ...
\end{equation}
is a cochain complex. 

Now suppose that also Assumptions \ref{A:C}, \ref{A:CE} and \ref{A:abscont} hold. For any integer $p\geq 0$
let $\mathcal{D}(\delta_p,U_{\alpha_0\cdots\alpha_q}^{(p+1)})$ be the domain of the closure of $(\delta_p, \mathcal{C}^p(U_{\alpha_0\cdots\alpha_q}^{(p+1)}))$ in $L^2(U_{\alpha_0\cdots\alpha_q}^{(p+1)}\setminus D_p,J_p)$ as constructed in Theorem \ref{T:closable}; the closability follows similarly as there. Then also the sequence 
\begin{equation}\label{E:HilbertcomplexU}
0\longrightarrow \mathcal{D}(\delta_0,U_{\alpha_0\dots \alpha_q})\stackrel{\delta_0}{\longrightarrow}\mathcal{D}(\delta_1,U_{\alpha_0\dots \alpha_q}^{(2)})\stackrel{\delta_1}{\longrightarrow} ... \stackrel{\delta_{p-1}}{\longrightarrow} \mathcal{D}(\delta_p,U_{\alpha_0\dots \alpha_q}^{(p+1)})  \stackrel{\delta_p}{\longrightarrow} ...
\end{equation}
is a cochain complex.

Let $r$ be the natural restriction that takes an element $F$ of $\mathcal{C}^p(N_p)$ into the element 
$rF=((rF)_\alpha)_{\alpha\in I}$ of the direct product $\prod_{\alpha\in I} \mathcal{C}^p(U_\alpha^{(p+1)})$ with components 
\[(rF)_\alpha=F|_{U_\alpha^{(p+1)}};\]
note that $\mathcal{C}^p(N_p)|_{U_{\alpha}^{(p+1)}}\subset \mathcal{C}^p(U_{\alpha}^{(p+1)})$ by (\ref{E:defasrestriction}).
Assume that $<$ is a strict order on $I$. Given an ordered subset $\alpha_0<\dots<\alpha_q$ of $I$ and an element $F$ of the direct product $\prod_{\alpha_0<\dots<\alpha_q}\mathcal{C}^p(U_{\alpha_0\cdots \alpha_q}^{(p+1)})$ with components $F_{\alpha_0\dots\alpha_q}\in \mathcal{C}^p(U_{\alpha_0\cdots \alpha_q}^{(p+1)})$, we define an element $\check{\delta}F=((\check{\delta}F)_{\alpha_0\dots \alpha_{q+1}}))_{\alpha_0<\dots<\alpha_{q+1}}$ of $\prod_{\alpha_0<\dots<\alpha_{q+1}}\mathcal{C}^p(U_{\alpha_0\cdots \alpha_{q+1}}^{(p+1)})$ component-wise by
\[(\check{\delta}F)_{\alpha_0\dots \alpha_{q+1}}:=\sum_{i=0}^{q+1}(-1)^i F_{\alpha_0\dots \hat{\alpha}_i\dots\alpha_{q+1}}.\]
It is quickly seen that $\check{\delta}\circ\check{\delta}=0$. Under Assumptions \ref{A:C}, \ref{A:CE} and \ref{A:abscont} we can apply the operators $r$ and $\check{\delta}$ similarly to the spaces $\mathcal{D}(\delta_p,N_p)$ and $\prod_{\alpha_0<\dots<\alpha_q}\mathcal{D}(\delta_p,U_{\alpha_0...\alpha_q}^{(p+1)})$. Note that by (\ref{E:defasrestriction}) and Theorem \ref{T:closable} we have $\mathcal{D}(\delta_p,N_p)|_{U_{\alpha}^{(p+1)}}\subset \mathcal{D}(\delta_p,U_{\alpha}^{(p+1)})$ and so on.

With the convention that the interchange of two indices provokes a change of sign, 
\[F_{\alpha_0\cdots \alpha_k\cdots \alpha_i\cdots \alpha_q}:=- F_{\alpha_0\cdots \alpha_i\cdots \alpha_k\cdots \alpha_q},\] 
we may drop the requirement that $\alpha_0<\dots<\alpha_q$, see \cite[Exercise 8.4]{BT82}.

\begin{proposition}\label{P:MVlemma}
Let Assumption \ref{A:Plateaus} be in force. Let $\mathcal{V}$, $\mathcal{U}$ and $N_\ast$ be as above.
\begin{enumerate}
\item[(i)] For any integer $p\geq 0$ the sequence 
\begin{equation}\label{E:MVlemma}
0\longrightarrow\mathcal{C}^p(N_p)\stackrel{r}{\longrightarrow}\prod_{\alpha_0} \mathcal{C}^p(U_{\alpha_0}^{(p+1)})\stackrel{\check{\delta}}{\longrightarrow}\prod_{\alpha_0<\alpha_1} \mathcal{C}^p(U_{\alpha_0\alpha_1}^{(p+1)})\stackrel{\check{\delta}}{\longrightarrow} ... 
\end{equation}
is exact.
\item[(ii)] Suppose that also Assumptions \ref{A:C}, \ref{A:CE} and \ref{A:abscont} hold.
Then for any integer $p\geq 0$ the sequence 
\begin{equation}\label{E:MVlemmaclosures}
0\longrightarrow\mathcal{D}(\delta_p,N_p)\stackrel{r}{\longrightarrow}\prod_{\alpha_0} \mathcal{D}(\delta_p,U_{\alpha_0}^{(p+1)})\stackrel{\check{\delta}}{\longrightarrow}\prod_{\alpha_0<\alpha_1} \mathcal{D}(\delta_p,U_{\alpha_0\alpha_1}^{(p+1)})\stackrel{\check{\delta}}{\longrightarrow} ... 
\end{equation}
is exact. 
\end{enumerate}
\end{proposition}

For convenience we provide a suitable variant of the well-known classical arguments.

\begin{proof}
The injectivity of $r$ in (\ref{E:MVlemma}) and (\ref{E:MVlemmaclosures}) is clear. Let $F=(F_\alpha)_{\alpha\in I}\in \prod_\alpha \mathcal{C}^p(U_\alpha^{(p+1)})$ be such that $\check{\delta}F=0$. Then $F_{\alpha_0}=F_{\alpha}$ on $U_{\alpha\alpha_0}^{(p+1)}$ for any distinct $\alpha,\alpha_0 \in I$. Let $\chi_{p,\alpha}$, $\alpha\in I$, be as in Lemma \ref{L:partofunity} and set $G:=\sum_{\alpha\in I}\chi_{p,\alpha} F_\alpha$. Since each $F_\alpha$ is a restriction to the corresponding open set $U_\alpha^{(p+1)}$ of an element of $\mathcal{C}^p$ and $\chi_{p,\alpha}$ is in $\mathcal{C}^{s,p}$ and supported in $U_\alpha^{(p+1)}$, we have $\chi_{p,\alpha} F_\alpha\in \mathcal{C}^p(N_p)$ by (\ref{E:symmaltmult}) and the comments following it. Consequently $G\in \mathcal{C}^p(N_p)$. It is quickly seen that $rG=F$, and we can conclude that $\im\check{\delta}=\ker\check{\delta}$ for $\check{\delta}$ acting along the third arrow in (\ref{E:MVlemma}). Under the hypotheses in (ii) the same argument applies to $F=(F_\alpha)_{\alpha\in I}\in \prod_\alpha \mathcal{D}(\delta_p,U_\alpha^{(p+1)})$, because for each $\alpha$ we have $\chi_{p,\alpha}F_\alpha \in \mathcal{D}(\delta_p,N_p)$ by Corollary \ref{C:Sobomult}. If $q\geq 1$ and $F=(F_{\alpha_0\cdots\alpha_q})_{\alpha_0<...<\alpha_q}$ is an element of $\prod_{\alpha_0<...<\alpha_q}\mathcal{C}^p(U_{\alpha_0...\alpha_q}^{(p+1)})$ satisfying $\check{\delta}F=0$, then  
\[F_{\alpha_0\cdots \alpha_q}=\sum_{i=0}^q(-1)^iF_{\alpha\alpha_0\cdots\hat{\alpha}_i\cdots\alpha_q}\quad \text{on $U_{\alpha\alpha_0\cdots\alpha_q}^{(p+1)}$}\]
for any distinct $\alpha, \alpha_0,...\alpha_q\in I$. Let $G=(G_{\alpha_0\cdots\alpha_{q-1}})_{\alpha_0<...<\alpha_{q-1}}$ be defined by  
\[G_{\alpha_0\cdots \alpha_{q-1}}:=\sum_{\alpha\in I}\chi_{p,\alpha}F_{\alpha\alpha_0\cdots \alpha_{q-1}}.\]
Again by (\ref{E:symmaltmult}) this defines an element $G$ of $\prod_{\alpha_0<...<\alpha_{q-1}}\mathcal{C}^p(U_{\alpha_0...\alpha_{q-1}}^{(p+1)})$. Since 
\[(\check{\delta}G)_{\alpha_0\cdots\alpha_q}=\sum_{i=0}^q(-1)^i G_{\alpha_0\cdots \hat{\alpha}_i\cdots \alpha_q}=\sum_{\alpha\in I}\chi_{p,\alpha}\sum_{i=0}^q(-1)^i F_{\alpha\alpha_0\cdots \hat{\alpha}_i\cdots\alpha_q}=F_{\alpha_0\cdots\alpha_q}\]
we see that  $\check{\delta}G=F$. This shows $\im\check{\delta}=\ker\check{\delta}$ along the remaining arrows.
As before the arguments remain valid for $F\in \prod_{\alpha_0<...<\alpha_q}\mathcal{D}(\delta_p,U_{\alpha_0...\alpha_q}^{(p+1)})$ if the hypotheses of (ii) are assumed; note that in this case $\chi_{p,\alpha} F_{\alpha\alpha_0\cdots \alpha_{q-1}}\in \mathcal{D}(\delta_p,U_{\alpha_0\cdots\alpha_{q-1}}^{(p+1)})$ by Corollary \ref{C:Sobomult} and therefore $G\in \prod_{\alpha_0<...<\alpha_{q-1}}\mathcal{D}(\delta_p,U_{\alpha_0\cdots\alpha_{q-1}}^{(p+1)})$.
\end{proof}

The complexes (\ref{E:elemcomplexU}) and (\ref{E:MVlemma}), together with the operators $\delta_\ast$ and $r$ respectively $\check{\delta}$, define an augmented \emph{\v{C}ech-Alexander type} bicomplex of the form
\begin{equation}\label{E:bicomplex}
	\begin{tikzcd}[column sep = scriptsize]
		& \vdots & \vdots & \vdots & \\
		0 \ar[r] & {\mathcal{C}^2(N_2)} \ar[u, "\delta_2"]\ar[r, "r"] & {\smash{\prod\limits_{\alpha_0}}\:\mathcal{C}^2(U_{\alpha_0}^{(3)})} \ar[u, "\delta_2"]\ar[r, "\check{\delta}"] & {\smash{\prod\limits_{\alpha_0 <\alpha_1}}\mathcal{C}^2(U^{(3)}_{\alpha_0\alpha_1})} \ar[u, "\delta_2"]\ar[r, "\check{\delta}"] & \dotsb \\
		0 \ar[r] & {\mathcal{C}^1(N_1)} \ar[u, "\delta_1"]\ar[r, "r"] & {\smash{\prod\limits_{\alpha_0}}\:\mathcal{C}^1(U_{\alpha_0}^{(2)})} \ar[u, "\delta_1"]\ar[r, "\check{\delta}"] & {\smash{\prod\limits_{\alpha_0 <\alpha_1}}\mathcal{C}^1(U_{\alpha_0\alpha_1}^{(2)})} \ar[u, "\delta_1"]\ar[r, "\check{\delta}"] & \dotsb \\		
		0 \ar[r] & {\mathcal{C}^0(N_0)} \ar[u, "\delta_0"]\ar[r, "r"] & {\smash{\prod\limits_{\alpha_0}}\:\mathcal{C}^0(U_{\alpha_0})} \ar[u, "\delta_0"]\ar[r, "\check{\delta}"] & {\smash{\prod\limits_{\alpha_0 <\alpha_1}}\mathcal{C}^0(U_{\alpha_0\alpha_1})} \ar[u, "\delta_0"]\ar[r, "\check{\delta}"] & \dotsb 
	\end{tikzcd}
\end{equation}

Given integers $p,q\geq 0$, let $K^{q, p}_\mathcal{C}(\mathcal{U},N_\ast) := \prod_{\alpha_0< \dots <\alpha_q} \mathcal{C}^p(U_{\alpha_0\dotsb\alpha_q}^{(p+1)})$. Defining 
\[\mathcal{K}_\mathcal{C}^\ell(\mathcal{U},N_\ast):=\bigoplus_{p+q=\ell} K^{q,p}_\mathcal{C}(\mathcal{U},N_\ast)\quad \text{for any $\ell\geq 0$ and}\quad D_\ell:=\check{\delta}+(-1)^q\delta_p\quad \text{on $K_\mathcal{C}^{q,p}(\mathcal{U},N_\ast)$},\] 
we obtain a cochain complex $(\mathcal{K}_\mathcal{C}^\ast(\mathcal{U},N_\ast),D_\ast)$.

Analogously the complexes (\ref{E:HilbertcomplexU}) and (\ref{E:MVlemmaclosures}) define a bicomplex with $\mathcal{D}(\delta_p,U_{\alpha_0\dotsb\alpha_q}^{(p+1)})$ in place of $\mathcal{C}^p(U_{\alpha_0\dotsb\alpha_q}^{(p+1)})$. We write $K^{q, p}_\mathcal{D}(\mathcal{U},N_\ast)$ and $\mathcal{K}_\mathcal{D}^\ell(\mathcal{U},N_\ast)$ for the corresponding counterparts of the above spaces; this gives a  complex $(\mathcal{K}_\mathcal{D}^\ast(\mathcal{U},N_\ast),D_\ast)$.

We write 
\[H^\ell\mathcal{K}_\mathcal{C}^\ast(\mathcal{U},N_\ast):=\ker D_\ell|_{\mathcal{K}_\mathcal{C}^\ell(\mathcal{U},N_\ast)}/\im D_{\ell-1}|_{\mathcal{K}_\mathcal{C}^{\ell-1}(\mathcal{U},N_\ast)}\] 
and $H^\ell\mathcal{K}_\mathcal{D}^\ast(\mathcal{U},N_\ast)$, defined similarly, to denote the respective $\ell$-th cohomologies. 

Recall (\ref{E:elemcoho}) and (\ref{E:closedcoho}). By Proposition \ref{P:MVlemma} and well-known abstract arguments, \cite[Proposition 8.8]{BT82}, \cite[Lemma 4]{BSSS12}, one can see that $r$ is a cochain map from $\mathcal{C}^\ast(N_\ast)$ to $\mathcal{K}_\mathcal{C}^\ast(\mathcal{U},N_\ast)$ respectively from $\mathcal{D}(\delta_\ast, N_\ast)$ to $\mathcal{K}_\mathcal{D}^\ast(\mathcal{U},N_\ast)$ and induces an isomorphism in cohomology.  

\begin{corollary}\label{C:MVlemma}
Let Assumption \ref{A:Plateaus} be in force. Let $\mathcal{V}$, $\mathcal{U}$ and $N_\ast$ be as above.
\begin{enumerate}
\item[(i)] For any integer $\ell\geq 0$ the spaces $H^\ell\mathcal{C}^\ast(N_\ast)$ and $H^\ell\mathcal{K}_\mathcal{C}^\ast(\mathcal{U},N_\ast)$ are isomorphic.
\item[(ii)] Suppose that also Assumptions \ref{A:C}, \ref{A:CE} and \ref{A:abscont} hold.
Then for any integer $\ell\geq 0$ the spaces $H^\ell\mathcal{D}(\delta_\ast, N_\ast)$ and $H^\ell\mathcal{K}_\mathcal{D}^\ast(\mathcal{U},N_\ast)$ are isomorphic.
\end{enumerate}
\end{corollary}

\subsection{Poincar\'e lemma}

We prove versions of Poincar\'e's lemma for specific finite covers by open balls. We start with a small observation.
\begin{lemma}\label{L:smallobs}
Given $f_0,...,f_p\in \mathcal{C}\oplus\mathbb{R}$, we have
\begin{equation}\label{E:smallobs}
\Alt_{p+1}(f_0\otimes\cdots\otimes f_p)=\frac{1}{p+1}\sum_{k=0}^p(-1)^k f_k\otimes \Alt_p(f_0\otimes \cdots \otimes \hat{f}_k\otimes \cdots \otimes f_p).
\end{equation}
\end{lemma}

\begin{proof}
The evaluation at $(x_0,...,x_p)$ of the left hand side in (\ref{E:smallobs}) equals
\begin{align}
\frac{1}{(p+1)!}&\sum_{\sigma\in \mathcal{S}_{p+1}}\sgn \sigma \:f_{\sigma(0)}(x_0)f_{\sigma(1)}(x_1)\cdots f_{\sigma(p)}(x_p)\notag\\
&=\frac{1}{(p+1)!}\Big\lbrace \sum_{\sigma\in \mathcal{S}_{p+1}:\ \sigma(0)=0}\sgn \sigma\:f_0(x_0)f_{\sigma(1)}(x_1)\cdots f_{\sigma(p)}(x_p)\notag\\
&\hspace{50pt} + \sum_{\sigma\in \mathcal{S}_{p+1}:\ \sigma(0)=1}\sgn \sigma\:f_1(x_0)f_{\sigma(1)}(x_1)\cdots f_{\sigma(p)}(x_p)\notag\\
&\hspace{50pt} + \dots\notag\\
&\hspace{50pt} + \sum_{\sigma\in \mathcal{S}_{p+1}:\ \sigma(0)=p}\sgn \sigma\:f_p(x_0)f_{\sigma(1)}(x_1)\cdots f_{\sigma(p)}(x_p)\Big\rbrace.\notag
\end{align}
Varying $\sigma$ in the first summand, $(\sigma(1),...,\sigma(p))$ runs through all permutations of the tupel $(1,...,p)$, and if 
$(\sigma(1),...,\sigma(p))=(1, ...,p)$, then clearly $\sgn\sigma=1$. Varying $\sigma$ in the second summand, $(\sigma(1),...,\sigma(p))$ runs through all permutations of $(0,2,...,p)$, and if $(\sigma(1),...,\sigma(p))=(0,2, ...,p)$ (that is, if the natural order is preserved), we have $\sgn\sigma=-1$. The other summands behave similarly. Therefore the preceding equals
\begin{align}
&\frac{1}{(p+1)!}\Big\lbrace \sum_{\pi\in \mathcal{S}_{p}(1,...,p)}\sgn \pi\:f_0(x_0)f_{\pi(1)}(x_1)\cdots f_{\pi(p)}(x_p)\notag\\
&\hspace{50pt} - \sum_{\pi\in \mathcal{S}_{p}(0,2,...,p)}\sgn \pi\:f_1(x_0)f_{\pi(0)}(x_1)\cdots f_{\pi(p)}(x_p)\notag\\
&\hspace{50pt} + \dots\notag\\
&\hspace{50pt} +(-1)^p \sum_{\pi\in \mathcal{S}_{p}(0,...,p-1)}\sgn \pi\:f_p(x_0)f_{\pi(0)}(x_1)\cdots f_{\pi(p-1)}(x_p)\notag\Big\rbrace\\
&=\frac{1}{p+1}\sum_{k=0}^p(-1)^k f_k(x_0)\Alt_p(f_0\otimes \cdots \otimes\hat{f}_k\otimes \cdots \otimes f_p)(x_1,...,x_p);\notag
\end{align}
here $\mathcal{S}_{p}(1,...,p)$ means that the permutations are applied to $(1,...,p)$, and so on.
\end{proof}

Now let $\varepsilon>0$ and let $\hat{N}_\ast(\varepsilon)=(\hat{N}_p(\varepsilon))_{p\geq 0}$ be as in (\ref{E:boundedrangehat}). We follow \cite[Section 9]{BSSS12} and, given $p\geq 1$ and $(x_0,...,x_{p-1})\in \hat{N}_{p-1}(\varepsilon)$, define the \emph{slice} 
\[\hat{S}_p(\varepsilon;(x_0,...,x_{p-1})):=\{t\in X:\ (t,x_0,...,x_{p-1})\in \hat{N}_{p}(\varepsilon)\}.\]
The following assumption is a variant of \cite[Hypothesis ($\ast$), p. 34]{BSSS12}. It is a condition on $X$, $\varepsilon$ and an integer $K\geq 0$. It guarantees that for a suitable finite cover of $X$ by open balls we can prove a Poincar\'e lemma up to order $K$, compatible with Proposition \ref{P:MVlemma}. 

\begin{assumption}\label{A:hypostar}
There is a number $\eta>0$ such that for any nonempty intersection $B_{\alpha_0\cdots \alpha_q}=\bigcap_{k=0}^q B_{\alpha_k}$ of finitely many open balls $B_{\alpha_k}$ of radius $\varepsilon+2\eta$ we can find a nonempty open set $W_{\alpha_0\cdots \alpha_q}$ satisfying 
\begin{equation}\label{E:setW}
W_{\alpha_0\cdots \alpha_q}\subset B_{\alpha_0\cdots \alpha_q}\cap \Big(\bigcap_{(x_0,...,x_{p-1})\in B_{\alpha_0\cdots \alpha_q}^{(p)}} \hat{S}(\varepsilon;(x_0,...,x_{p-1}))\Big)\quad \text{for any $p\leq K+1$}.
\end{equation}
Here $B_{\alpha_0\cdots \alpha_q}^{(p)}$ is defined as in (\ref{E:intersection}). 
\end{assumption}

Suppose that $X$, $\varepsilon$ and $K$ satisfy Assumption \ref{A:hypostar} and let $\eta$ be as stated there. Let $\{B(y_\alpha,\eta)\}_{\alpha\in I}$ be a finite cover of $X$ by open balls $B(y_\alpha,\eta)$ having centers $y_\alpha\in X$ and common radius $\eta$.
Since $X$ is assumed to be compact, we can always find such a cover. Then clearly also $\mathcal{V}:=\{B(y_\alpha,\varepsilon+\eta)\}_{\alpha\in I}$ covers $X$, and it is easily seen that the system $\hat{N}_\ast(\varepsilon)$ as defined in (\ref{E:boundedrangehat}) is dominated by the system $N_\ast(\mathcal{V})$ as defined in (\ref{E:comesfromcover}) respectively (\ref{E:comesfromcover2}) in the sense that 
\begin{equation}\label{E:systemdominates}
\hat{N}_\ast(\varepsilon)\prec N_\ast(\mathcal{V}).
\end{equation}

We consider the finite open cover 
\begin{equation}\label{E:Uviaepseta}
\mathcal{U}:=\{B_\alpha\}_{\alpha\in I}\quad \text{consisting of the open balls $B_\alpha:=B(y_\alpha,\varepsilon+2\eta)$, $\alpha\in I$.}
\end{equation}
Here we add one more $\eta$ to enlarge the radius so that each $B_\alpha$ contains the closure of the concentric smaller ball from $\mathcal{V}$. Note that Lemma \ref{L:partofunity} and Proposition \ref{P:MVlemma} may be applied with these ball covers $\mathcal{V}$, $\mathcal{U}$ and with $\hat{N}_\ast(\varepsilon)$ in place of $N_\ast$.

Suppose that $\mu$ is a finite nonnegative Borel measure on $X$ with full support and that Assumption \ref{A:hypostar} is in force. Given $p,q\geq 0$ and a bounded Borel function $F:B_{\alpha_0\cdots \alpha_q}^{(p+1)}\to \mathbb{R}$, where $B_{\alpha_0\cdots \alpha_q}$ is a nonempty intersection, we can define a bounded Borel function 
\[\Psi_{B_{\alpha_0\cdots \alpha_q}} F:B_{\alpha_0\cdots \alpha_q}^{(p)}\to \mathbb{R}\] 
by 
\begin{equation}\label{E:PsiU}
\Psi_{B_{\alpha_0\cdots \alpha_q}}F(x_0,...,x_{p-1}):=\frac{1}{\mu(W_{\alpha_0\cdots \alpha_q})}\int_{W_{\alpha_0\cdots \alpha_q}}F(t,x_0,...,x_{p-1})\:\mu(dt),\ (x_0,...,x_{p-1})\in B_{\alpha_0\cdots \alpha_q}^{(p)},
\end{equation}
where $W_{\alpha_0\cdots \alpha_q}$ is as in (\ref{E:setW}). For brevity we suppress $p$ from notation, although $\Psi_{B_{\alpha_0\cdots \alpha_q}}$ clearly depends on it. Note that $\Psi_{B_{\alpha_0\cdots \alpha_q}}$ is well-defined, since by (\ref{E:setW}) we have 
\begin{equation}\label{E:productiscontained}
W_{\alpha_0\cdots \alpha_q}\times B_{\alpha_0\cdots \alpha_q}^{(p)}\subset B_{\alpha_0\cdots \alpha_q}^{(p+1)}.
\end{equation} 
If instead $F$ is a $\mu^{p+1}$-class of functions on $B_{\alpha_0\cdots \alpha_q}^{(p+1)}$, then -- subject to suitable integrability conditions -- also $\Psi_{B_{\alpha_0\cdots \alpha_q}}F$ is well-defined as a $\mu^{p}$-class of functions on $B_{\alpha_0\cdots \alpha_q}^{(p)}$.

\begin{proposition}\label{P:Poincare}
Let $\mu$ be a finite nonnegative Borel measure on $X$ with full support. Let $\varepsilon>0$ and $K\geq 0$ be such that Assumption \ref{A:hypostar} holds, let $\eta>0$ be as there and $\mathcal{U}$ as in (\ref{E:Uviaepseta}). Suppose that $B_{\alpha_0\cdots \alpha_q}$ is a nonempty intersection.
\begin{enumerate}
\item[(i)] For any $1\leq p\leq K+1$ identity (\ref{E:PsiU}) defines a linear map 
\[\Psi_{B_{\alpha_0\cdots \alpha_q}}:\mathcal{C}^p(B_{\alpha_0\cdots \alpha_q}^{(p+1)})\to \mathcal{C}^{p-1}(B_{\alpha_0\cdots \alpha_q}^{(p)}).\] 
Moreover, for any $F\in \mathcal{C}^p(B_{\alpha_0\cdots \alpha_q}^{(p+1)})$ we have
\begin{equation}\label{E:homotopy}
(\Psi_{B_{\alpha_0\cdots \alpha_q}}\circ \delta_p+\delta_{p-1}\circ \Psi_{B_{\alpha_0\cdots \alpha_q}})F=F.
\end{equation}
\item[(ii)] Suppose in addition that Assumptions \ref{A:C}, \ref{A:CE} and \ref{A:abscont} are satisfied and that
\begin{equation}\label{E:jlowerbound}
\inf_{(x,y)\in \hat{N}_1(\varepsilon)\setminus D_1}j(x,y)\geq c
\end{equation}
with a constant $c>0$. Then for any $1\leq p\leq K+1$ identity (\ref{E:PsiU}) defines a bounded linear map 
\[\Psi_{B_{\alpha_0\cdots \alpha_q}}:L^2_a(B_{\alpha_0\cdots \alpha_q}^{(p+1)}\setminus D_p,J_p)\to L^2_a(B_{\alpha_0\cdots \alpha_q}^{(p)}\setminus D_{p-1},J_{p-1})\] 
and there is a constant $c>0$ such that 
\begin{equation}\label{E:PsiUL2est}
\left\|\Psi_{B_{\alpha_0\cdots \alpha_q}}F\right\|_{L^2(B_{\alpha_0\cdots \alpha_q}^{(p)}\setminus D_{p-1},J_{p-1})}\leq c\:\left\|F\right\|_{L^2(B_{\alpha_0\cdots \alpha_q}^{(p+1)}\setminus D_{p},J_{p})},\quad F\in L^2_a(B_{\alpha_0\cdots \alpha_q}^{(p+1)}\setminus D_p,J_p).
\end{equation}
Moreover, for any $F\in \mathcal{D}(\delta_p,B_{\alpha_0\cdots \alpha_q}^{(p+1)})$ we have $\Psi_{B_{\alpha_0\cdots \alpha_q}}F\in \mathcal{D}(\delta_p,B_{\alpha_0\cdots \alpha_q}^{(p)})$ and (\ref{E:homotopy}) holds on $B_{\alpha_0\cdots \alpha_q}^{(p+1)}$ in the $J_p$-a.e. sense.
\end{enumerate}
\end{proposition}

\begin{proof}
By Assumption \ref{A:hypostar} the map $\Psi_{B_{\alpha_0\cdots \alpha_q}}$ in (i) is well-defined. Since it is obviously linear, we may assume that $F=\Alt_{p+1}(f_0\otimes...\otimes f_p)$ with $f_0\in \mathcal{C}\oplus \mathbb{R}$ and $f_1,...,f_p\in \mathcal{C}$. Let $W:=W_{\alpha_0\cdots \alpha_q}$ be as in (\ref{E:setW}) and (\ref{E:PsiU}). By Lemma \ref{L:smallobs} 
\begin{align}
&\Psi_{B_{\alpha_0\cdots \alpha_q}}F(x_0,...,x_{p-1})\notag\\
&=\frac{1}{\mu(W)}\int_W \Alt_{p+1}(f_0\otimes ... \otimes f_p)(t,x_0,...,x_{p-1})\mu(dt)\notag\\
&=\frac{1}{p+1}\sum_{k=0}^p(-1)^k\left(\frac{1}{\mu(W)}\int_Wf_k(t)\mu(dt)\right)\Alt_p(f_0\otimes...\otimes\hat{f}_k\otimes...\otimes f_p)(x_0,...,x_{p-1}),\notag
\end{align}
and consequently $\Psi_{B_{\alpha_0\cdots \alpha_q}}F \in \mathcal{C}^{p-1}(B_{\alpha_0\cdots \alpha_q}^{(p)})$. Identity (\ref{E:homotopy}) is straightforward: Given $F\in \mathcal{C}^p(B_{\alpha_0\cdots \alpha_q}^{(p+1)})$, 
\begin{align}
(&\Psi_{B_{\alpha_0\cdots \alpha_q}}\circ \delta_p+\delta_{p-1}\circ \Psi_{B_{\alpha_0\cdots \alpha_q}})F(x_0,..,x_p)\notag\\
&=\frac{1}{\mu(W)}\int_W\delta_pF(t,x_0,...,x_p)\mu(dt) + \delta_{p-1}\Psi_{B_{\alpha_0\cdots \alpha_q}}F(x_0,...,x_p)\notag\\
&=F(x_0,...,x_p)+\sum_{k=0}^p (-1)^{k+1}\frac{1}{\mu(W)}\int_W F(t,x_0,...,\hat{x}_k,...,x_p)\mu(dt)\notag\\
&\hspace{50pt}+\sum_{k=0}^p(-1)^k\Psi_{B_{\alpha_0\cdots \alpha_q}}F(x_0,...,\hat{x}_k,...,x_p)\notag\\
&=F(x_0,...,x_p).\notag
\end{align}
This completes the proof of (i).

To see the first part of (ii), let $F\in L^2_a(B_{\alpha_0\cdots \alpha_q}^{(p+1)}\setminus D_p,J_p)$. By Jensen's inequality we have 
\[(\Psi_{B_{\alpha_0\cdots \alpha_q}}F (x_0,...,x_{p_1}))^2\leq \frac{1}{\mu(W)}\int_W F(t,x_0,...,x_{p-1})^2\mu(dt).\]
Using (\ref{E:Jpac}), (\ref{E:productiscontained}) and (\ref{E:jlowerbound}), 
\begin{align}
&\int_{B_{\alpha_0\cdots \alpha_q}^{(p)}\setminus D_{p-1}}(\Psi_{B_{\alpha_0\cdots \alpha_q}}F(x_0,...,x_{p-1}))^2 J_{p-1}(d(x_0,...,x_{p-1}))\notag\\
&\leq \frac{1}{p\:\mu(W)}\sum_{k=0}^{p-1}\int_{B_{\alpha_0\cdots \alpha_q}^{(p)}\setminus D_{p-1}}\int_WF(t,x_0,...,x_{p-1})^2\mu(dt)\prod_{\ell\neq k}j(x_k,x_\ell)\mu(dx_0)\cdots \mu(dx_{p-1})\notag\\
&\leq \frac{1}{c^{p-1}\:p\:\mu(W)}\sum_{k=-1}^{p-1}\int_{B_{\alpha_0\cdots \alpha_q}^{(p+1)}\setminus D_{p}}F(x_{-1},x_0,...,x_{p-1})^2\prod_{\ell\neq k}j(x_k,x_\ell)\mu(dx_{-1})\mu(dx_0)\cdots \mu(dx_{p-1})\notag\\
&=\frac{p+1}{c^{p-1}\:p\:\mu(W)}\int_{B_{\alpha_0\cdots \alpha_q}^{(p+1)}\setminus D_{p}} F(x_0,...,x_p)^2J_p(d(x_0,...,x_p)),\notag
\end{align}
where $c$ is as in (\ref{E:jlowerbound}).
To see the second part of (ii), let $F\in \mathcal{D}(\delta_p,B_{\alpha_0\cdots \alpha_q}^{(p+1)})$ and let $(F_k)_k\subset \mathcal{C}^p(B_{\alpha_0\cdots \alpha_q}^{(p+1)})$ be such that 
\begin{equation}\label{E:approxF}
\lim_{k\to \infty} F_k=F\quad \text{in $L^2(B_{\alpha_0\cdots \alpha_q}^{(p+1)}\setminus D_p,J_p)$}
\end{equation}
and 
\begin{equation}\label{E:approxdF}
\lim_{k\to \infty} \delta_pF_k=\delta_pF\quad \text{in $L^2(B_{\alpha_0\cdots \alpha_q}^{(p+2)}\setminus D_{p+1},J_{p+1})$.}
\end{equation}
By (\ref{E:homotopy}) we have 
\begin{equation}\label{E:homotopy2}
\Psi_{B_{\alpha_0\cdots \alpha_q}}\delta_pF_k+\delta_{p-1}\Psi_{B_{\alpha_0\cdots \alpha_q}}F_k=F_k\quad \text{on $B_{\alpha_0\cdots \alpha_q}^{(p+1)}$ $J_p$-a.e.}
\end{equation}
for all $k$. From (\ref{E:PsiUL2est}) and (\ref{E:approxdF}) it follows that 
\[\lim_{k\to\infty}\Psi_{B_{\alpha_0\cdots \alpha_q}}\delta_pF_k=\Psi_{B_{\alpha_0\cdots \alpha_q}}\delta_pF\quad \text{in $L^2(B_{\alpha_0\cdots \alpha_q}^{(p+1)}\setminus D_p,J_p)$,}\]
and taking into account (\ref{E:approxF}) and (\ref{E:homotopy2}), we see that $(\delta_{p-1}\Psi_{B_{\alpha_0\cdots \alpha_q}}F_k)_k$ is a Cauchy sequence in $L^2(B_{\alpha_0\cdots \alpha_q}^{(p+1)}\setminus D_p,J_p)$. Since by (\ref{E:PsiUL2est}) and (\ref{E:approxF}) we have 
\[\lim_{k\to \infty} \Psi_{B_{\alpha_0\cdots \alpha_q}}F_k=\Psi_{B_{\alpha_0\cdots \alpha_q}}F\] 
in $L^2(B_{\alpha_0\cdots \alpha_q}^{(p+1)}\setminus D_p,J_p)$  
and the operator $(\delta_p,\mathcal{D}(\delta_p,B_{\alpha_0\cdots \alpha_q}^{(p+1)}))$ is closed, the function $\Psi_{B_{\alpha_0\cdots \alpha_q}}F$ is in $\mathcal{D}(\delta_p,B_{\alpha_0\cdots \alpha_q}^{(p+1)})$ and 
\[\lim_{k\to\infty} \delta_{p-1}\Psi_{B_{\alpha_0\cdots \alpha_q}}F_k=\delta_{p-1}\Psi_{B_{\alpha_0\cdots \alpha_q}}F\quad\text{in $L^2(B_{\alpha_0\cdots \alpha_q}^{(p+1)}\setminus D_p,J_p)$.}\]
Taking limits in (\ref{E:homotopy2}) we obtain the desired variant of (\ref{E:homotopy}).
\end{proof}

As a consequence we have the following variant of Poincar\'e's lemma: Functions, closed on a given nonempty intersection, are also exact there. 

\begin{corollary}\label{C:Poincare}
Let $\mu$ be a finite nonnegative Borel measure on $X$ with full support. Let $\varepsilon>0$ and $K\geq 0$ be such that Assumption \ref{A:hypostar} holds, let $\eta>0$ be as there and $\mathcal{U}$ as in (\ref{E:Uviaepseta}). Suppose that $B_{\alpha_0\cdots \alpha_q}$ is a nonempty intersection and $1\leq p\leq K+1$. 
\begin{enumerate}
\item[(i)] For any $F\in \mathcal{C}^p(B_{\alpha_0\cdots \alpha_q}^{(p+1)})$ with $\delta_pF=0$ on $B_{\alpha_0\cdots \alpha_q}^{(p+2)}$ there is some $G\in \mathcal{C}^{p-1}(B_{\alpha_0\cdots \alpha_q}^{(p)})$ such that $\delta_{p-1}G=F$.
\item[(ii)] Suppose in addition that Assumptions \ref{A:C}, \ref{A:CE} and \ref{A:abscont} are satisfied and that (\ref{E:jlowerbound}) holds. Then for any $F\in \mathcal{D}(\delta_p,B_{\alpha_0\cdots \alpha_q}^{(p+1)})$ with $\delta_pF=0$ $J_{p+1}$-a.e. on $B_{\alpha_0\cdots \alpha_q}^{(p+2)}$ there is some $G\in \mathcal{D}(\delta_{p-1},B_{\alpha_0\cdots \alpha_q}^{(p)})$ such that $\delta_{p-1}G=F$ on $B_{\alpha_0\cdots \alpha_q}^{(p+1)}$ $J_p$-a.e.
\end{enumerate}
\end{corollary}

For any integer $q\geq 0$ let $\check{C}^q(\mathcal{U},\mathbb{R})$ be the space of \v{C}ech cochains of order $q$ with real coefficients associated with the ball cover $\mathcal{U}$ as defined in (\ref{E:Uviaepseta}). We make the following additional assumption.

\begin{assumption}\label{A:separatecomps}
For any distinct $\alpha_0,...,\alpha_q\in I$ the intersection $B_{\alpha_0\cdots \alpha_q}$ is connected or consists of 
finitely many connected components that have positive minimal distance.
\end{assumption}

If Assumptions \ref{A:Plateaus} and \ref{A:separatecomps} are satisfied, then for any distinct $\alpha_0,...,\alpha_q\in I$ the locally constant functions on $B_{\alpha_0\cdots \alpha_q}$ are contained in $\mathcal{C}^0(B_{\alpha_0\cdots\alpha_q})$. As a consequence, there are natural inclusions $i$ taking elements of $\check{C}^q(\mathcal{U},\mathbb{R})$  into elements of $K^{q,0}_{\mathcal{C}}(\mathcal{U},\hat{N}_\ast(\varepsilon)):=\prod_{\alpha_0<\dots<\alpha_q} \mathcal{C}^0(B_{\alpha_0\cdots\alpha_q})$. The complex (\ref{E:bicomplex}) can be augmented further by a bottom line,
\begin{equation}\label{E:bicomplex2}
	\begin{tikzcd}
		& \vdots & \vdots & \vdots \\
		0 \ar[r] & {\mathcal{C}^0(\hat{N}_0(\varepsilon))} \ar[u, "\delta_0"]\ar[r, "r"] & {\prod\limits_{\alpha_0}\:\mathcal{C}^0(B_{\alpha_0})} \ar[u, "\delta_0"]\ar[r, "\check{\delta}"] & {\prod\limits_{\alpha_0 < \alpha_1}\mathcal{C}^0(B_{\alpha_0\alpha_1})} \ar[u, "\delta_0"]\ar[r, "\check{\delta}"] & \dotsb \\ 
			 & 0 \ar[r]\ar[u] & {\check{C}^0(\mathcal{U}, \mathbb{R})} \ar[u, "i"]\ar[r, "\check{\delta}"] & {\check{C}^1(\mathcal{U}, \mathbb{R})} \ar[u, "i"]\ar[r, "\check{\delta}"] &  \dotsb \\ 
			 & & 0 \ar[u] & 0 \ar[u] & 
	\end{tikzcd}
\end{equation}
Note that a function $f\in \mathcal{C}^0(B_{\alpha_0\cdots \alpha_q})$ with $\delta_0f=0$ on $B_{\alpha_0\cdots \alpha_q}^{(2)}$ is locally constant on $B_{\alpha_0\cdots \alpha_q}$.

If in addition Assumptions \ref{A:C}, \ref{A:CE} and \ref{A:abscont} are satisfied, then similar inclusions $i$ take elements of $\check{C}^q(\mathcal{U},\mathbb{R})$, interpreted in the sense of $\mu$-equivalence classes, into elements of $K^{q,0}_{\mathcal{D}}(\mathcal{U},\hat{N}_\ast(\varepsilon)):=\prod_{\alpha_0<\dots<\alpha_q} \mathcal{D}(\delta_0,B_{\alpha_0\cdots\alpha_q})$. This interpretation is unambiguous, because each $\mu$-equivalence class of locally constant functions contains exactly one representative locally constant in the strict (everywhere) sense. Any $f\in \mathcal{D}(\delta_0, B_{\alpha_0\cdots \alpha_q})$ with $\delta_0f=0$ $J_1$-a.e. is $\mu$-a.e. locally constant on $B_{\alpha_0\cdots \alpha_q}$. 

We write $H^\ell\check{C}^\ast(\mathcal{U},\mathbb{R})$ for the $\ell$-th \v{C}ech cohomology of the ball cover $\mathcal{U}$.
Corollaries \ref{C:MVlemma} and \ref{C:Poincare}, together with similar abstract arguments as before, \cite[Theorem 8.9]{BT82}, \cite[Corollary 3]{BSSS12}, give the following.

\begin{corollary}\label{C:Poincare2}
Let Assumption \ref{A:Plateaus} be in force. Suppose that $\varepsilon>0$ and $K\geq 0$ are such that Assumption \ref{A:hypostar} holds, $\eta>0$ is as there, $\mathcal{U}$ as in (\ref{E:Uviaepseta}) and $\mu$ as specified above. Suppose that also Assumption \ref{A:separatecomps} holds.
\begin{enumerate}
\item[(i)] For any integer $0\leq \ell\leq K$ the spaces $H^\ell\mathcal{C}^\ast(\hat{N}_\ast(\varepsilon))$, $H^\ell\mathcal{K}_\mathcal{C}^\ast(\mathcal{U},\hat{N}_\ast(\varepsilon))$ and $H^\ell\check{C}^\ast(\mathcal{U},\mathbb{R})$
are isomorphic.
\item[(ii)] Suppose that also Assumptions \ref{A:C}, \ref{A:CE} and \ref{A:abscont} are satisfied and that (\ref{E:jlowerbound}) holds for any nonempty intersection of sets from $\mathcal{U}$.
Then for any integer $0\leq \ell\leq K$ the spaces $H^\ell\mathcal{D}(\delta_\ast, \hat{N}_\ast(\varepsilon))$, $H^\ell\mathcal{K}_\mathcal{D}^\ast(\mathcal{U},\hat{N}_\ast(\varepsilon))$ and $H^\ell\check{C}^\ast(\mathcal{U},\mathbb{R})$ are isomorphic, and they are also isomorphic to $H^\ell\mathcal{C}^\ast(\hat{N}_\ast(\varepsilon))$ and $H^\ell\mathcal{K}_\mathcal{C}^\ast(\mathcal{U},\hat{N}_\ast(\varepsilon))$.
\end{enumerate}
\end{corollary}

\begin{remark}
In the case of local complexes on manifolds one typically uses smoothing methods to show that the cohomologies defined in terms of cores and the cohomologies defined in terms of their closures are isomorphic, see for instance \cite[Theorems 2.12 and 3.5]{BL92} or \cite[Section 8]{Cheeger}. In the non-local case the implementation of similar smoothing arguments does not seem straightforward,  but one can pass directly from the core to its closure as done in Propositions \ref{P:MVlemma} and \ref{P:Poincare} to obtain $H^\ell\mathcal{D}(\delta_\ast, \hat{N}_\ast(\varepsilon))\cong H^\ell\mathcal{C}^\ast(\hat{N}_\ast(\varepsilon))$ as in Corollary \ref{C:Poincare2} (ii).
\end{remark}

\subsection{Recovering deRham cohomology}\label{SS:deRham}

Suppose that $X=M$ is a compact smooth Riemannian manifold of dimension $d$. Let $\mathcal{C}:=\lip(M)$. Then clearly Assumption \ref{A:Plateaus} is satisfied. Let $\mathrm{r}_c(M)$ be the convexity radius of $M$, \cite[Section IX.6]{Chavel}.
Recall that if $d\geq 2$, then $\mathrm{r}_c(M)\geq \min \{\frac{\inj(M)}{2},\frac{\pi}{2\sqrt{k}}\}$, where $\inj(M)>0$ denotes the injectivity radius of $M$ and $k>0$ is an upper bound on its sectional curvatures, \cite[Theorem IX.6.1]{Chavel}. Assume that 
\begin{equation}\label{E:boundforeps2eta}
0<\varepsilon<\min\{\mathrm{r}_c(M), \frac{\pi}{2\sqrt{k}}\}\quad \text{if $d\geq 2$}\quad \text{and}\quad 0<\varepsilon<\mathrm{r}_c(M)\quad \text{if $d=1$.}
\end{equation}
Let $\hat{N}_\ast(\varepsilon)$ be as in (\ref{E:boundedrangehat}). Then a combination of results from \cite{BSSS12}, classical theorems and Corollary \ref{C:Poincare2} shows that the complexes $\mathcal{C}^\ast(\hat{N}_\ast(\varepsilon))$ and $\mathcal{D}(\delta_\ast,\hat{N}_\ast(\varepsilon))$ can be used to recover the deRham cohomology of $M$. By $H^\ell_{dR}\Omega^\ast(M)$ we denote the $\ell$-th deRham cohomology and by $H^\ell\check{C}^\ast(M,\mathbb{R})$ the $\ell$-th \v{C}ech cohomology of $M$. 
To a measure $\mu$ with strictly positive and smooth density (with respect to the Riemannian volume) we refer as a \emph{smooth measure}.

\begin{theorem}\label{T:deRham}
Let $M$ be a compact smooth Riemannian manifold of dimension $d$ and $\mu$ a smooth measure on $M$. Suppose that $\varepsilon$ is as in (\ref{E:boundforeps2eta}) and $N_\ast(\varepsilon)$ as defined in (\ref{E:boundedrangehat}). 
\begin{enumerate}
\item[(i)] For any integer $\ell\geq 0$ the spaces $H^\ell\mathcal{C}^\ast(\hat{N}_\ast(\varepsilon))$, $H^\ell\check{C}^\ast(M,\mathbb{R})$ and $H^\ell_{dR}\Omega^\ast(M)$ are finite dimensional and isomorphic; they are trivial for $\ell>d$.
\item[(ii)] If $j$ is a kernel such that Assumptions \ref{A:C}, \ref{A:CE} and \ref{A:abscont} are satisfied and (\ref{E:jlowerbound}) holds, then for any integer $\ell\geq 0$ also $H^\ell\mathcal{D}(\delta_\ast,\hat{N}_\ast(\varepsilon))$ is isomorphic to the spaces in (i). 
\end{enumerate}
\end{theorem}

\begin{examples}
If $j$ is as in (\ref{E:frackerneld}), then the hypotheses in (ii) are satisfied.
\end{examples}

\begin{proof}
Let $K\geq 0$ be an integer and let $\varepsilon$ be as in (\ref{E:boundforeps2eta}). In the case $d\geq 2$ \cite[Theorem 11 and Propositions 20 and 21]{BSSS12} show that Assumption \ref{A:hypostar} holds with $\eta>0$ satisfying $\varepsilon+2\eta<\min\{\mathrm{r}_c(M), \frac{\pi}{2\sqrt{k}}\}$. The proof of \cite[Proposition 21]{BSSS12} uses a well-known consequence of the Rauch comparison theorem, \cite[Theorem 2.7.6]{Klingenberg}. In the case $d=1$ Assumption \ref{A:hypostar} holds with $\eta>0$ such that $\varepsilon+2\eta<\mathrm{r}_c(M)$ by Lemma \ref{L:onedim} below. 

Let $\mathcal{U}=\{B_\alpha\}_{\alpha\in I}$ be a finite cover of $M$ by open balls $B_\alpha$ of radius $\varepsilon+2\eta$. Since the $B_\alpha$ have radii less or equal to $\mathrm{r}_c(M)$, $\mathcal{U}$ is a good cover, that is, for any $q\geq 0$ and distinct $\alpha_0,...,\alpha_q\in I$ the intersection $B_{\alpha_0\cdots\alpha_q}$ is diffeomorphic to $\mathbb{R}^d$. Consequently for any $\ell\geq 0$ the cohomologies $H^\ell\check{C}^\ast(\mathcal{U},\mathbb{R})$, $H^\ell_{dR}\Omega^\ast(M)$ and $H^\ell\check{C}^\ast(M,\mathbb{R})$ are isomorphic, \cite[Theorem 8.9 and Proposition 10.6]{BT82}. The space $H^\ell_{dR}\Omega^\ast(M)$ is finite dimensional, and it is trivial for $\ell>d$. 

Clearly any good cover satisfies Assumption \ref{A:separatecomps}, so that Corollary \ref{C:Poincare2} now yields statements (i) and (ii).
\end{proof}

\begin{remark}\mbox{}
\begin{enumerate}
\item[(i)] The space $H^\ell\check{C}^\ast(\mathcal{U},\mathbb{R})$ with $\mathcal{U}$ as in the proof, can be added to Theorem \ref{T:deRham} (i).
\item[(ii)] Theorem \ref{T:deRham} result also holds with $\mathcal{C}=C^\infty(M)$.
\end{enumerate}
\end{remark}

\begin{lemma}\label{L:onedim}
Let $M$ be a smooth Riemannian manifold of dimension one. Suppose that $\varepsilon$ is as in (\ref{E:boundforeps2eta}) and $N_\ast(\varepsilon)$ as defined in (\ref{E:boundedrangehat}). Then Assumption \ref{A:hypostar} holds.
\end{lemma} 

\begin{proof}
Let $0<2\eta<\min\{\varepsilon/8,\mathrm{r}_c(M)-\varepsilon\}$. A nonempty finite intersection $B_{\alpha_0\cdots \alpha_q}$ of balls of radius $\varepsilon +2\eta$ is a geodesically convex arc with midpoint $y_{\alpha_0\cdots \alpha_q}$, say.

Suppose first that $q=0$. Given $(x_0,...,x_{p-1})\in B_{\alpha_0}^p\cap \hat{N}_{p-1}(\varepsilon)$, we can find a closed ball $B_{\alpha_0}'$ of radius $\varepsilon$ containing $x_0,...,x_{p-1}$, and shifting it along $M$, we may assume that $B_{\alpha_0}'\subset B_{\alpha_0}$. Then $B(y_{\alpha_0},\varepsilon-4\eta)\subset B_{\alpha_0}'$. This means that for any $t\in 
B(y_{\alpha_0},\varepsilon-4\eta)$ and any $(x_0,...,x_{p-1})\in B_{\alpha_0}^{(p)}$ we can identify a closed ball $B_{\alpha_0}'$ of radius $\varepsilon$ containing $t,x_0,...,x_{p-1}$, in other words,
\[B(y_{\alpha_0},\varepsilon-4\eta)\subset \hat{S}_p(\varepsilon;(x_0,...,x_{p-1}))\quad \text{for all $(x_0,...,x_{p-1})\in B_{\alpha_0}^{(p)}$.} \]
Now suppose that $B_{\alpha_0\cdots \alpha_q}$ is a general nonempty intersection. If its length is less than or equal to $2\varepsilon$, then 
\[B_{\alpha_0\cdots \alpha_q}\subset \hat{S}_p(\varepsilon;(x_0,...,x_{p-1}))\quad \text{for all $(x_0,...,x_{p-1})\in B_{\alpha_0\cdots \alpha_q}^{(p)}$.} \]
If not, then its length is $2\varepsilon+2\varepsilon_{\alpha_0\cdots \alpha_q}$ with some $\varepsilon_{\alpha_0\cdots \alpha_q}>0$; we may always assume that $\varepsilon_{\alpha_0\cdots \alpha_q}<2\eta$ (otherwise the balls are identical). Given $(x_0,...,x_{p-1})\in B_{\alpha_0\cdots \alpha_q}^p\cap \hat{N}_{p-1}(\varepsilon)$, we can again find a closed ball $B_{\alpha_0\cdots \alpha_q}'$ of radius $\varepsilon$ containing $x_0,...,x_{p-1}$, and we may again assume it is a subset of $ B_{\alpha_0\cdots \alpha_q}$. Similarly as before this gives $B(y_{\alpha_0\cdots \alpha_q},\varepsilon-2\varepsilon_{\alpha_0\cdots \alpha_q})\subset B_{\alpha_0\cdots \alpha_q}'$, and therefore 
\[B(y_{\alpha_0\cdots \alpha_q},\varepsilon-2\varepsilon_{\alpha_0\cdots \alpha_q})\subset \hat{S}_p(\varepsilon;(x_0,...,x_{p-1}))\quad \text{for all $(x_0,...,x_{p-1})\in B_{\alpha_0\cdots \alpha_q}^{(p)}$.} \]
\end{proof}

\begin{remark}\label{R:onedim}
The proof shows that under the hypotheses of Lemma \ref{L:onedim},
$H^\ell\mathcal{C}^\ast(B_{\alpha}^{(\ast+1)})=\{0\}$, $\ell\geq 1$, for all $\alpha\in I$. This follows using Corollary \ref{C:Poincare} (i). Under the additional assumptions in Corollary \ref{C:Poincare} (ii) also $H^\ell\mathcal{D}(\delta_\ast, B_{\alpha}^{(\ast+1)})=\{0\}$, $\ell\geq 1$, holds for all $\alpha\in I$.
\end{remark}

\section{Some basic examples}\label{SS:Examples} We provide some very basic examples. 

\begin{examples} Let $(X,\varrho)$ be a metric space and $\mathcal{C}=\lip_b(X)$ the algebra of bounded Lipschitz functions on $X$. Let $N_\ast(\varepsilon)$ be as in Example \ref{Ex:N} (ii). If $\varepsilon >\diam(X)$, then $H^0\mathcal{C}^\ast(N_\ast(\varepsilon))\cong \mathbb{R}$ and $H^\ell\mathcal{C}^\ast(N_\ast(\varepsilon))=\{0\}$ for all $\ell\geq 1$; note that if  $\ell\geq 1$ and $F\in \mathcal{C}^\ell(N_\ell(\varepsilon))$ is such that $\delta_\ell F=0$, then for any $y\in X$ the function 
\[G(x_0,...,x_{\ell-1}):=F(y,x_0,...,x_{\ell-1})\] 
is in $\mathcal{C}^{\ell-1}(N_{\ell-1}(\varepsilon))$ and satisfies $\delta_{\ell-1}G=F$. If $(X,\varrho)$ is locally compact, $\mu$ and $j$ satisfy Assumptions \ref{A:C}, \ref{A:CE} and \ref{A:abscont}, then an averaged version of this argument gives again 
$H^0\mathcal{D}(\delta_\ast, N_\ast(\varepsilon))\cong\mathbb{R}$ and $H^\ell\mathcal{D}(\delta_\ast, N_\ast(\varepsilon))=\{0\}$, $\ell\geq 1$. 
\end{examples}

\begin{examples}
Suppose that the metric space $(X,\varrho)$ has two connected components $X_1$ and $X_2$, $\mathcal{C}=\lip_b(X)$ and that $N_\ast(\varepsilon)$ is as in Example \ref{Ex:N} (ii). If $\varepsilon>\dist(X_1,X_2)$, then $H^0\mathcal{C}^\ast(N_\ast(\varepsilon))\cong \mathbb{R}$; if
$\varepsilon\leq \dist(X_1,X_2)$, then $H^0\mathcal{C}^\ast(N_\ast(\varepsilon))\cong \mathbb{R}^2$. If $(X,\varrho)$ is locally compact, both $X_1$ and $X_2$ have nonempty interior and $\mu$ and $j$ are such that Assumptions \ref{A:C}, \ref{A:CE} and \ref{A:abscont} hold, then a similar observation is true for $H^0\mathcal{D}(\delta_\ast, N_\ast(\varepsilon))$. 
\end{examples}

\begin{examples}\label{Ex:S1}
Let $X$ be the unit circle $S^1:=\{e^{i\theta}: 0\leq \theta<2\pi\}$, let $\mathcal{C}=\lip(S^1)$ be the algebra of Lipschitz functions on $S^1$ and $\mu$ the Riemannian volume (Haar measure) on $S^1$. Assume that $j$ satisfies (\ref{E:frackerneld}) with $d=1$ and some $\alpha\in (0,2)$. An application of Theorem \ref{T:deRham} to $M=S^1$ with $0<\varepsilon<\pi/2$ shows that $H^0\mathcal{C}^\ast(\hat{N}_\ast(\varepsilon))$ and $H^0\mathcal{D}(\delta_\ast,\hat{N}_\ast(\varepsilon))$ are both one-dimensional and also 
\begin{equation}\label{E:H1S1}
\text{ $H^1\mathcal{C}^\ast(\hat{N}_\ast(\varepsilon))$ and $H^1\mathcal{D}(\delta_\ast,\hat{N}_\ast(\varepsilon))$ are both one-dimensional.}
\end{equation}
The results for order zero could also be concluded directly from the fact that $\ker \delta_0=\mathbb{R}$. An alternative way to see see (\ref{E:H1S1}) is to (repeatedly) inspect long exact sequences as in \cite[Example 2.6]{BT82}: 
We can patch single balls $B_\alpha$ together to larger and larger open arcs, the fact that their first cohomologies are trivial (cf. Remark \ref{R:onedim}) propagates from the $B_\alpha$ to the larger arcs by the exactness of the long sequence. Eventually we cover all of $S^1$ by two arcs whose intersection has two connected components, and exactness gives (\ref{E:H1S1}). A generating element for $H^1\mathcal{C}^\ast(\hat{N}_\ast(\varepsilon))$ can be constructed as in \cite[Example 2.6]{BT82} or using a 'spiral with constant slope'.

\end{examples}

\begin{examples}\label{Ex:interval}
Let $X$ be the unit interval $[0,1]$, let $\mathcal{C}=\lip([0,1])$ and let $\mu$ be the one-dimensional Lebesgue measure, restricted to $[0,1]$. Assume that $\varepsilon \ll 1$ and that $j$ satisfies Assumptions \ref{A:C}, \ref{A:CE} and \ref{A:abscont}. Then $H^0\mathcal{C}^\ast(\hat{N}_\ast(\varepsilon))$ and $H^0\mathcal{D}(\delta_\ast,\hat{N}_\ast(\varepsilon))$ are both one-dimensional. We can cover $[0,1]$ by finitely many open balls of radius slightly larger than $\varepsilon$, on each of these balls a closed $1$-function is exact (Remark \ref{R:onedim}), and we can use long exact sequences to see that $H^1\mathcal{C}^\ast(\hat{N}_\ast(\varepsilon))$ and $H^1\mathcal{D}(\delta_\ast,\hat{N}_\ast(\varepsilon))$ are trivial.
\end{examples}

\begin{examples}\label{Ex:pointedinterval}
Let again $X=[0,1]$. Let $\mathring{X}:=[0,1]\setminus \{1/2\}$ and let $\mathring{\mathcal{C}}$ the ideal in $\lip([0,1])$ consisting of functions that vanish in a neighborhood of $1/2$. Let $\mu$ be the one-dimensional Lebesgue measure, restricted to $[0,1]$, $\varepsilon \ll 1$ and $j$ a kernel satisfying (\ref{E:frackerneld}) with $d=1$ and some $\alpha\in (0,2)$. Let $\mathring{D}(\delta_\ast,\hat{N}_\ast(\varepsilon))$ be as explained before (\ref{E:fullcoincidence}). If $\alpha\leq 1$, then $H^0\mathring{\mathcal{D}}(\delta_\ast,\hat{N}_\ast(\varepsilon))$ is one-dimensional and $H^1\mathring{\mathcal{D}}(\delta_\ast,\hat{N}_\ast(\varepsilon))$ is trivial by Example \ref{Ex:interval}, Corollary \ref{C:catalogue}, Remark \ref{R:completeness} and Remark \ref{R:adapted}. If $\alpha>1$, then by Corollary \ref{C:catalogue} and Remark \ref{R:adapted} the complex $\mathring{\mathcal{D}}(\delta_\ast,\hat{N}_\ast(\varepsilon))$ is known to be different from $\mathcal{D}(\delta_\ast,\hat{N}_\ast(\varepsilon))$. In fact, we have $H^0\mathring{\mathcal{D}}(\delta_\ast,\hat{N}_\ast(\varepsilon))=\{0\}$ since $\mathring{\mathcal{D}}(\delta_0,\hat{N}_0(\varepsilon))$ does not contain nonzero constants. We can cover $[0,1]$ by finitely many open balls of radius slightly larger than $\varepsilon$, and proceeding similarly as in Lemma \ref{L:onedim}, we can see that Assumption \ref{A:hypostar} holds. These balls have trivial first cohomology (Remark \ref{R:onedim} resp. Corollary \ref{C:Poincare})). A variant of \cite[Proposition 2.3]{BT82} remains true, and using long exact sequences we find that also $H^1\mathring{\mathcal{D}}(\delta_\ast,\hat{N}_\ast(\varepsilon))=\{0\}$. 
\end{examples}

\begin{examples}
Consider the unit circle $X=S^1$ with $\mathcal{C}=\lip(S^1)$, the Riemannian volume $\mu$ and with $j$ satisfying (\ref{E:frackerneld}) with $d=1$ and some $\alpha\in (0,2)$. Let $x\in S^1$ and $\mathring{S}^1:=S^1\setminus \{x\}$. Let  $0<\varepsilon<\pi/2$. If $\alpha\leq 1$, then $H^0\mathring{\mathcal{D}}(\delta_\ast,\hat{N}_\ast(\varepsilon))\cong \mathbb{R}$ and $H^1\mathring{\mathcal{D}}(\delta_\ast,\hat{N}_\ast(\varepsilon))\cong \mathbb{R}$ by Example \ref{Ex:S1}, Corollary \ref{C:catalogue}, Remark \ref{R:completeness} and Remark \ref{R:adapted}. If $\alpha>1$, then $\mathring{\mathcal{D}}(\delta_\ast,\hat{N}_\ast(\varepsilon))\neq \mathcal{D}(\delta_\ast,\hat{N}_\ast(\varepsilon))$. Again $H^0\mathring{\mathcal{D}}(\delta_\ast,\hat{N}_\ast(\varepsilon))=\{0\}$ since nonzero constants are lost. Similarly as in Example \ref{Ex:pointedinterval} a small enough open arc containing a neighborhood of $x$ has trivial zero and first cohomologies. Patching together covering arcs eventually produces two arcs that cover $\mathring{S}^1$. If $x$ is contained in only one of these two arcs, their overlap has two connected components, both arcs have zero $H^1$, one has trivial $H^0$, the other one-dimensional $H^0$. Using a long exact sequence as in Example \ref{Ex:S1} (resp. \cite[Example 2.6]{BT82}),  $H^1\mathring{\mathcal{D}}(\delta_\ast,\hat{N}_\ast(\varepsilon))$ is seen to be one-dimensional. If $x$ is contained in the overlap, a slightly different calculation gives the same result.
\end{examples}

\begin{examples}
Consider the unit sphere $X=S^2$, let $\mathcal{C}=\lip(S^2)$ and $N_\ast=N_\ast(\varepsilon)$ with small $\varepsilon>0$. Suppose that $\mu$ is the Riemannian volume and that $j$ satisfies (\ref{E:frackerneld}) with $d=2$ and some $\alpha\in (0,2)$.
Then $H^0\mathcal{D}(\delta_\ast,N_\ast)\cong \mathbb{R}$. Now let $\Sigma\subset S^2$ be a closed set of Hausdorff dimension $0\leq \beta< 2$ and consider $\mathring{X}:=S^2\setminus \Sigma$. If $\alpha< 2-\beta$, then $H^0\mathring{D}(\delta_\ast,N_\ast)\cong \mathbb{R}$ by Corollary \ref{C:catalogue}. This remains true if $\alpha=2-\beta$ and $\mathcal{H}^\beta(\Sigma)<+\infty$, for instance if $\beta=1$ and $\Sigma$ is the equator. If $\alpha>2-\beta$, then $H^0\mathring{D}(\delta_\ast,N_\ast)=\{0\}$, because nonzero constants are no longer in $\mathring{D}(\delta_0,N_0)$.
\end{examples}


\begin{thebibliography}
\normalsize

\bibitem{AbatangeloValdinoci14}
N. Abatangelo, E. Valdinoci, \emph{A notion of nonlocal curvature}, Numer. Funct. Anal. Optim.
35 (7--9) (2014), 793--815.

\bibitem{AH96}
D.R. Adams, L.I. Hedberg, \emph{Function Spaces and Potential Theory}, Grundlehren der math. Wiss. 314,
Springer, New York, 1996.

\bibitem{Alexander35}
J.W. Alexander, \emph{On the ring of a compact metric space}, Proc. Acad. Sci. {\bf 21} (1935), 509--512.

\bibitem{Allain75}
G. Allain, \emph{Sur la repr\'esentation des formes de Dirichlet},
Ann. Inst. Fourier {\bf 25} (1975), 1--10.

\bibitem{AGS14}
L. Ambrosio, N. Gigli, G. Savar\'e, \emph{Metric measure spaces with Riemannian Ricci curvature bounded from
below}, Duke Math. J. {\bf 163} (7) (2014), 1405--1490.

\bibitem{BCW}
H. Barcelo, V. Capraro, J.A. White, \emph{Discrete homology theory for metric spaces},
Bull. London Math. Soc. {\bf 46} (5) (2014), 

\bibitem{BarlowGrigKumagai09}
M.T. Barlow, A. Grigoryan, T. Kumagai, \emph{Heat kernel upper bounds for jump processes
and the first exit time}, J. reine angew. Math. {\bf 626} (2009), 135--157.

\bibitem{BSSS12}
L. Bartholdi, T. Schick, N. Smale, S. Smale, \emph{Hodge theory on metric spaces},
Found. Comp. Math. {\bf 12} (2012), 1--48.

\bibitem{Bertoin96}
J. Bertoin, \emph{L\'{e}vy Processes}, Cambridge Tracts in Mathematics, Vol. 121, Cambridge
University Press, Cambridge, 1996.

\bibitem{BD58}
A. Beurling, J. Deny, \emph{Espaces de Dirichlet: I. Le cas \'el\'ementaire}, Acta Math. {\bf 99} (1958), 203--224.

\bibitem{BD59}
A. Beurling, J. Deny, \emph{Dirichlet space}, Proc. Nat. Acad. Sci. U.S.A. {\bf 45} (1959), 208--215.

\bibitem{BlooreRoberts}
F.J. Bloore, G. Roberts, \emph{The structure of $\mathcal{F}$-tensorial cochains of differential operators},
Diff. Geom. and Appl. {\bf 10} (1999), 295--301.

\bibitem{BT82}
R. Bott, L.W. Tu, \emph{Differential Forms in Algebraic Topology}, Grad. Texts in Math. vol. 82, Springer, Berlin, 1982. 

\bibitem{BH91}
N. Bouleau, F. Hirsch, \emph{Dirichlet Forms and Analysis on Wiener Space},
deGruyter Studies in Math. 14, deGruyter, Berlin, 1991.

\bibitem{Bourbaki}
N. Bourbaki, \emph{Elements of Mathematics, General Topology, Part 1 }, Springer, Berlin, 1995.

\bibitem{BL92}
J. Br\"uning, M. Lesch, \emph{Hilbert complexes},
J. Funct. Anal. {\bf 108} (1992), 88-132.

\bibitem{BucurValdinoci16}
C. Bucur, E. Valdinoci, \emph{Nonlocal Diffusion and Applications}, Lecture Notes of the Unione Matematica Italiana, Springer Intl. Publ., Switzerland, 2016.

\bibitem{CaffarelliRoquejoffreSavin10}
L. Caffarelli, J.-M. Roquejoffre, O. Savin, \emph{Nonlocal minimal surfaces}, Comm. Pure Appl. Math.
{\bf 63} (9) (2010), 1111--1144.

\bibitem{CaffarelliSilvestre07}
L. Caffarelli, L. Silvestre, \emph{An extension problem related to the fractional Laplacian}, Comm. Partial Differential Equations {\bf 32} (2007), 1245--1260.


\bibitem{Carlsson}
G. Carlsson, \emph{Topology and Data}, Bull. Amer. Math. Soc. {\bf 46} (2) (2009), 255--308.

\bibitem{ChangGonzales11}
S.-Y.A. Chang, M. del Mar Gonz\'{a}les, \emph{Fractional Laplacian in conformal geometry}, Adv. Math. {\bf 226} (2) (2011), 1410--1432.

\bibitem{Chavel}
I. Chavel, \emph{Riemannian Geometry. A Modern Introduction}, Cambridge Univ. Press, Cambridge, 2006.

\bibitem{ChazaldeSilvaOudot}
F. Chazal, V. deSilva, S. Oudot, \emph{Persistence stability for geometric complexes},
Geom. Dedicata {\bf 173} (2014), 193--214.

\bibitem{Cheeger}
J. Cheeger, \emph{On the Hodge theory of Riemannian pseudomanifolds}, In: Geometry of the
Laplace operator, Proc. Sympos. Pure Math. Vol. 36, Amer. Math. Soc., Providence, 1980, pp. 91--146.

\bibitem{ChKu03}
Z.-Q. Chen, T. Kumagai, \emph{Heat kernel estimates for stable like processes on $d$-sets},
Stoch. Proc. Appl. {\bf 108} (2003), 27-62.

\bibitem{ChKu08}
Z.-Q. Chen, T. Kumagai, \emph{Heat kernel estimates for jump processes of mixed types
on metric measure spaces}, Probab. Theory Relat. Fields {\bf 140} (2008), 277--317.

\bibitem{CS03}
F. Cipriani, J.-L. Sauvageot, \emph{Derivations as square roots of Dirichlet forms},
J. Funct. Anal. {\bf201} (2003), 78--120.

\bibitem{CM90}
A. Connes, H. Moscovici, \emph{Cyclic cohomology, the Novikov conjecture and hyperbolic groups},
Topology {\bf 29} (3) (1990), 345--388.

\bibitem{DaLioRiviere11}
F. Da Lio, T. Rivi\`{e}re, \emph{Three-term commutator estimates and the regularity of $\frac12$-harmonic maps
into spheres}, Anal. PDE {\bf 4} (1) (2011), 149--190.

\bibitem{dRh60}
G. deRham, \emph{Vari\'et\'es Differentiables}, Hermann, 1960, Paris.

\bibitem{dRh-K50}
G. deRham, K. Kodaira, \emph{Harmonic Integrals}, Inst. Adv. Study, Princeton, 1950.

\bibitem{DiCastroKuusiPalatucci14}
A. Di Castro, T. Kuusi, G. Palatucci, \emph{Nonlocal Harnack inequalities}, J. Funct. Anal. {\bf 267} (6) (2014),
1807--1836.

\bibitem{diNezza12}
E. Di Nezza, G. Palatucci and E. Valdinoci,\emph{Hitchhiker’s guide to the fractional Sobolev spaces},
Bull. Sci. Math. {\bf 136} (5) (2012), 521--573.

\bibitem{DipierroRos-OtonValdinoci17}
S. Dipierro, X. Ros-Oton, E. Valdinoci, \emph{Nonlocal problems with Neumann boundary conditions},
Rev. Mat. Iberoam. {\bf 33} (2) (2017) 377--416.

\bibitem{Eckmann44}
B. Eckmann, \emph{Harmonische Funktionen und Randwertaufgaben in einem Komplex}, Comment. Math.
Helv. {\bf 17} (1945), 240--255.

\bibitem{ELZ}
H. Edelsbrunner, D. Letscher, A. Zomorodian, \emph{Topological persistence and simplification}, Discr.
Comput. Geom. {\bf 28} (4) (2002), 511--533.

\bibitem{EKS15}
M. Erbar, K. Kuwada, K-Th. Sturm, \emph{On the equivalence of the entropic curvature-dimension condition
and Bochner's inequality on metric measure spaces}, Invent. Math. {\bf 201} (2015), 993--1071.

\bibitem{Federer}
H. Federer, \emph{Geometric Measure Theory}, Grundlehren math. Wiss. 153, Springer, Berlin, 1996.

\bibitem{FOT94}
M. Fukushima, Y. Oshima and M. Takeda, \emph{Dirichlet forms and symmetric Markov processes},
deGruyter, Berlin, New York, 1994.

\bibitem{Gaffney55}
M.P. Gaffney, \emph{Hilbert space methods in the theory of harmonic integrals},
Trans. Amer. Math. Soc. {\bf 78} (2) (1955), 426--444.

\bibitem{Garofalo19}
N. Garofalo, \emph{Fractional thoughts}, In: New Developments in the Analysis of Nonlocal Operators, D. Danielli,  A. Petrosyan, C.A. Pop, editors, Contemp. Math. {\bf 273} (2019), pp. 1--136.

\bibitem{Genton}
L. Genton, \emph{Scaled Alexander-Spanier Cohomology and Lqp Cohomology
for Metric Spaces}, TH\`ESE No. 6330 (2014), PhD-thesis, EPFL, Lausanne, 2014.

\bibitem{Gigli17}
N. Gigli, \emph{Non-smooth differential geometry - an approach tailored for spaces with Ricci curvature bounded from below}, Mem. Amer. Math. Soc. {\bf 251} (11) (2017). 

\bibitem{GT10}
V. Gol’dshtein, M. Troyanov, \emph{A conformal deRham complex}, J. Geom. Anal. {\bf 20} (2010), 651--669.

\bibitem{GVF}
J. Gracia-Bondia, J. V\'arilly, H. Figueroa, \emph{Elements of Noncommutative Geometry}, Birkh\"auser, 2001.

\bibitem{Grigoryan2009}
A. Grigor'yan, \emph{Heat kernel and analysis on manifolds}, AMS/IP Studies in Advanced Mathematics 47, Boston, MA, 2009.

\bibitem{GrigHuHu18}
A. Grigoryan, E. Hu, J. Hu, \emph{Two-sided estimates of heat kernels of jump type
Dirichlet forms}, Adv. Math. {\bf 330} (2018), 433--515.

\bibitem{GrigHuLau14}
A. Grigoryan, J. Hu, K.-S. Lau, \emph{Estimates of heat kernels for non-local regular Dirichlet forms},
Trans. Amer. Math. Soc. {\bf 366} (12) (2014), 6397--6441.


\bibitem{Gromov87}
M. Gromov, \emph{Hyperbolic groups}, In: Essays in Group Theory, Math. Sci. Res. Inst. Pub. vol. 8, Springer, New York, 1987, pp. 75--263.

\bibitem{Grubb15}
G. Grubb, \emph{Fractional Laplacians on domains, a development of H\"ormander’s theory
of $\mu$-transmission pseudodifferential operators}, Adv. Math. {\bf 268} (2015), 478--528.

\bibitem{Hausmann95}
J.-C. Hausmann, \emph{On the Vietoris-Rips complexes and a cohomology theory for metric
spaces}, Prospects in Topology: Proceedings of a Conference in honour of William
Browder, Princeton, 1995, pp. 175--188.

\bibitem{Heinonen}
J. Heinonen, \emph{Lecture on Analysis on Metric Spaces},
Universitext, Springer, New York, 2001.

\bibitem{H15}
M. Hinz, \emph{Magnetic energies and Feynman-Kac-It\^o formulas for symmetric Markov processes}, 
Stoch. Anal. Appl. {\bf 33}(6) (2015), 1020-1049.


\bibitem{HinzKommer22+}
M. Hinz, J. Kommer, \emph{Non-local approximations to local differential complexes},
in preparation.

\bibitem{HMS22}
M. Hinz, J. Masamune, K. Suzuki, \emph{Removable sets and $L^p$-uniqueness on manifolds and metric measure spaces}, preprint (2022), arXiv:2204.01378.


\bibitem{HRT13}
M. Hinz, M. R\"ockner, A. Teplyaev, \emph{Vector analysis for  Dirichlet forms and quasilinear PDE and SPDE on metric measure spaces}, Stoch. Proc. Appl. {\bf 123} (12) (2013), 4373--4406.

\bibitem{HTams}
M. Hinz,   A. Teplyaev, \emph{Local Dirichlet forms, Hodge theory, and
the Navier-Stokes equations on topologically one-dimensional
fractals}, Trans. Amer. Math. Soc. {\bf 367}  (2015),
1347--1380, Corrigendum in Trans. Amer. Math. Soc. {\bf 369} (2017), 6777--6778.


\bibitem{HinzTeplyaev18}
M. Hinz,  A. Teplyaev, \emph{Densely defined non-closable curl on carpet like metric measure spaces}, Math. Nachr. {\bf 291} (11-12) (2018), 1743--1756.

\bibitem{Hochschild45}
G. Hochschild, \emph{On the cohomology groups of an associative algebra}, Ann. Math. {\bf 46} (1) (1945), 58--67.

\bibitem{Hodge41}
W.V.D. Hodge, \emph{The Theory and Applications of Harmonic Integrals}, Cambridge Univ.
Press, Cambridge, 1941.

\bibitem{Hutchinson}
J.E. Hutchinson, \emph{Fractals and self-similarity}, Indiana Univ. Math. J. {\bf 30} (1980), 713--747.

\bibitem{IRT12}
M. Ionescu, L. Rogers, A. Teplyaev, \emph{Derivations, Dirichlet forms and spectral analysis},
J. Funct. Anal. {\bf 263} (2012), no. 8, 2141--2169. 

\bibitem{Jacob01}
N. Jacob, \emph{Pseudo-Differential Operators and Markov processes. Volume I: Fourier Analysis and Semigroups}, Imperial College Press, London, 2001.

\bibitem{KassmannMimica17}
M. Kassmann, A. Mimica, \emph{Intrinsic scaling properties for nonlocal operators}, J. Eur. Math. Soc. (JEMS), {\bf 19} (4) (2017), 983--1011.

\bibitem{Klingenberg}
W. Klingenberg, \emph{Riemannian Geometry}, deGruyter, Berlin, 1995.

\bibitem{Kodaira49}
K. Kodaira, \emph{Harmonic fields in Riemannian manifolds (generalized potential theory)}, Ann. Math. {\bf 50} (1949), 587--665.

\bibitem{Kolmogorov36}
A.N. Kolmogorov, \emph{\"Uber die Dualitat im Aufbau der
kombinatorischen Topologie}, Rec. Math. (Mat. Sbornik) N.S. {\bf 1} (43) (1) (1936), 97--102.

\bibitem{Kolmogorov36a}
A.N. Kolmogorov, \emph{Homologiering des Komplexes und des lokal-bikompakten Raumes}, 
Rec. Math. (Mat. Sbornik) N.S. {\bf 1} (43) (5) (1936), 701--706.

\bibitem{Latschev01}
J. Latschev, \emph{Vietoris-Rips complexes of metric spaces near a closed Riemannian manifold},
Arch. Math. {\bf 77} (2001), 522--528.


\bibitem{LeJan78}
Y. LeJan, \emph{Mesures associ\'ees \`a une forme de Dirichlet. Applications.}, Bull. Soc. Math. France {\bf 106}
(1978), 61--112.

\bibitem{Lueck}
W. L\"uck, \emph{$L^2$-Invariants: Theory and Applications to Geometry and K-Theory}, Ergebnisse der Mathematik und ihrer Grenzgebiete 44, Springer, Berlin, 2002.

\bibitem{MasamuneUemura2011}
J. Masamune, T. Uemura, \emph{Conservation property of symmetric jump processes}, Ann. Inst. H. Poincar\'e Probab. Stat.
{\bf 47} (3) (2011), 650--662.

\bibitem{Massey78}
W.M. Massey, \emph{Homology and Cohomology Theory. An approach based on Alexander-Spanier cochains}, Monographs and textbooks in pure and applied mathematics, vol. 46, Marcel Dekker, Inc., New York, 1978.

\bibitem{Mazya85}
V. Maz'ya, \emph{Sobolev Spaces}, Springer, New York, 1985.

\bibitem{Nazarov22}
A.I. Nazarov, \emph{Variety of fractional Laplacians}, preprint (2021), arXiv:2108.12924, to appear in Proc. ICM2022.

\bibitem{PalatucciSavinValdinoci13}
G. Palatucci, O. Savin, E. Valdinoci, \emph{Local and global minimizers for a variational energy involving a fractional norm}, Ann. Mat. Pura Appl. {\bf 192} (4) (2013), 673--718.

\bibitem{PalatucciKuusi18}
G. Palatucci, T. Kuusi (editors), \emph{Recent Developments in Nonlocal Theory}, deGruyter Open Poland, 2017.

\bibitem{Pansu96}
P. Pansu, \emph{Introduction to $L^2$ Betti numbers}, In: Riemannian geometry, Waterloo, 1993. Fields
Inst. Monogr. 4, Amer. Math. Soc., Providence, 1996, pp. 53--86.

\bibitem{Pansu04}
P. Pansu, \emph{Cohomologie $L^p$: invariance sous quasiisom\'etries}, preprint (2004).

\bibitem{Plaut18}
C. Plaut, \emph{Spectra related to the length spectrum}, Asian J. Math. {\bf 25} (4) (2021), 521--550.

\bibitem{PlautWilkins12}
C. Plaut, J. Wilkins, \emph{Discrete homotopies and the fundamental group}, Adv. Math. {\bf 232} (2013), 271--294.

\bibitem{RS80}
M. Reed, B. Simon, \emph{Methods of Modern Mathematical Physics I: Functional Analysis}, Academic Press, San Diego, 1980.

\bibitem{Ros-OtonSerra12}
X. Ros-Oton, J. Serra, \emph{The Dirichlet problem for the fractional Laplacian: regularity up to the boundary}, J. Math. Pures Appl. {\bf 101} (2012), 275--302.


\bibitem{Sauvageot89}
J.-L. Sauvageot, \emph{Tangent bimodule and locality for dissipative operators on $C^\ast$-algebras},
Quantum Probability and Applications IV, Lect. Notes Math. 1396, Springer, Berlin, 1989, pp. 322--338.

\bibitem{Sauvageot90}
J.-L. Sauvageot, \emph{Quantum differential forms, differential calculus and semi-groups}, Quantum
Probability and Applications V, Lect. Notes Math. 1442, Springer, Berlin, 1990,
pp. 334--346.

\bibitem{Schikorra15}
A. Schikorra, \emph{Integro-Differential Harmonic Maps into Spheres}, Comm. Partial Diff. Equations {\bf 40} (3) (2015),  
506--539.

\bibitem{Silvestre05}
L. Silvestre, \emph{H\"older estimates for solutions of integro-differential equations like the fractional
Laplace}, Indiana Univ. Math. J. {\bf 55} (3) (2006), 1155--1174.

\bibitem{SS12}
N. Smale, S. Smale, \emph{Abstract and classical Hodge DeRham theory},
Anal. Appl. {\bf 10}(1) (2012), 91--111.

\bibitem{Spanier48}
E.H. Spanier, \emph{Cohomology for general spaces}, Ann. Math. {\bf 49} (2) (1948), 407--427.

\bibitem{Spanier66}
E. H. Spanier, \emph{Algebraic Topology}, Springer, New York, 1966.

\bibitem{SpenerWeberZacher}
A. Spener, F. Weber, R. Zacher, \emph{The fractional Laplacian has infinite dimension}, Comm. Partial Diff.
Eq. {\bf 45} (1) (2020), 57--75.

\bibitem{Stos00}
A. St\'os, \emph{Symmetric $\alpha$-stable processes on $d$-sets}, Bull. Pol. Acad. Sci. Math. {\bf 48} (2000), 237--245.

\bibitem{Vietoris27}
L. Vietoris, \emph{\"Uber den h\"oheren Zusammenhang kompakter R\"aume und eine Klasse von zusammenh\"angenden Abbildungen}, Math. Ann. {\bf 91} (1927), 454--472.

\bibitem{Warner}
F.W. Warner, \emph{Foundations of Differentiable manifolds and Lie Groups}, Grad. Texts in Math. 94, Springer, New York, 1983.

\bibitem{W00}
N. Weaver, \emph{Lipschitz algebras and derivations II. Exterior differentiation}, J. Funct. Anal. {\bf 178} (2000), 64-112.

\bibitem{Weidmann}
J. Weidmann, \emph{Linear Operators in Hilbert Spaces}, Grad. Texts in Math. 68, Springer, New York, 1980.

\bibitem{Weil52}
A. Weil, \emph{Sur les th\'eor\`emes de deRham}, Comment. Math. Helv. {\bf 26} (1952), 119--145.
\end{thebibliography}
\end{document}